\documentclass[3p,12pt]{elsarticle}
\usepackage{lineno,hyperref}
\usepackage{amsmath,bm,geometry}
\usepackage{amsfonts}
\usepackage{color}
\usepackage{nomencl,etoolbox,ragged2e,siunitx}
\usepackage{mathtools}
\usepackage{multirow}
\usepackage{caption}
\usepackage{parskip}
\usepackage{calligra}
\usepackage{makeidx}
\usepackage{tikz}
\usetikzlibrary{matrix}
\makeindex
\DeclareMathAlphabet{\mathcalligra}{T1}{calligra}{m}{n}
\newcommand*{\defeq}{\stackrel{\text{def}}{=}}
\begin{document}
\title{Unified Gas-kinetic Wave-Particle methods I: Continuum and Rarefied Gas Flow }
\author[ad1]{Chang Liu}
\ead{cliuaa@connect.ust.hk}
\author[ad3]{Yajun Zhu}
\ead{zhuyajun@mail.nwpu.edu.cn}
\author[ad1,ad2]{Kun Xu\corref{cor1}}
\ead{makxu@ust.hk}
\address[ad1]{Department of Mathematics, Hong Kong University of Science and Technology, Hong Kong}
\address[ad2]{Department of Mechanical and Aerospace Engineering, Hong Kong University of Science and Technology, Hong Kong, China}
\address[ad3]{National Key Laboratory of Science and Technology on Aerodynamic Design and Research, Northwestern Polytechnical University, Xi'an, Shaanxi 710072, China}
\cortext[cor1]{Corresponding author}

\begin{abstract}
The unified gas-kinetic scheme (UGKS) provides a framework for simulating multiscale transport with
the updates of both gas distribution function and macroscopic flow variables on the
cell size and time step scales.
The multiscale dynamics in UGKS is achieved through the coupled particle transport and collision in the particle
evolution process within a time step.
In this paper, under the UGKS framework, we propose an efficient multiscale unified gas-kinetic wave-particle (UGKWP) method.
The gas dynamics in UGKWP method is described by the individual particle movement coupled with the evolution of the probability density function (PDF).
During a time step, the trajectories of simulation particles are tracked until collision happens,
and the post-collision particles are evolved collectively through the evolution of the corresponding distribution function.
The evolution of simulation particles and distribution function is guided by evolution of macroscopic variables.
The two descriptions on a gas particle, i.e. wave and particle, switch dynamically with time.
A new concept of multiscale multi-efficiency preserving (MMP) method is introduced, and the UGKWP method is shown to be an MMP scheme.
Multiscale preserving means UGKWP method preserves the flow regime from collisionless regime to hydrodynamic regime without requiring the cell size and time step to be less than the mean free path and collision time.
Multi-efficiency preserving means the computational cost of UGKWP method including the computational time and memory cost is on the same level as the particle methods in the rarefied regime, and becomes comparable to the hydrodynamic solvers in continuum regime.
The UGKWP method is specially efficient for hypersonic flow simulation in all regimes in comparison with the wave-type discrete ordinate methods,
and presents a much lower stochastic noise in the continuum flow regime in comparison with the particle-based Monte Carlo methods.
Numerical tests for flows over a wide range of Mach and Knudsen numbers are presented. The examples include mainly the hypersonic flow passing a circular cylinder at Mach numbers $20$ and $30$ and Knudsen numbers $1$ and $10^{-4}$, low speed lid-driven cavity flow, and laminar boundary layer.
These results validate the accuracy, efficiency, and multiscale property of UGKWP method.
\end{abstract}

\maketitle

\section{Introduction}
The Boltzmann equation is a fundamental equation for gas dynamics, which resolves the particle mean free path and collision time scale gas flow physics.
Theoretically, from the Boltzmann equation the flow physics in all Knudsen regimes can be captured through the accumulation of
molecular dynamic evolution on the kinetic scale.
The numerical methods for solving kinetic Boltzmann equations can be categorized into two groups: the stochastic methods and the deterministic methods.
For stochastic methods, the evolution of velocity distribution function is represented by the motion of simulation particles.
Such kind of Lagrangian-type schemes achieve high computational efficiency in rarefied and hypersonic flow simulation.
The real gas effects, such as quantum effect, ionization, and chemical reaction, can be properly handled on the particle level.
Moreover, the particle methods are robust and not sensitive to mesh quality.
However, the particle methods suffer from statistical noise that greatly reduces the computational efficiency for low speed flow simulation.
At the same time, the cell size and time step of conventional particle methods are usually restricted to be less than the mean free path and collision time
due to the splitting treatment of particle transport and collision.
Therefore, in the near continuum regime, with the increase of collision rate, the computational cost will become very high.
The deterministic methods are constructed on a discretized space and time.
Compared to the stochastic method, the deterministic method usually achieves high accuracy.
However, it suffers the ray effect for a relative long time simulation of rarefied flow,
and the computational cost is high for the calculation of hypersonic and multidimensional flow due to the direct discretization of particle velocity space.
Similar to the stochastic methods, many deterministic methods also have the constraints on the cell size and time step for an accurate description of flow physics.

The direct simulation Monte Carlo (DSMC) method is one of the most popular methods for the simulation of rarefied hypersonic flow  \cite{bird1963approach},
which effectively captures the real molecular physics through the separate modeling of particle transport and collision \cite{bird1994molecular}.
Suffering from stochastic noise, DSMC has low efficiency in low speed flow simulation like all particle methods.
Progresses has been made to reduce the noise, improve the efficiency, and extend its applicable regimes.
To deal with the low-signal flow, the information preserving (IP) method \cite{shen2001information,sun2002direct},
as well as an efficient low variance DSMC method \cite{baker2005variance,homolle2007low} has also been developed.
Since the cell size and time step are restricted to be less than the mean free path and collision time, DSMC is highly expensive in the simulation of continuum flow.
In order to extend DSMC to the continuum flow simulation, the asymptotic preserving Monte Carlo methods \cite{pareschi2000asymptotic,ren2014asymptotic},
moment-guided Monte Carlo method \cite{degond2011moment}, low diffusion particle method \cite{burt2008low}, as well as hybrid methods between DSMC and CFD methods \cite{schwartzentruber2006hybrid} have been constructed.
The stochastic methods based on the kinetic model equations, for example the BGK/ES-BGK equation \cite{macrossan2001particle,fei2018particle} and the kinetic Fokker-Planck equation \cite{jenny2010solution} have been proposed to reduce the computational cost of DSMC.
Among those stochastic kinetic methods, the particle Fokker-Planck method \cite{jenny2010solution} and the stochastic BGK method  \cite{fei2018particle} are applicable over different flow regimes.

The deterministic discrete ordinate methods (DOM) for Boltzmann and kinetic equations have been extensively studied in the last several decades \cite{chu1965,JCHuang1995,Mieussens2000,tcheremissine2005direct,Kolobov2007,LiZhiHui2009,ugks2010,wu2015fast,aristov2012direct},
which have great advantages for the simulation of low speed microflow \cite{huang2013,wu2014solving}.
In order to improve the efficiency and remove time step limitations of kinetic methods, asymptotic preserving schemes \cite{jin1999efficient}, implicit schemes \cite{chen2017unified}, and kinetic-fluid hybrid methods \cite{degond2010multiscale} have been proposed and developed.
Following the direct modeling methodology \cite{xu-book}, an effective multiscale unified gas-kinetic scheme (UGKS) has been
proposed and developed \cite{ugks2010,huang2012,liu2016}.
In the construction of UGKS, the particle transport and collision are closely coupled for both flux transport and inner cell evolution,
and the scheme is effectively applicable from rarefied to continuum flows.
In the continuum flow regime, the UGKS is able to accurately capture the Navier-Stokes (NS) solutions, such as capturing the laminar boundary layer,
even with the cell size being much larger than the local particle mean free path and the time step being larger than the local collision time.
In order to improve the computational efficiency, the implicit and multigrid techniques have been incorporated into the UGKS \cite{zhu2016implicit,zhu2017implicit,zhu2018implicit}, and the scheme becomes highly efficient and accurate for flow simulations over a wide range of Knudsen and Mach numbers.
For example, for a space vehicle computation at hypersonic speed in the near space flight, the local Knudsen number around the
flying vehicle can be changed greatly over a fifth-order magnitude and the UGKS can present accurate solution with a large variation of the
ratio between the local cell size and local particle mean free path \cite{jiang}.
The UGKS has also been successfully extended to radiative transfer \cite{sun2015asymptotic,sun2017,sun2018,li2018ugkp}, plasma transport \cite{liu2017}, and multiphase flow \cite{liu2018unified}.

The UGKS framework plays an important role in the construction of unified gas-kinetic wave-particle (UGKWP) method.
Therefore a brief review of UGKS is given here, and readers can refer to \cite{wang2015unified} for detailed formulae.
The UGKS is a finite volume scheme for the update of both gas distribution function $f(\vec{x},t,\vec{v})$ and macroscopic flow variables $\vec{W}$ in physical and velocity space.
For UGKS, the time step is only limited by the CFL condition,
and the cell size and time step are not restricted to be less than the mean free path and collision time to accurately capture the flow physics.
The modeled dynamics depends on the ratios of the time step over the particle collision time and the cell size over the particle mean free path.
For a physical control volume $\Omega_i$ and velocity control volume $\Omega_j$,
the evolution equation of cell averaged gas distribution function $f_{ij}$ in the control volume $\Omega_{ij}=\Omega_i\bigcap\Omega_j$
from time step $t^n$ to $t^{n+1}$ is modeled as
\begin{equation} \label{kinetic-f} \footnotesize
  f_{ij}^{n+1}=
  f_{ij}^n-\frac{1}{|\Omega_{ij}|}\int_{t^n}^{t^{n+1}}\oint_{\partial \Omega_i}
  \vec{v_j}\cdot \vec{n} f_{\partial \Omega_{ij}}(t,\vec{v}_j) ds dt
  +\frac{\Delta t}{2}   (Q^{n}_{ij} + Q^{n+1}_{ij}),
\end{equation}
where $\vec{v}_j$ is the particle velocity,
$f_{\partial \Omega_{ij}}(t,\vec{v}_j)$ is time-dependent solution at a cell interface $\partial\Omega_i$ for the flux evaluation,
and $Q_{ij}$ is cell averaged collision term.
The above evolution equation for the velocity distribution function is coupled with the evolution of cell averaged macroscopic flow variables $\vec{W}_{i}$,
\begin{equation}\label{kinetic-m}
  \vec{W}_i^{n+1}=\vec{W}_i^{n}-\frac{1}{|\Omega_i|}\int^{t^{n+1}}_{t^n}\int\oint_{\partial \Omega_i}
  \psi \vec{v}\cdot \vec{n} f_{\partial \Omega_i}(t,\vec{v}) ds d\Xi dt,
\end{equation}
where $d\Xi = dudvdwd\vec{\xi}$, with $\vec{v}=(u,v,w)$ the particle velocity and $\vec{\xi}$ the internal variable.
In order to close the above two discretized governing equations and capture multiscale gas evolution,
the time-dependent interface gas distribution function $f_{\partial \Omega_i}$
and the cell averaged collision term $Q_{ij}$ need to be modeled.
Since the time step $\Delta t=t^{n+1} - t^n$ can be large in comparison with the particle collision time $\tau$, the particle transport and collision effect have
to be included in the modeling of the interface evolution of the distribution function.
One of the important ingredients in UGKS is to use a time-dependent interface distribution function,
i.e. the integral equation of the kinetic Shakhov model \cite{shakhov1968generalization}, which reads
\begin{equation}\label{integral-solution-ugks}
f_{\partial \Omega_i}(t,\vec{v})=\frac{1}{\tau}\int_{t^n}^{t} f^+(\vec{x}^\prime,t^\prime,\vec{v})\mathrm{e}^{-(t-t^\prime)/\tau}dt^\prime+
\mathrm{e}^{-t/\tau}f_0(\vec{x}_{\partial \Omega_i}-\vec{v}t,\vec{v}),
\end{equation}
where $\tau$ is the local relaxation parameter, $f^+$ is post collision distribution function, $\vec{x}'=\vec{x}_{\partial \Omega_i}-\vec{v}(t-t')$ is the particle trajectory,
and $f_0$ is the distribution function at time $t^n$.
Above integral equation explicitly provides the solution of velocity distribution at cell interface
once the initial distribution is given and the post collision distribution is constructed.
And this time evolution solution couples the particle free transport and collision.
When the time step $\Delta t$ is used in the time evolution solution,
it determines the flow dynamics from the initial distribution to equilibrium evolution, and provides a multiscale flux construction for UGKS.
Based on above UGKS formulation, in the continuum flow regime with $\Delta t \gg \tau$, the integral solution converges to the
Chapman-Enskog type Navier-Stokes distribution function, and the UGKS will automatically give the Navier-Stokes solution, which is the same
solution obtained from the gas-kinetic scheme (GKS) \cite{xu2001gas}.
The GKS is designed for solving NS equations where the Chapman-Enskog expansion is directly used to reconstruct the initial cell interface distribution function $f_0$
from macroscopic flow variables.

One important concept in theoretical and computational fluid dynamics is the multiscale modeling.
In the field of theoretical fluid dynamics, multiscale equation refers to the equation which can recover multiple physical scale flow phenomena.
The concept is trivial because the equations constructed on a specific physical scale always hold for the scales above,
and the multiscale research in theoretical fluid dynamics focuses on deriving concise equations on large scale.
In the field of computational fluid dynamics,
multiscale scheme or multiple physical-numerical scale scheme mean the numerical scheme can capture multiple physical flow phenomena on variable numerical scales.
Different from the theoretical fluid dynamics, the numerical scheme constructed on a small scale may not be valid on large scale.
For example, the DOM-based direct Boltzmann solver is constructed on the scale of mean free path, and its numerical valid scale is fixed to be the mean free path scale.
Even the direct Boltzmann solver is used to simulate the continuum flow, the cell size is required to be less than the local mean free path \cite{chen2015comparative}.
Due to the coupled treatment of particle transport and collision,
the UGKS is a multiscale scheme captures the flow physics from rarefied regime to continuum regime without the small cell size limitation.

Analogous to quantum mechanics, the gas particle can be described in terms of not only particles, but also waves or its probability density function.
At current stage, most numerical schemes use either particle description or wave description, for example the particle description based Monte Carlo methods and the wave description based discrete ordinate methods.
The particle methods suffer from statistical noise and the discrete ordinate methods suffer large computational cost especially for the simulation of hypersonic flow.
A wave-particle formulation is proposed in this paper, and the framework can be applied to simulate many other transport process.

The purpose of this paper is to apply the wave-particle formulation to solve the equations of UGKS.
In unified gas-kinetic wave-particle method, the gas particles are divided into the hydro-particles, collisional particles, and collisionless particles,
and the definition of three particles will be given in Section \ref{method}.
The hydro-particles are described by the probability density function, while the collisional and collisionless particles are described by the simulation particles.
The integral equation Eq.\eqref{integral-solution-ugks} is solved through the simulation particles, which is coupled with the evolution of macroscopic variables.
The update of macroscopic flow variables is the same as Eq.\eqref{kinetic-m}, but partially use the simulation particles to evaluate the interface flux.
Due to the Lagrangian formulation through particle transport to get the solution of gas distribution function in Eq.\eqref{integral-solution-ugks} directly inside each cell, the Eulerian formulation for the separate update of gas distribution function in Eq.\eqref{kinetic-f} within each cell is not necessary anymore.
One of the distinguishable points in UGKWP method is that the particles are divided into hydro-particles, collisional particles, and collisionless particles.
The dynamics of hydro-particles can be described analytically and its computation is very cost-effective.
The proportion of three kinds of particles varies dynamically in different flow regimes.
Physically, the collisionless particles are mainly used for the description of non-equilibrium transport dynamics and the hydro-particles for the equilibrium one.
The three kinds of particles have a circular relation from collisionless to collisional to hydro-particles and then go back to collisionless particles.
In the continuum flow regime, the number of collisional and collisionless particles will be greatly reduced and
the UGKWP method will automatically converge to the GKS with the same amount of computational cost.
In other words, the numerical flux in Eq.\eqref{kinetic-m} converges to GKS flux in the hydrodynamic regime.
For the simulation of hypersonic flow, the UGKWP method will be much more efficient than the original UGKS due to the use of simulation particles,
which has a nature adaption in particle velocity space.
The computational cost of UGKWP method is similar to the particle methods in rarefied regime and reduces to the hydrodynamic solver in continuum regime,
and such property is summarized as the multi-efficiency preserving scheme in Section \ref{discussion}.

The rest of the paper is organised as following.
In Section \ref{method}, in order to understand the UGKWP method we will first introduce a fully particle description based unified gas-kinetic particle (UGKP) method,
and then the wave-particle description based unified gas-kinetic wave-particle method will be proposed.
Both UGKP and UGKWP methods are multiscale methods for all flow regimes, but UGKWP method is more cost-efficient than UGKP method due to the
analytical formulation for the hydro-particles.
The definition of multiscale multi-efficiency method will be introduced in Section \ref{discussion},
as well as the analysis of the asymptotic properties of UGKWP method.
Numerical tests for flows over a wide range of Mach and Knudsen numbers are shown in Section \ref{numericaltest}.
Section \ref{conclusion} is the conclusion.

\section{Unified Gas-kinetic Wave-Particle method}\label{method}
In this section, two novel multiscale numerical schemes will be introduced under UGKS framework,
i.e. the unified gas-kinetic particle method and the unified gas-kinetic wave-particle method.
The UGKP method is a particle based method and its computational cost and statistical noise keep the same level in different flow regimes.
The UGKWP method improves the UGKP method by decomposing the simulation particles into the particle-described collisional
and collisionless particles, and the wave-described hydro-particles which can be described by an analytical distribution function.
Both methods are built on a discretized physical space $\sum_i \Omega_i\subset\mathcal{R}^3$ and discretized time $t^n\in\mathcal{R^+}$.

\subsection{Unified gas-kinetic particle method}\label{method1}
The particle dynamics in UGKP method is constructed based on the kinetic BGK equation \cite{BGK1954},
\begin{equation}\label{BGK}
\frac{\partial f}{\partial t}+ \vec{v}\cdot\nabla_{\vec{x}} f = \frac{g-f}{\tau},
\end{equation}
where $f(\vec{x},t,\vec{v})$ is the velocity distribution function of gas particle, $\tau$ is the local relaxation parameter which is determined by $\tau=\mu/p$
with the gas pressure $p$ and dynamic viscosity $\mu$.
The local equilibrium Maxwellian distribution $g(\vec{x},t,\vec{v})$ has the form
\begin{equation}\label{maxwellian}
  g(\vec{x},t,\vec{v})=\rho\left(\frac{\lambda}{\pi}\right)^\frac{K+3}{2} \exp(-(\vec{v}-\vec{U})^2+\vec{\xi}^2),
\end{equation}
with density $\rho$, velocity $\vec{U}$, internal degree of freedom $K$, and the internal variable $\vec{\xi}$.
The main idea of UGKP method is to track the particle trajectory until the collision happens.
Once the particle collide with other particles, it will be numerically merged into the macroscopic flow quantities,
and get re-sampled from the updated macroscopic flow variables at the beginning of next time step.
The evolution of particles will be given in Section \ref{section211},
which is coupled with the evolution of macroscopic quantities presented in Section \ref{section212}.

\subsubsection{Evolution of particles}\label{section211}
\begin{figure}
\centering
\includegraphics[width=0.7\textwidth]{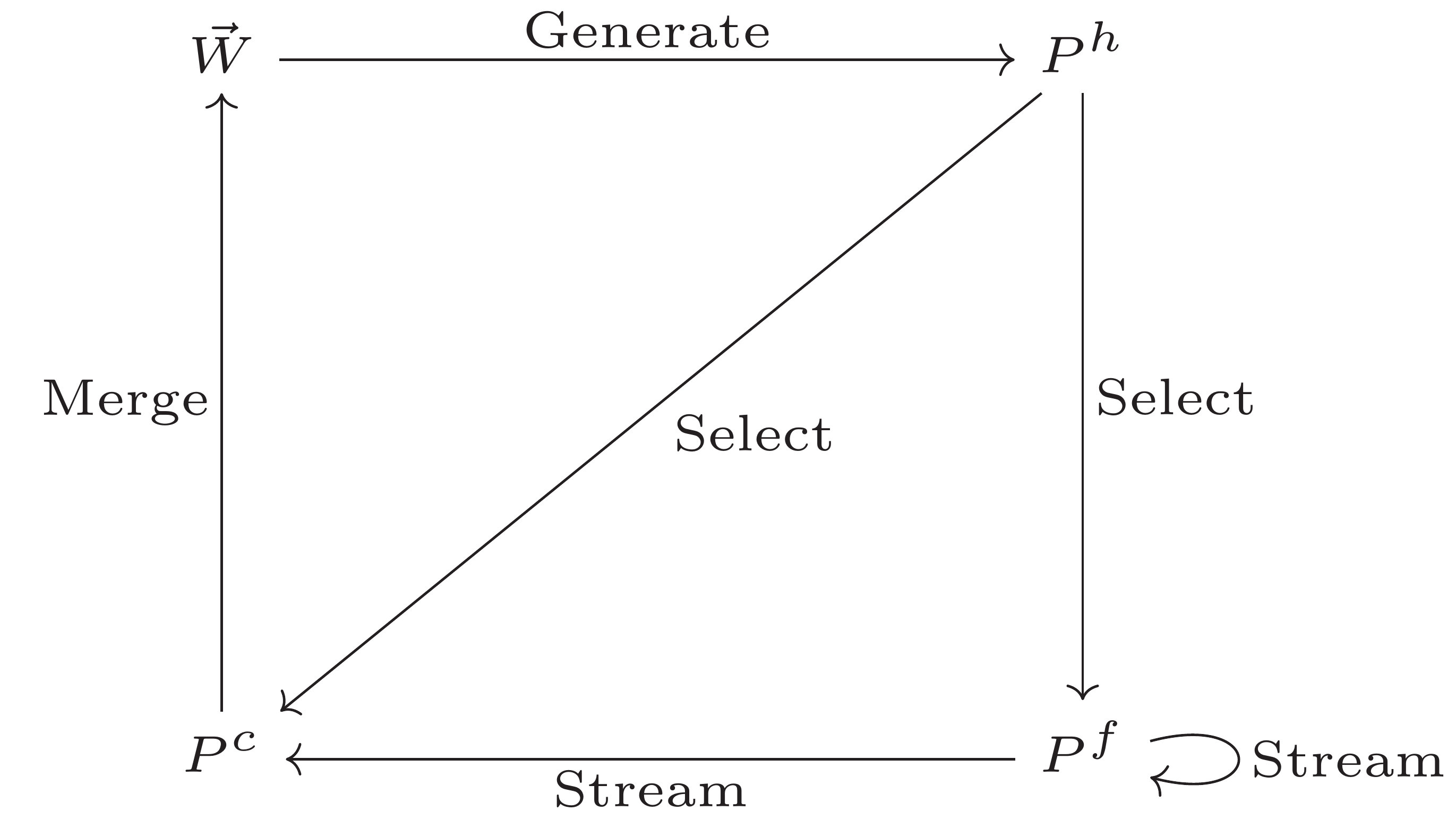}
\caption{Dynamical circulation of particle-described collisionless particle $P^f$, collisional particle $P^c$, wave-described hydro-particle $P^h$, and macroscopic variables for UGKP method.}
\label{dynamical1}
\end{figure}

The simulation particle $P_k(m_k,\vec{x}_k,\vec{v}_k,e_k)$ is represented by its weight $m_k$,
position coordinate $\vec{x}_k$, velocity coordinate $\vec{v}_k$, and internal energy $e_k$,
whose evolution follows the integral form of the kinetic BGK equation,
\begin{equation}\label{integral-solution}
  f(\vec{x},t,\vec{v})=\frac{1}{\tau}\int_{t^n}^t e^{-(t-t')/\tau} g(\vec{x}',t',\vec{v}) dt'+e^{-t/\tau}f_0(\vec{x}_0,\vec{v}),
\end{equation}
where $f_0$ is the initial distribution function at $t=t^n$, and $g$ is the local equilibrium distribution function.
The equilibrium distribution is integrated along the characteristics $\vec{x}'=\vec{x}+\vec{v}(t'-t)$.
Numerically, the equilibrium distribution function can be expanded as
\begin{equation}
  g(\vec{x}',t',\vec{v})=g(\vec{x},t,\vec{v})+\nabla_{\vec{x}} g(\vec{x},t,\vec{v})\cdot(\vec{x}'-\vec{x})+\partial_tg(\vec{x},t,\vec{v})t'+O((\vec{x}'-\vec{x})^2,t'^2),
\end{equation}
following which the integral solution can be expressed as
\begin{equation}\label{particle}\\
  f(\vec{x},t,\vec{v})=(1-e^{-t/\tau}) g^+(\vec{x},t,\vec{v})+e^{-t/\tau}f_0(\vec{x}_0,\vec{v}).
\end{equation}
The first order expansion of $g$ implies
\begin{equation}\label{1st-particle}
  g^+(\vec{x},t,\vec{v})=g(\vec{x},t,\vec{v}),
\end{equation}
and the second order expansion gives
\begin{equation}\label{2nd-particle}
  g^+(\vec{x},t,\vec{v})=g(\vec{x},t,\vec{v})
  +\frac{e^{-t/\tau}(t+\tau)-\tau}{1-e^{-t/\tau}}(\partial_tg(\vec{x},t,\vec{v})+\vec{v}\cdot\nabla_{\vec{x}}g(\vec{x},t,\vec{v})).
\end{equation}
Above $g^+$ is named the hydrodynamic distribution function with analytical formulation.
For UGKP method, the first order expansion of $g$ is used for a simple particle-sampling algorithm \cite{bird1994molecular}.
The particle evolution equation Eq.\eqref{particle} means that the simulation particle has a probability of $e^{-t/\tau}$ to free stream, and has a probability of $(1-e^{-t/\tau})$ to collide with other particles and follow the velocity distribution $g^+(\vec{x},t,\vec{v})$.
The time, when one simulation particle stops free streaming and follows the distribution $g^+$, is defined as its `first collision time' $t_c$.
The cumulative distribution function of the first collision time is
\begin{equation}\label{tc-distribution}
  F(t_c<t)=1-\exp(-t/\tau),
\end{equation}
from which $t_c$ can be sampled as $t_c=-\tau\ln(\eta)$ with $\eta$ a uniform distribution $\eta\sim(0,1)$.
From a simulation time step $t^{n}$ to $t^{n+1}$, all simulation particles in UGKP method can be categorized into two groups: the `collisionless particle' \index{collisionless particle} $P^f$ and the `collisional particle' \index{collisional particle} $P^c$.
The categorization is based on the relation between the first collision time $t_c$ and the time step $\Delta t$.
More specifically, the collisionless particle is defined as the particle whose first collision time $t_c$ greater than or equal to the time step $\Delta t$,
and the collisional particle is defined as the particle whose first collision time $t_c$ smaller than $\Delta t$.
For the collisionless particle, its trajectory is fully tracked during the whole time step.
For collisional particle, the particle trajectory is tracked till $t_c$.
Then the particle's mass, momentum, and energy are merged into the macroscopic quantities in that cell and the simulation particle gets eliminated.
Those eliminated particles will get re-sampled once the updated macroscopic quantities $\vec{W}^{n+1}$ are obtained.
As shown in Eq.\eqref{particle}, the re-sampled particles follow the hydrodynamic distribution $g^+$ and therefore they are defined as `hydro-particle' \index{hydro-particle} $P^h$. The macroscopic quantities corresponding to the hydro-particles are defined as `hydro-quantities' \index{hydro-quantities} $\vec{W}^h$.
The hydro-particles will be sampled at the beginning of each time step and become the candidates for collisionless/collisional particles
again in the next time step evolution according to their newly-sampled $t_c$.
The dynamical circulation of particle-described collisionless particle $P^f$, collisional particle $P^c$, wave-described hydro-particle $P^h$, and macroscopic variables is shown in Fig. \ref{dynamical1}, and the algorithm for the evolution of particles is presented as following.

\begin{description}
  \item[Step 1] Sample the first collision time $t_c$ for all particles $P_k$. For example, if $P_k\in\Omega_i$, then $t_{c,k}=-\tau_i\ln(\eta)$, where $t_{c,k}$ is the first collision time for $P_k$ and $\tau_i$ is calculated from the cell averaged macroscopic quantities $\vec{W}_i$ in cell $i$,
      \begin{equation}
      \vec{W}_i\defeq\frac{1}{|\Omega_i|}\int_{\Omega_i} \vec{W} dx;
      \end{equation}
  \item[Step 2] Steam collisionless particles $P^{f,n}$ to $P^{f,n+1}$ by $\vec{x}_k^{n+1}=\vec{x}_k^n+\vec{v}_k (t^{n+1}-t^n)$;
  \item[Step 3] Update the cell averaged variables $\vec{W}^{n+1}_i$ by Eq.\eqref{kinetic-m}, and calculate the hydro-quantities
  $\vec{W}_{i}^h=(\rho_i^h,\rho_i^hU_i^h,\rho_i^hV_i^h,\rho_i^hW_i^h,\rho_i^hE_i^h)$ by
  \begin{equation}
  \vec{W}_i^h=\vec{W}_i^{n+1}-\left(\vec{W}_i^{P^f}\right)^{n+1},
  \end{equation}
  Here $\left(\vec{W}_i^{P^f}\right)^{n+1}$ is the cell averaged macroscopic quantities from all remaining collisionless particles in cell i,
  \begin{equation}\label{wf}
  \left(\vec{W}_i^{P^f}\right)^{n+1}=\frac{1}{|\Omega_i|}\sum_k \left( m^f_k, m^f_k u^f_k, m^f_k v^f_k,
   m^f_k w^f_k,\frac12 m^f_k \left(|\vec{v}^f_k|^2+e^f_k\right)\right)^T,
  \end{equation}
  where the index $k$ covers all collisionless particles $P_k^f$ in cell $i$.
  The detailed formulation of the evolution of macroscopic quantities will be given in next subsection. Note that the calculation of hydro-quantities will also be used in UGKWP method.
  \item[Step 4] Sample hydro-particles from a $\vec{W}^{n+1}_i$-based Maxwellian distribution with a total mass of $\Omega_i \rho_i^h$, and the physical coordinates of these hydro-particles are uniformly distributed in cell i.
\end{description}

\subsubsection{Evolution of macroscopic flow variables}\label{section212}
The evolution of simulation particles is coupled with the evolution of macroscopic quantities
$\vec{W}=\left(\rho,\rho U,\rho V,\rho W,\rho E\right)^T$, where $\rho$ is density, $\vec{U}=(U,V,W)^T$ is macroscopic velocity, and $E$ is energy per unit mass.
The cell averaged macroscopic variables $\vec{W}_i$ are evolved by the macroscopic governing equation of UGKS
\begin{equation}\label{update-w}
  \vec{W}_i^{n+1}=\vec{W}_i^{n}-\frac{1}{|\Omega_i|} \sum_{l_s\in\partial\Omega} |l_s|\vec{F}_{s},
\end{equation}
where the UGKS flux for the macroscopic variables are
\begin{equation}\label{flux-w-ugks}
  \vec{F}_{s}=\int_{t^n}^{t^{n+1}}\int\bigg[\frac1\tau\int_0^{t} e^{(t'-t)/\tau}g(\vec{x}'_s,t',\vec{v})dt'+
  e^{-t/\tau}f_0(\vec{x}_s-\vec{v}t,\vec{v})\bigg] \vec{v}\cdot \vec{n}_s \vec{\psi} d\Xi dt,
\end{equation}
with the outward normal $\vec{n}_s$ and the characteristics $\vec{x}'_s=\vec{x}_s+\vec{v}(t'-t)$.
The flux terms related to the Maxwellian distribution are denoted as $F_{g,s}$
\begin{equation}\label{Fg}
  \vec{F}_{g,s}\defeq\int_{t^n}^{t^{n+1}}\int\frac1\tau\int_0^{t} e^{(t'-t)/\tau}g(\vec{x}'_s,t',\vec{v})dt' \vec{v}\cdot\vec{n}_s\vec{\psi} d\Xi dt,
\end{equation}
and the flux terms related to the initial distribution are $\vec{F}_{f,s}$
\begin{equation}\label{Ff}
  \vec{F}_{f,s}\defeq\int_{t^n}^{t^{n+1}}\int e^{-t/\tau} f_0(\vec{x}_s-\vec{v}t,\vec{v}) \vec{v}\cdot \vec{n}_s\vec{\psi} d\Xi dt,
\end{equation}
where $\vec{\psi}$ is the vector of conservative moments
\begin{align*}
  \vec{\psi}=\left(1,u,v,w,\frac12(\vec{v}^2+\vec{\xi}^2)\right)^T.
\end{align*}
The calculation of the Maxwellian-related terms are the same as UGKS \cite{ugks2010}.
Assume the interface is located at $\vec{x}_0$ with a local coordinate $(\vec{e}_1,\vec{e}_2,\vec{e}_3)$ and $\vec{e}_1$ is the outward unit normal $\vec{n}_0$.
The Maxwellian distribution is expanded around $\vec{x}_0$ as
\begin{equation}\label{g-expansion1}
\begin{aligned}
  g(\vec{x},t,\vec{v})=&g_0(\vec{x}_0,\vec{v})
  +(1-H[\bar{x}])\frac{\partial^l}{\partial\vec{e}_1} g_0(\vec{x}_0,\vec{v})\bar{x}
  +H[\bar{x}]\frac{\partial^r}{\partial\vec{e}_1}  g_0(\vec{x}_0,\vec{v})\bar{x}\\
  &+\frac{\partial}{\partial\vec{e}_2}g_0(\vec{x}_0,\vec{v})\bar{y}
  +\frac{\partial}{\partial\vec{e}_3}g_0(\vec{x}_0,\vec{v})\bar{z}
  +\frac{\partial}{\partial t}g_0(\vec{x}_0,\vec{v})(t-t^n) \\
  =&g_0(\vec{x}_0,\vec{v})\left[1+(1-H[\bar{x}])a^l\bar{x}+H[\bar{x}]a^r\bar{x}+b\bar{y}+c\bar{z}+A(t-t^n)\right],
\end{aligned}
\end{equation}
where $g_0(\vec{x}_0,\vec{v})=g(\vec{x}_0,t^n,\vec{v})$, and $\bar{x}=\Delta x\cdot \vec{e}_1$, $\bar{y}=\Delta y\cdot \vec{e}_2$, and $\bar{z}=\Delta z\cdot \vec{e}_3$.
The derivative functions of Maxwellian distribution $a^l$, $a^r$, $b$, $c$, and $A$ have the following form
\begin{align*}
  a^l&=a^l_1+a^l_2u+a^l_3v+a^l_4w+a^l_5\frac12(\vec{v}^2+\vec{\xi}^2)=a_\alpha^l\psi_\alpha,\\
  a^r&=a^r_1+a^r_2u+a^r_3v+a^r_4w+a^r_5\frac12(\vec{v}^2+\vec{\xi}^2)=a_\alpha^r\psi_\alpha,\\
  &\qquad\qquad\qquad...\\
  A&=A_1+A_2u+A_3v+A_4w+A_5\frac12(\vec{v}^2+\vec{\xi}^2)=A_\alpha\psi_\alpha.
\end{align*}
The heaviside function $H[x]$ is
\begin{equation}\nonumber
H[x]=\left\{
\begin{aligned}
  1 \quad x>0,\\
  0 \quad x\le0.
\end{aligned}\right.
\end{equation}
The Maxwellian at $\vec{x}_0$ and its derivative functions can be obtained from the reconstructed macroscopic variables.
In this paper, the van Leer limiter is used for reconstruction,
\begin{equation}
  s=(\text{sign}(s_l)+\text{sign}(s_r))\frac{|s_l||s_r|}{|s_l|+|s_r|},
\end{equation}
where $s$, $s_l$, and $s_r$ are the slopes of macroscopic variables.
The Maxwellian distribution at cell interface can be obtained from
\begin{equation}
  \vec{W}_{0}=\int \vec{\psi} \left(g^l_{0}H[\bar{u}]+g^r_{0}(1-H[\bar{u}])\right)d\Xi,
\end{equation}
where $\vec{W}_{0}$ is the macroscopic variables at $x_0$ corresponding to $g_0$, and $\bar{u}=\vec{u}\cdot\vec{e}_1$.
The derivative functions $a^l,a^r,b,c,A$ are calculated from the spatial and time derivatives of $g_0$, Taking $A$ as an example,
\begin{align}
A&=\frac{1}{g_0}\bigg(\frac{\partial g_0}{\partial \vec{W}_0}\bigg)\bigg(\frac{\partial \vec{W}_0}{\partial t}\bigg)_{t=t^n},
\end{align}
and
\begin{align}
&A_5=\frac{\rho}{3p^2}\left(2\frac{\partial \rho E}{\partial t}+\left(U_iU_i-\frac{3p}{\rho}\right)\frac{\partial \rho}{\partial t}-2U_i\frac{\partial \rho U_i}{\partial t}\right),\\
&A_{i+1}=\frac{1}{p}\left(\frac{\partial \rho U_i}{\partial t}-U_i\frac{\partial \rho}{\partial t}\right)-U_iA_5 \quad (i=1,2,3),\\
&A_1=\frac{1}{\rho}\frac{\partial \rho}{\partial t}-U_ia_{i+1}-\frac12\left(U_iU_i+\frac{3p}{\rho}\right)A_5,
\end{align}
where the macroscopic quantities are those at $(\vec{x}_0,t^n)$. The time derivatives of macroscopic variables are determined by the conservative moments requirements on the first order Chapman-Enskog expansion \cite{chapman1990mathematical}, which gives
\begin{equation}
\left(\frac{\partial\vec{W}_0}{\partial t}\right)=-\int \left(a^l\bar{u}H[\bar{u}]+a^r\bar{u}(1-H[\bar{u}])+b\bar{v}+c\bar{w}\right)g_0 \vec{\psi} d\Xi.
\end{equation}
Readers can refer to \cite{wang2015unified} for derivation and detailed formulae.
Once the Maxwellian distribution at cell interface and its derivative functions are determined, the partial flux function Eq.\eqref{Fg} can be obtained using the expansion Eq.\eqref{g-expansion1} for the interface distribution function, which gives
\begin{equation}\label{Fg1}
\begin{aligned}
\vec{F}_{g,s}=\int \vec{v}\cdot \vec{n}_s \vec{\psi}\bigg\{&\left(\tau e^{-\Delta t/\tau}+\Delta t-\tau\right)g_0(\vec{x}_0,\vec{v})\\
&+\tau\left(-e^{-\Delta t/\tau}(\Delta t+2\tau)-\Delta t+2\tau\right)
\left[a^lH[\bar{u}]+a^r(1-H[\bar{u}])\right]\bar{u} g_0(\vec{x}_0,\vec{v})\\
&+\tau\left(e^{-\Delta t/\tau}(\Delta t+2\tau)-\Delta t+2\tau\right)(b\bar{v}+c\bar{w})g_0(\vec{x}_0,\vec{v})\\
&+\left(-\tau^2e^{-\Delta t/\tau}+\Delta t^2/2-\tau\Delta t+\tau^2\right)Ag_0(\vec{x}_0,\vec{v})\bigg\} d\Xi.
\end{aligned}
\end{equation}

The free transport terms in UGKS flux are calculated from the simulation particles.
The net flux of cell $i$ contributed by the free transport terms are
\begin{equation}
 \vec{F}_{f,i}=\left(\vec{W}_i^{P}\right)^{n+1}-\left(\vec{W}_i^{P}\right)^{n},
\end{equation}
where $\vec{W}_i^{P}$ is the vector of the macroscopic quantities of all particles in cell $i$.
The updated $\left(\vec{W}_i^{P}\right)^{n+1}$,
\begin{equation}
\left(\vec{W}_i^{P}\right)^{n+1}=\sum_k\left( m_k,m_k u_k, m_k v_k,
 m_k w_k, \frac12 m_k \left(|\vec{v}_k|^2+e_k\right)\right)^T,
\end{equation}
only counts the collisionless particles at the end of this time step.
Note that those particles differ from the total particles at the beginning of next time step which will also include the newly sampled hydro-particles.
For UGKP method, we have $\vec{W}_i^n=\left(\vec{W}_i^{P}\right)^{n}$, and the evolution equation of the macroscopic variables follows
\begin{equation}
\begin{aligned}
\vec{W}^{n+1}=&\vec{W}^n+\frac{1}{|\Omega_i|} \left(\sum_{l_s\in\partial \Omega}|l_s| \vec{F}_{g,s}+\vec{F}_{f,i}\right)\\
=&\left(\vec{W}_i^{P}\right)^{n+1}+\frac{1}{|\Omega_i|} \sum_{l_s\in\partial \Omega}|l_s| \vec{F}_{g,s}.
\end{aligned}
\end{equation}

\subsection{Unified gas-kinetic wave-particle method}\label{method2}
The UGKWP method improves UGKP method mainly in the following two aspects:
\begin{itemize}
  \item The free transport terms in numerical flux contributed by the hydro-particles are evaluated analytically;
  \item Only collisionless hydro-particles are sampled.
\end{itemize}

Firstly, since the distribution of the hydro-quantities is known as $g^+$, we can analytically evaluate the flux contributed by the free transport of hydro-particles, which gives
\begin{equation}
\sum_s |l_s| \vec{F}_{f,s}=\sum_s |l_s| \vec{F}_{f,s}^h+\vec{F}_{f,i}^p,
\end{equation}
where $\vec{F}_{f,s}^h$ is the free transport flux contributed by the hydro-quantities,
\begin{equation}\label{Ff1}
\vec{F}_{f,s}^h=\int_{t^n}^{t^{n+1}}\int\bigg\{e^{-t'/\tau}\left[g_0^{+,l}(\vec{x}-\vec{v}t',\vec{v})H[\bar{u}]
+g_0^{+,r}(\vec{x}-\vec{v}t',\vec{v})(1-H[\bar{u}])\right]\bigg\} \vec{v}\cdot\vec{n}_s \vec{\psi} d\Xi dt.
\end{equation}
Here the second order expansion of $g^+$ Eq.\eqref{2nd-particle} is used.
The numerical flux contributed by the streaming of collisionless and collisional particle, i.e. $\vec{F}_{f,i}^p$, will be given later.

Secondly, since $\vec{F}_{f,s}^h$ is analytically evaluated, there is no need to sample all hydro-particles.
Only the collisionless hydro-particles \index{collisionless hydro-particles} will be sampled.
Based on the cumulative distribution function of the first collision time Eq.\eqref{tc-distribution},
the collisionless hydro-particles are sampled from a $\vec{W}^{n+1}_i$-based Maxwellian distribution with the total mass of $e^{-t/\tau_i}\Omega_i\rho^h_i$.
Different from UGKP method, for the UGKWP method, simulation particles are categorized into three groups: the collisionless particle $P^f$, the collisional particle $P^c$, and the collisionless hydro-particle $P^{hf}$.
The flux contribution from collisionless hydro-particles can be evaluated analytically. The net flux $\vec{F}_{f,i}^p$
contributed by the collision and collisionless particles is
\begin{equation}\label{Fp}
 \vec{F}_{f,i}^P=\left(\vec{W}_i^{fc}\right)^{n+1}-\left(\vec{W}_i^{fc}\right)^{n},
\end{equation}
where $\vec{W}_i^{fc}$ is the macroscopic quantities of collisionless particle $P^f$ and collisional particle $P^c$ in cell $i$.
The evolution of the macroscopic variables for UGKWP method is
\begin{equation} \label{WF}
\vec{W}^{n+1}=\vec{W}^n+\frac{1}{|\Omega_i|} \left[\sum_{l_s\in\partial \Omega}|l_s| \left(\vec{F}_{g,s}+\vec{F}_{f,s}^h\right)+\vec{F}_{f,i}^P\right].
\end{equation}

\begin{figure}
\centering
\includegraphics[width=0.7\textwidth]{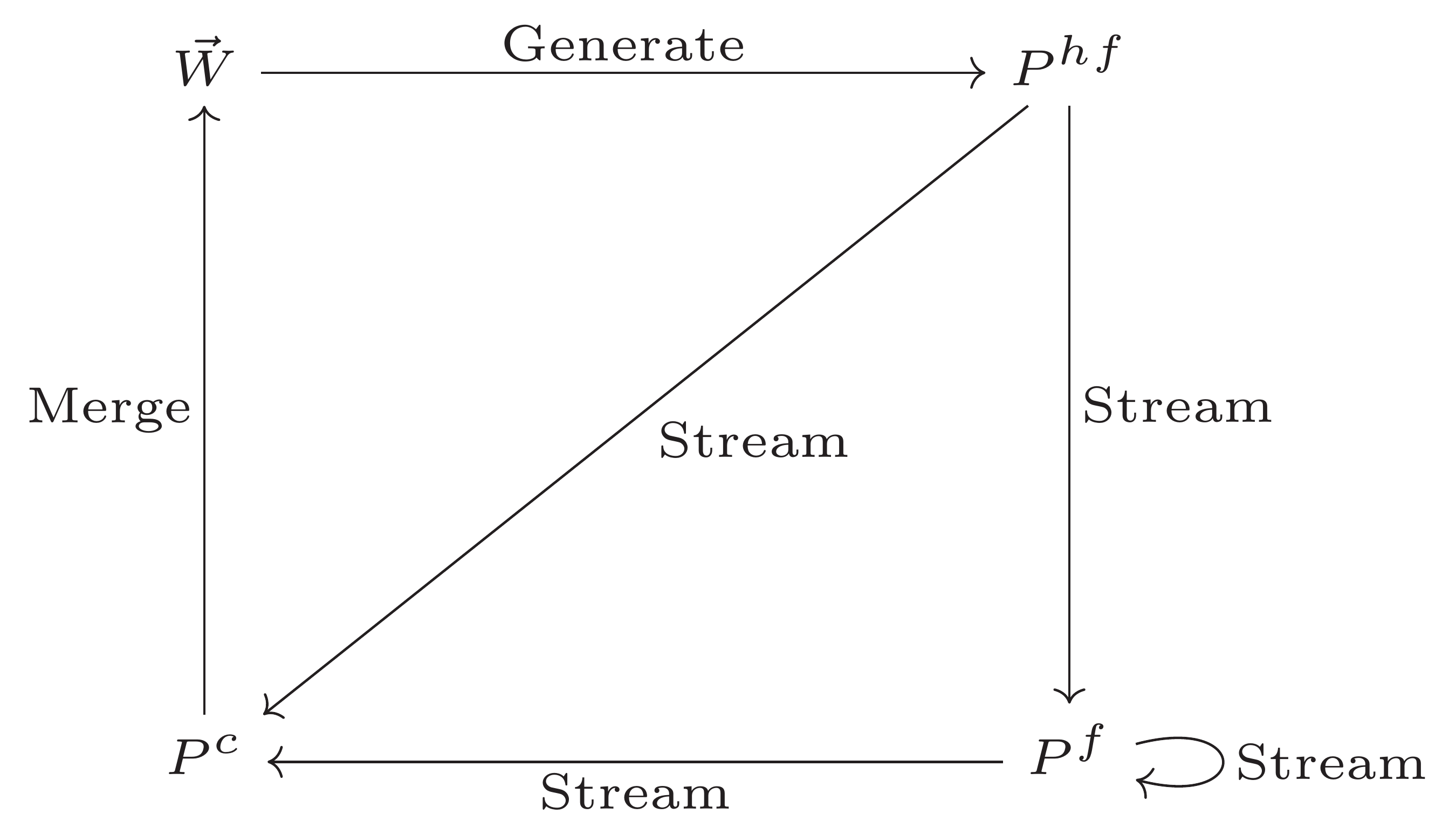}
\caption{Dynamical circulation of particle-described collisionless particle $P^f$, collisional particle $P^c$, wave-described free-transport hydro-particle $P^{hf}$, and macroscopic variables for UGKWP method.}
\label{dynamical2}
\end{figure}

\begin{figure}
\centering
\includegraphics[width=0.9\textwidth]{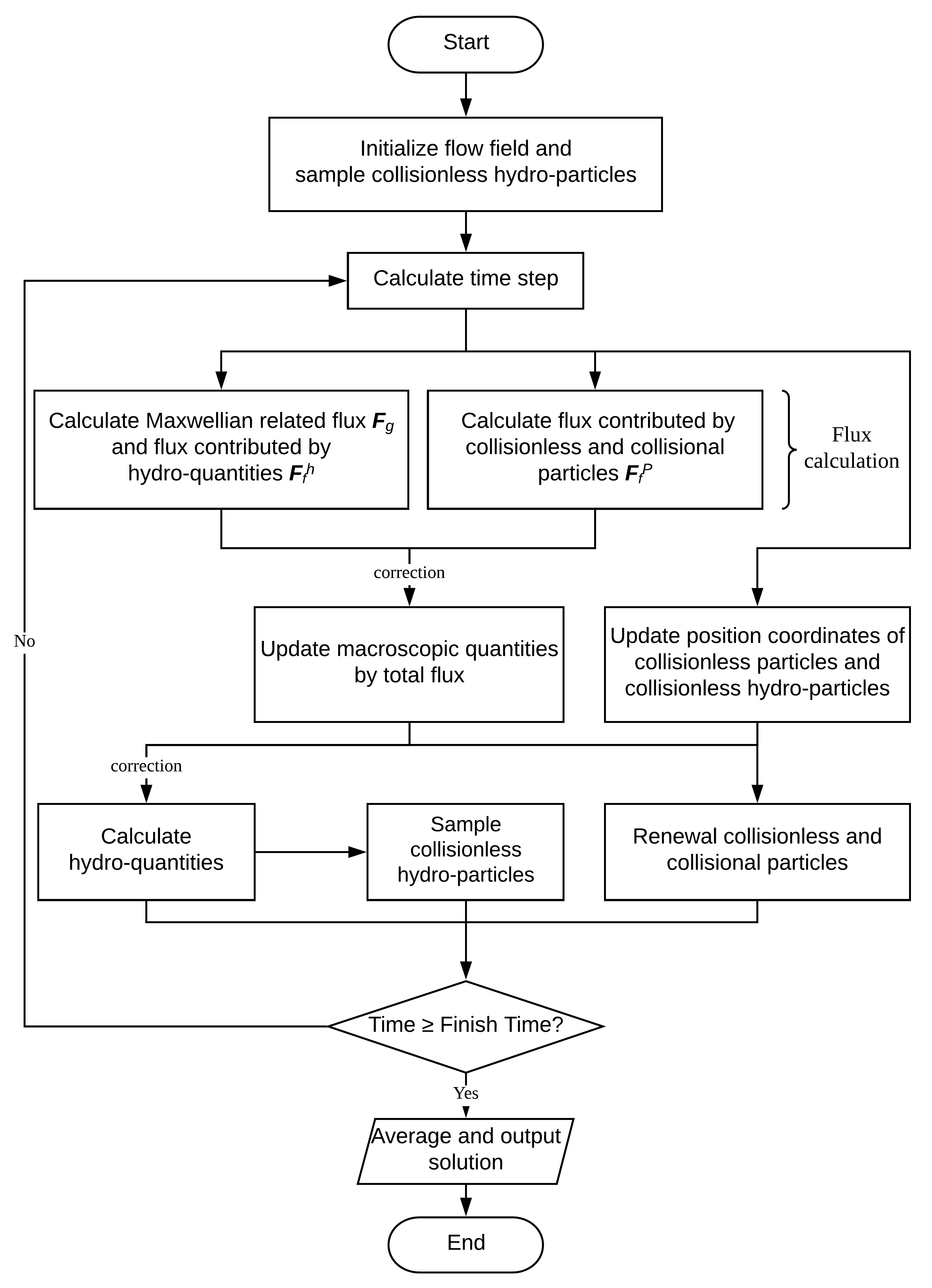}
\caption{Flow chart of the UGKWP method. In the flux calculation, the Maxwellian related flux $\vec{F}_{g,s}$ is calculated by Eq.\eqref{Fg1}; the flux contributed by hydro-quantities $\vec{F}_{f,s}^h$ is calculated by Eq.\eqref{Ff1}; and flux contributed by collisionless and collisional particles $\vec{F}_{f,i}^P$ is calculated by Eq.\eqref{Fp}. The two corrections will be given in the discussion of robustness in Section \ref{correction}.}
\label{flowchart}
\end{figure}

For the evolution of simulation particles, the physical coordinates of the collisionless particles $P^f$ and the collisionless hydro-particles $P^{hf}$ are updated by
\begin{equation}
 \left(\vec{x}^{f,h}\right)^{n+1}=\left(\vec{x}^{f,h}\right)^n+\vec{v}^{f,h}\Delta t,
\end{equation}
and the collisional particles are eliminated after flux calculation.
The updated collisionless particles and collisionless hydro-particles gather to be the candidates of collisionless and collisional particles in the next time step calculation according to their newly-sampled $t_c$.
In summary, the dynamical circulation of particle-described collisionless particle $P^f$, collisional particle $P^c$, wave-described free-transport hydro-particle $P^{hf}$, and macroscopic variables is shown in Fig \ref{dynamical2}, and the flow chart of the UGKWP method is given in Fig. \ref{flowchart}.

\section{Analysis and discussion}\label{discussion}
\subsection{Asymptotic behavior in continuum regime}
In this section, the asymptotic behavior and computational cost of the UGKWP method
will be analyzed in continuum regime when the time step $\Delta t$ is much larger than the local relaxation parameter $\tau$.
When $\Delta t \gg \tau$, for cell $i$, the total mass of the sampled collisionless hydro-particles is
$M^h=e^{-\Delta t/\tau_i}\Omega_i \rho_i^h$, and the total mass of the collisionless and collisional particles $M^p$ is proportional to $M^h$,
i.e. $M^p\sim O(e^{-\Delta t/\tau_i})$.
Therefore, the numerical flux contribution by collisionless and collisional particle streaming is $\vec{F}_{f,i}^p\sim O(e^{-\Delta t/\tau_i})$.
In the free transport flux of hydro-quantities $\vec{F}_{f,s}^h$ given by Eq.\eqref{Ff1},
 the hydrodynamic distribution functions $g_0^{+,l,r}(\vec{x},\vec{v})$ become
\begin{equation}\label{glimit}
\begin{aligned}
  g_0^{+,l,r}(\vec{x},\vec{v})=&g_0^{l,r}(\vec{x},\vec{v})-\tau
  \left(\partial_t g_0^{l,r}(\vec{x},\vec{v}),\vec{v}\cdot \nabla_{\vec{x}}g_0^{l,r}(\vec{x},\vec{v})\right)+O(e^{-\Delta t/\tau})\\
  =&f_{\text{NS}}^{l,r}+O(e^{-\Delta t/\tau}),
\end{aligned}
\end{equation}
where $f_{\text{NS}}^{l,r}$ is the local first order expansion in Chapman-Enskog asymptotic series.
Substituting Eq.\eqref{glimit} into Eq.\eqref{flux-w-ugks}, the total flux of the macroscopic variables in Eq.\eqref{WF}
converges to the flux of gas-kinetic scheme \cite{xu2001gas},
\begin{equation}\label{fluxlimit}
\vec{F}_{s}=\int_{t^n}^{t^{n+1}}\int \left[\frac{1}{\tau} \int_0^t e^{(t'-t)/\tau} g(\vec{x}'_s,t',\vec{v}) dt'+e^{-t/\tau} f_{\text{NS},0}(\vec{x}_s-\vec{v}t,\vec{v})\right]\vec{v}\cdot\vec{n}_s\vec{\psi}d\Xi dt+O(e^{-\Delta t/\tau}).
\end{equation}
For a well resolved flow region with $g_0^l=g_0^r$ and $\nabla_{\vec{x}}^l g=\nabla_{\vec{x}}^r g$, it can be derived from Eq.\eqref{fluxlimit} and Eq.\eqref{Fg1} that the total flux becomes
\begin{equation}
\vec{F}_{s}=\int_{t^n}^{t^{n+1}}\int \big[ g_0-\tau\left(\partial_t g_0+\vec{v}\cdot\nabla_{\vec{x}}g\right)\big] \vec{v}\cdot\vec{n}_s\vec{\psi}d\Xi dt,
\end{equation}
which gives consistent flux to the first order Chapman-Enskog expansion of kinetic equation for the NS solutions.
From the above analysis, it is concluded that the UGKWP method preserves the hydrodynamic Navier-Stokes equations in continuum regime for $\Delta t \gg \tau$.

In the continuum regime with $\Delta t \gg \tau$, for a fixed particle mass $m_p$, the number of sampled collisionless hydro-particles in cell $i$ is $e^{-\Delta t/\tau_i}\Omega_i \rho_i^h/m_p$. Therefore, the total simulation particle number $N_p$ in such regime decreases exponentially, $N_p\sim O(e^{-\Delta t/\tau})$, which means the computational cost of the UGKWP method becomes comparable to hydrodynamic NS solvers.

We define a numerical scheme as a `Multiscale Multi-efficiency Preserving' (MMP) \index{Multiscale multi-efficiency preserving method} scheme if the following two constraints are satisfied:
\begin{enumerate}
  \item The scheme preserves all flow regime solutions, i.e. from collisionless regime to hydrodynamic Navier-Stokes regime and Euler regime, and the cell size and time step are not constrained to be less than the mean free path and collision time.
  \item The computational efficiency of this scheme is comparable to the high efficient schemes in all flow regimes, for example the computational cost including the computational time and memory cost is comparable to the NS solvers in continuum regime. At the same time, for highly non-equilibrium hypersonic flow, the efficiency goes to the purely particle method.
\end{enumerate}
It is shown in above analysis that UGKWP method is a multiscale multi-efficiency preserving method.

\subsection{Consistent sampling}\label{consistency}
For both UGKP and UGKWP methods, we need to sample particles from a given Maxwellian distribution.
For cell $\Omega_i$ with sampling quantities $\vec{W}_{s,i}=\left(\rho_{s,i},U_{s,i},V_{s,i},W_{s,i},\rho_{s,i}E_{s,i}\right)^T$,
hydro-particles are sampled from the Maxwellian distribution
\begin{equation}
f_{s,i}=\rho_{s,i}\left(\frac{\lambda_{s,i}}{\pi}\right)^{\frac{K+3}{2}}
\exp\left\{-\lambda_{s,i}\left[\left(\vec{v}-\vec{U}_{s,i}\right)^2+\vec{\xi}^2\right]\right\},
\end{equation}
where $\lambda_{s,i}=\frac{\rho_{s,i}}{2(\gamma-1)\rho_{s,i}e_{s,i}}$, $\rho_{s,i}e_{s,i}=\rho_{s,i}E_{s,i}-\frac12\rho_{s,i}\vec{U}_{s,i}^2$, and $\gamma=\frac{K+5}{K+3}$.
The sampled particles $P^s_k$, $k=1,...,N_s$, follows
\begin{equation}
\begin{aligned}
&m_p^s=\frac{\rho_{s,i}|\Omega|_i}{N_s}, \quad \vec{x}_{k}^s\sim U(\omega_i), \quad e_k^s=e_{s,i}, \\
&\vec{v}_k^s=(-\ln (\vec{\eta_1})/\lambda_{s,i})^{1/2}\cos(\vec{\eta}_2), \quad \vec{\eta}_{1,2}\sim U(0,1)^3,
\end{aligned}
\end{equation}
where $U(\omega_i)$ is the uniform distribution on $\Omega_i$ and $U(0,1)^3$ is the uniform distribution on $(0,1)^3$.
A velocity transformation
\begin{equation}
\vec{v}_k^{s\prime}=b\left(\vec{v}_k^{s}-\vec{a}-\vec{U}_{s,i}\right)+\vec{U}_{s,i}
\end{equation}
is required to make the macroscopic quantities of the sampled particles $P_k^s$ consistent with the sampling quantities $\vec{W}_{s,i}$.
The parameters $\vec{a}$ and $b$ are solved from the consistent constraints
\begin{equation}
\left\{
\begin{aligned}
  &\sum_k m_k^s \vec{v}_k^{s\prime}= \rho_{s,i}\vec{U}_{s,i}|\Omega_i|,\\
  &\sum_k\left(\frac12 m_k^s \left(\vec{v}_k^{s\prime}\right)^2+m_k^se_k^s\right)=\rho_{s,i}E_{s,i}|\Omega_i|,
\end{aligned}
\right.
\end{equation}
which give
\begin{equation}
\left\{
\begin{aligned}
  &\vec{a}=\frac{1}{N_s}\sum_k\vec{v}_{k}^s-\vec{U}_{s,i},\\
  &b=\frac{-c_2\pm\sqrt{c_2^2-4c_1c_3}}{2c_1},
\end{aligned}
\right.
\end{equation}
where $c_1=\sum_k\frac12\left(\vec{v}_k^s-\vec{a}-\vec{U}_{s,i}\right)^2$,
$c_2=\sum_k\vec{U}_{s,i}\vec{v}_k^s-N_s(\vec{a}+\vec{U}_{s,i})\vec{U}_{s,i}$, and
$c_3=\frac12 N_s\vec{U}_{s,i}^2+\sum_k e_k^s-E_{s,i}$.

\subsection{Robustness and conservation}\label{correction}
Due to statistical noise, the density and temperature of the cell averaged macroscopic quantities $\vec{W}_i$ and hydro-quantities $\vec{W}_i^h$ may become negative.
Therefore two corrections are required. The first correction of
\begin{equation}\label{correct1}
\tilde{\rho}_i=\max\left(\rho_i,0\right),
\tilde{\rho}_i\tilde{E}_i=\max\left(\rho_iE_i,\frac12\rho_i\vec{U}_i^2\right)
\end{equation}
is put in the update of cell averaged macroscopic variables, and the second one
\begin{equation}\label{correct2}
\tilde{\rho}_i^h=\max\left(\rho_i^h,0\right),
\tilde{\rho}_i^h\tilde{E}_i^h=\max\left(\rho_i^hE_i^h,\frac12\rho_i^h(\vec{U}_i^h)^2\right)
\end{equation}
is put in the calculating of the hydro-quantities.
The UGKWP method is a finite volume scheme for the macroscopic variables. As long as Eq.\eqref{correct1} doesn't take effect, which means $\tilde{\rho}_i=\rho_i,\tilde{\rho}_i\tilde{E}_i=\rho_iE_i$, the mass, momentum, and energy conservation will be satisfied.
For all the numerical tests in this paper, Eq.\eqref{correct1} does not take effect at all.
It only gives a safe-guided warrantee once it is needed.

\section{Numerical tests}\label{numericaltest}
In this section, five numerical tests are calculated to demonstrate the multiscale property, and to present the computational efficiency of the UGKWP method.
The test cases are one dimensional Sod shock tube, normal shock wave, two dimensional flow passing a cylinder,
lid-driven cavity flow, and  flat plate boundary layer.
The solution of UGKWP method shows good agreement with the reference solution for both high speed and low speed flow in all flow regimes.
For hypersonic flow, the UGKWP method shows a much higher efficiency and lower memory cost than the conventional UGKS with a
direct discretization of particle velocity space.
In the near continuum flow regime, the UGKWP method shows a much higher efficiency and low statistical noise than conventional particle methods.
In the following numerical tests, the reference solution of kinetic equation is given by UGKS and the reference NS solution is given by GKS.
The code is sequential and the computation is carried out on a computer with Intel i7-8700K CPU, 64 GB memory.

\subsection{Sod shock tube problem}
We first calculate the Sod shock tube problem with difference Knudsen number to show the performance of the UGKWP method in different flow regime, and two set of simulation particle number is used to study the statistical noise.
In this calculation, the dimensionless quantities are used. The computational domain is [-0.5,0.5], with initial condition
\begin{equation}\nonumber
(\rho, u, p)=\left\{
\begin{aligned}
  &(1.0,0,1.0) \qquad  x\le0,\\
  &(0.125,0,0.1) \quad x>0.
\end{aligned}\right.
\end{equation}
The viscous coefficient is given as
\begin{equation}\label{vhs-vs1}
  \mu=\mu_{ref}\left(\frac{T}{T_0}\right)^{\omega},
\end{equation}
with the temperature dependency index $\omega=0.81$, and the reference viscosity
\begin{equation}\label{vhs-vs2}
\mu_{ref}=\frac{15\sqrt{\pi}}{2(5-2\omega)(7-2\omega)}\mathrm{Kn}.
\end{equation}

\begin{figure}
\centering
\includegraphics[width=0.48\textwidth]{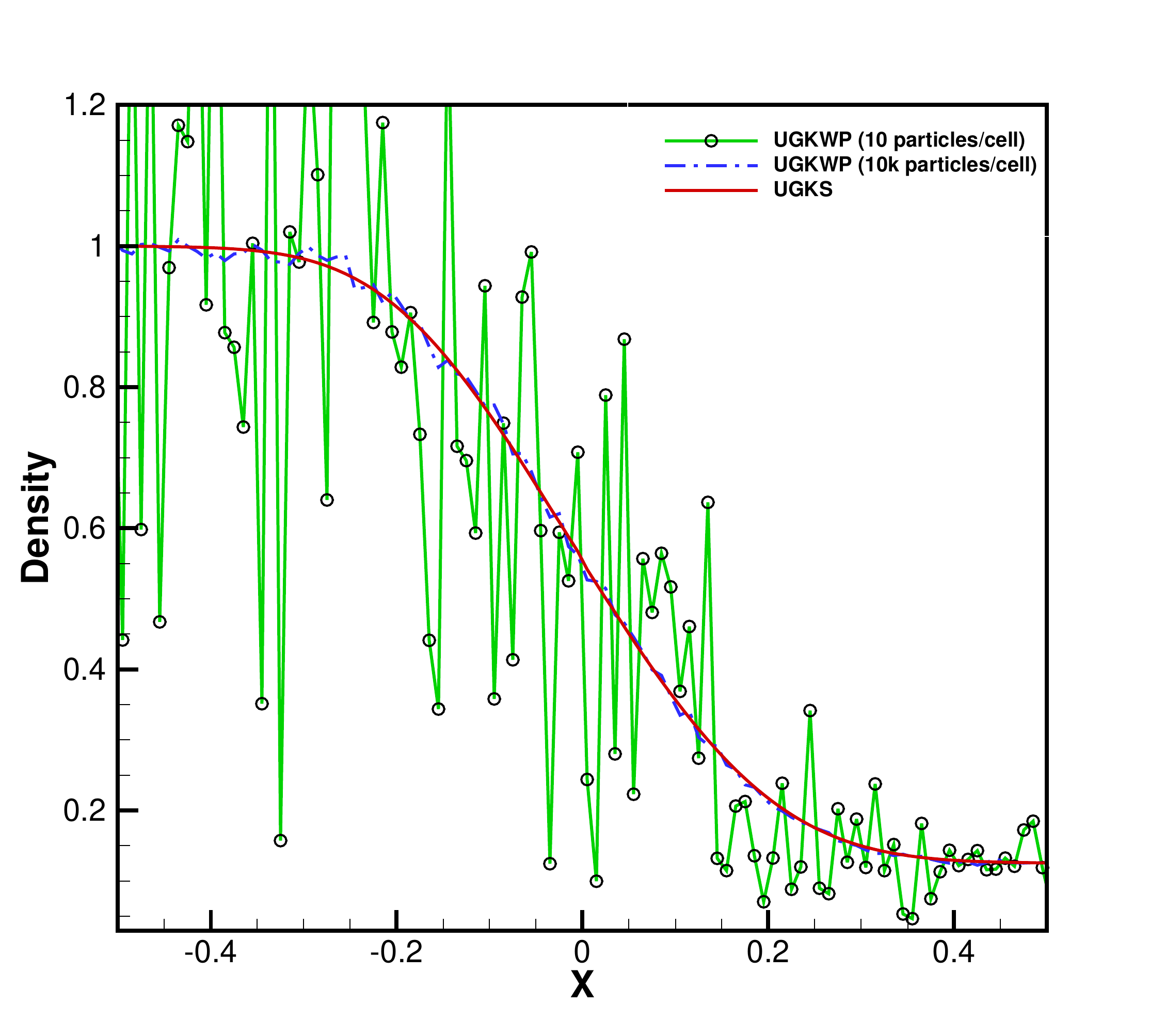}{a}
\includegraphics[width=0.48\textwidth]{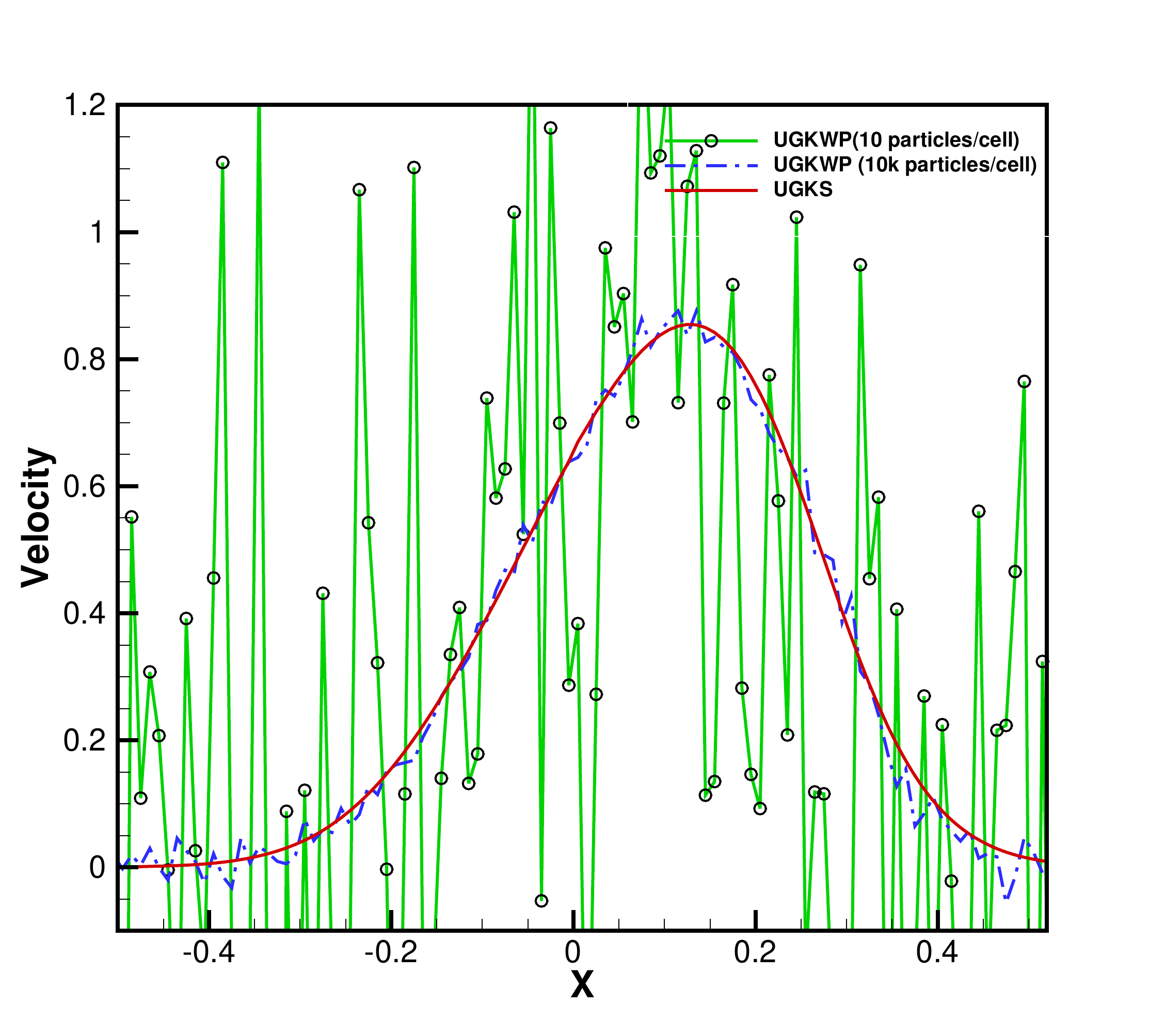}{b}
\includegraphics[width=0.48\textwidth]{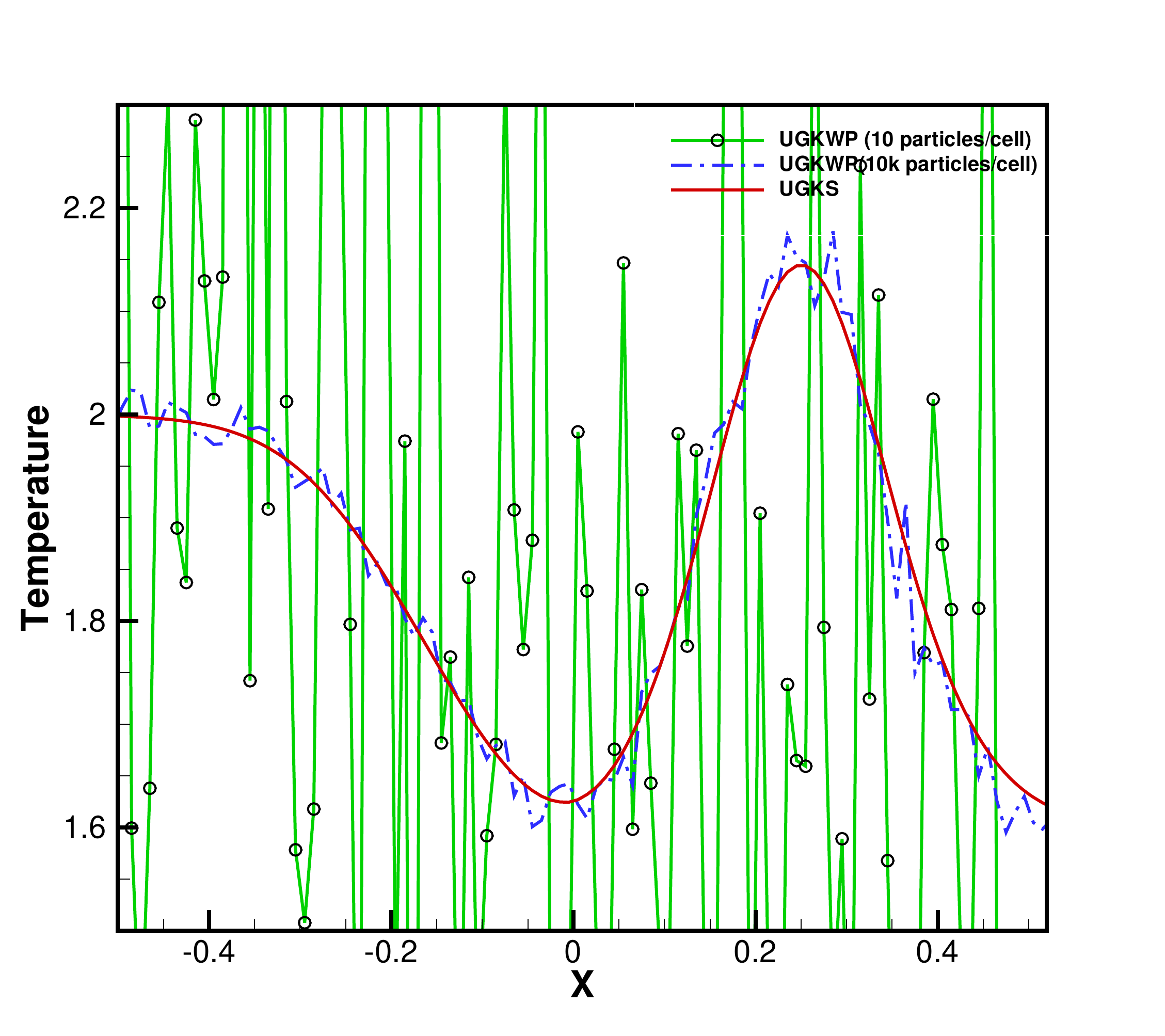}{c}
\includegraphics[width=0.48\textwidth]{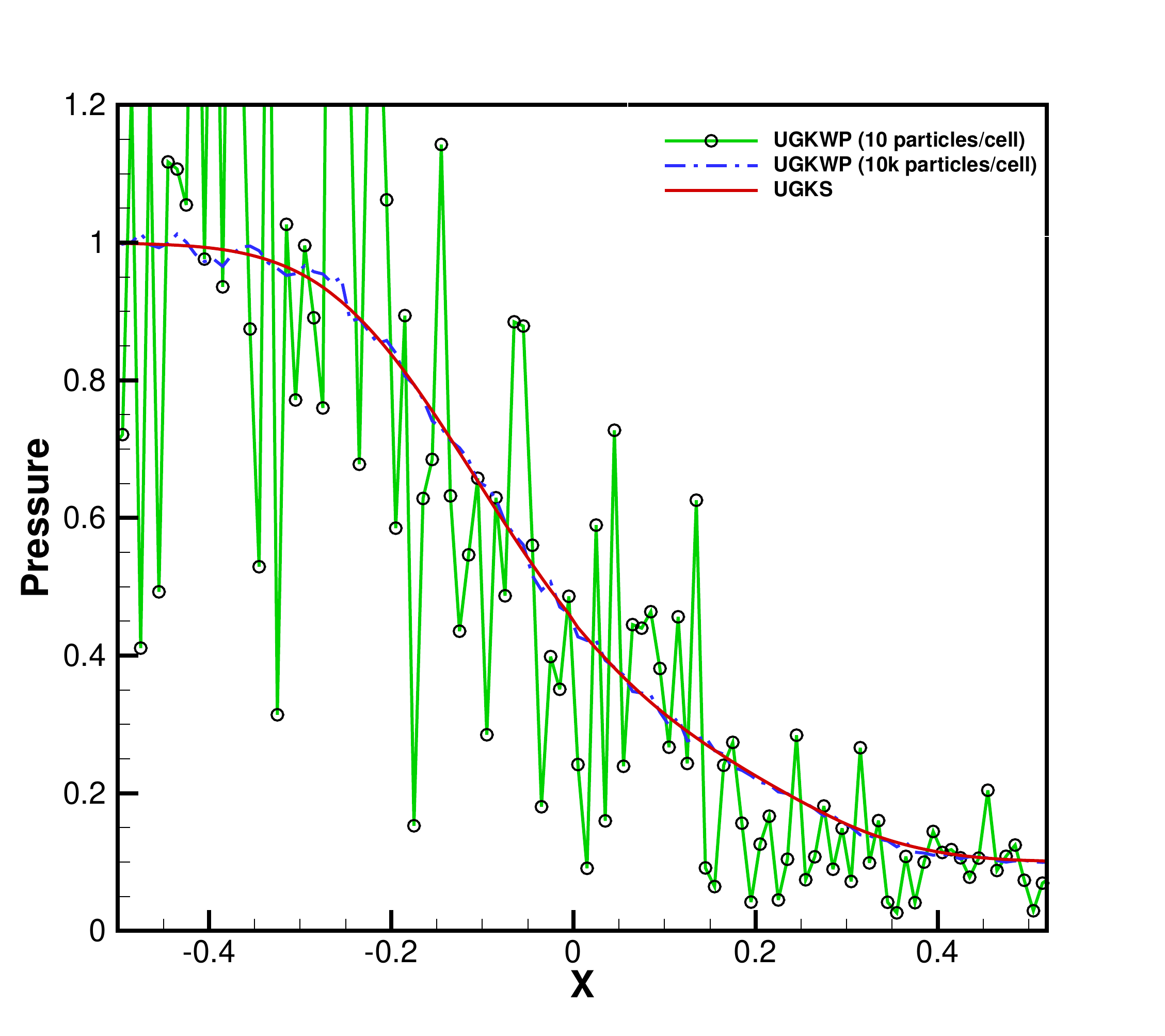}{d}
\caption{(a) Density, (b) velocity, (c) temperature, and (d) pressure profiles of Sod shock tube at $t=0.15$ with Knudsen number $\mathrm{Kn}=10^{-1}$. Symbol and green line are the UGKWP method result with 10 particles per cell; blue dashed line is the UGKWP method result with $10^{4}$ particles per cell; and red solid line is the UGKS solution.}
\label{sod1}
\end{figure}

\begin{figure}
\centering
\includegraphics[width=0.48\textwidth]{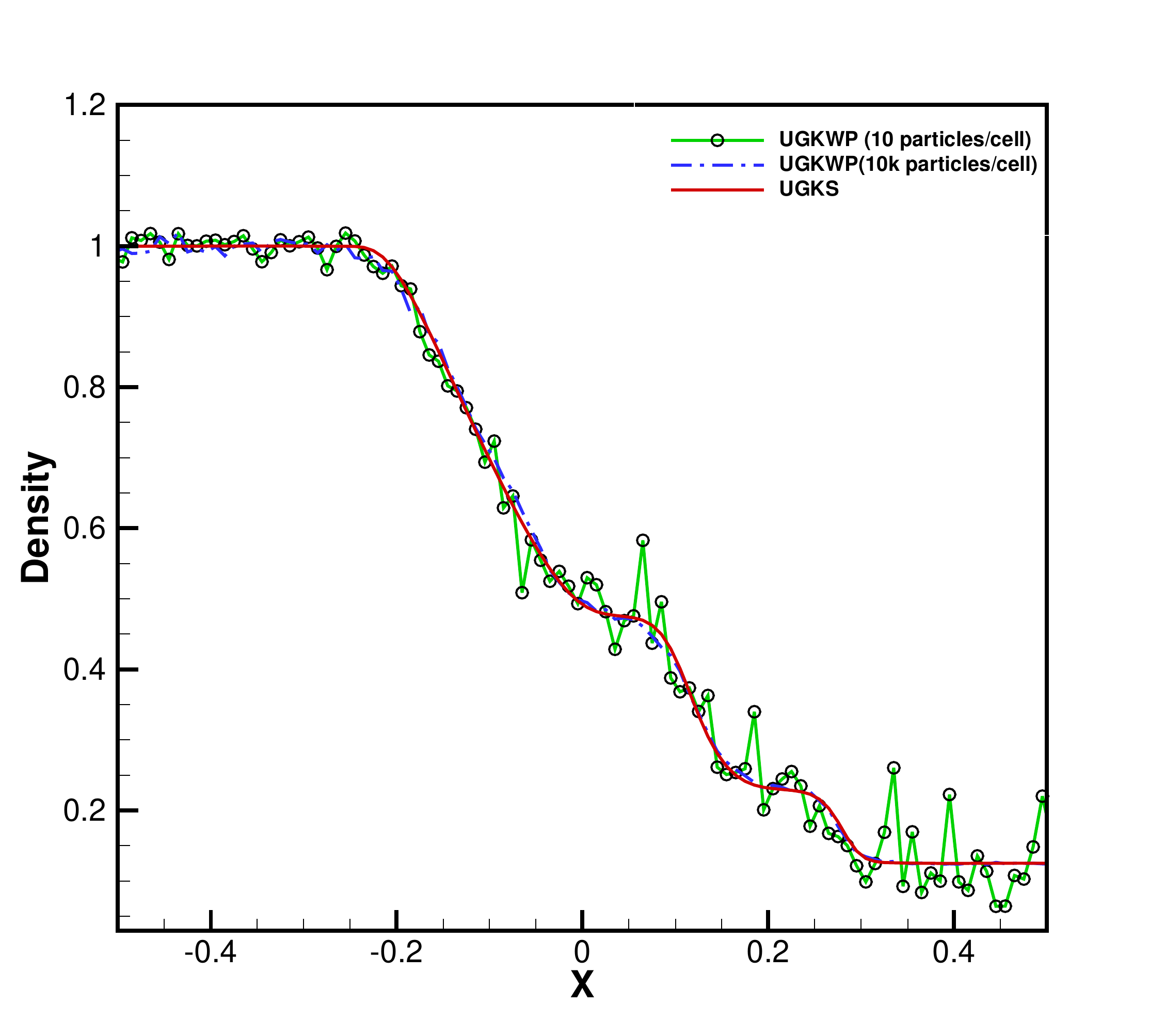}{a}
\includegraphics[width=0.48\textwidth]{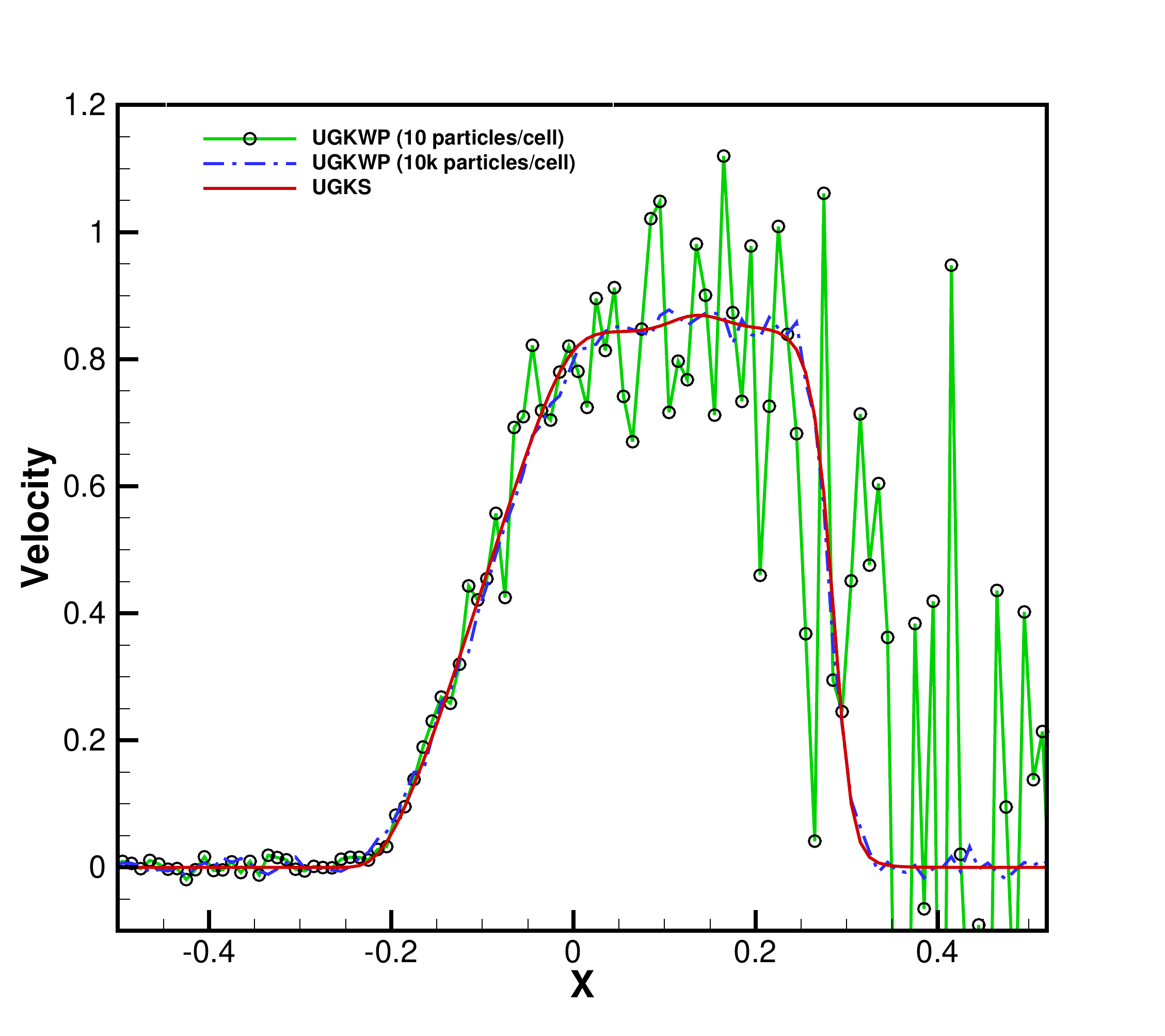}{b}
\includegraphics[width=0.48\textwidth]{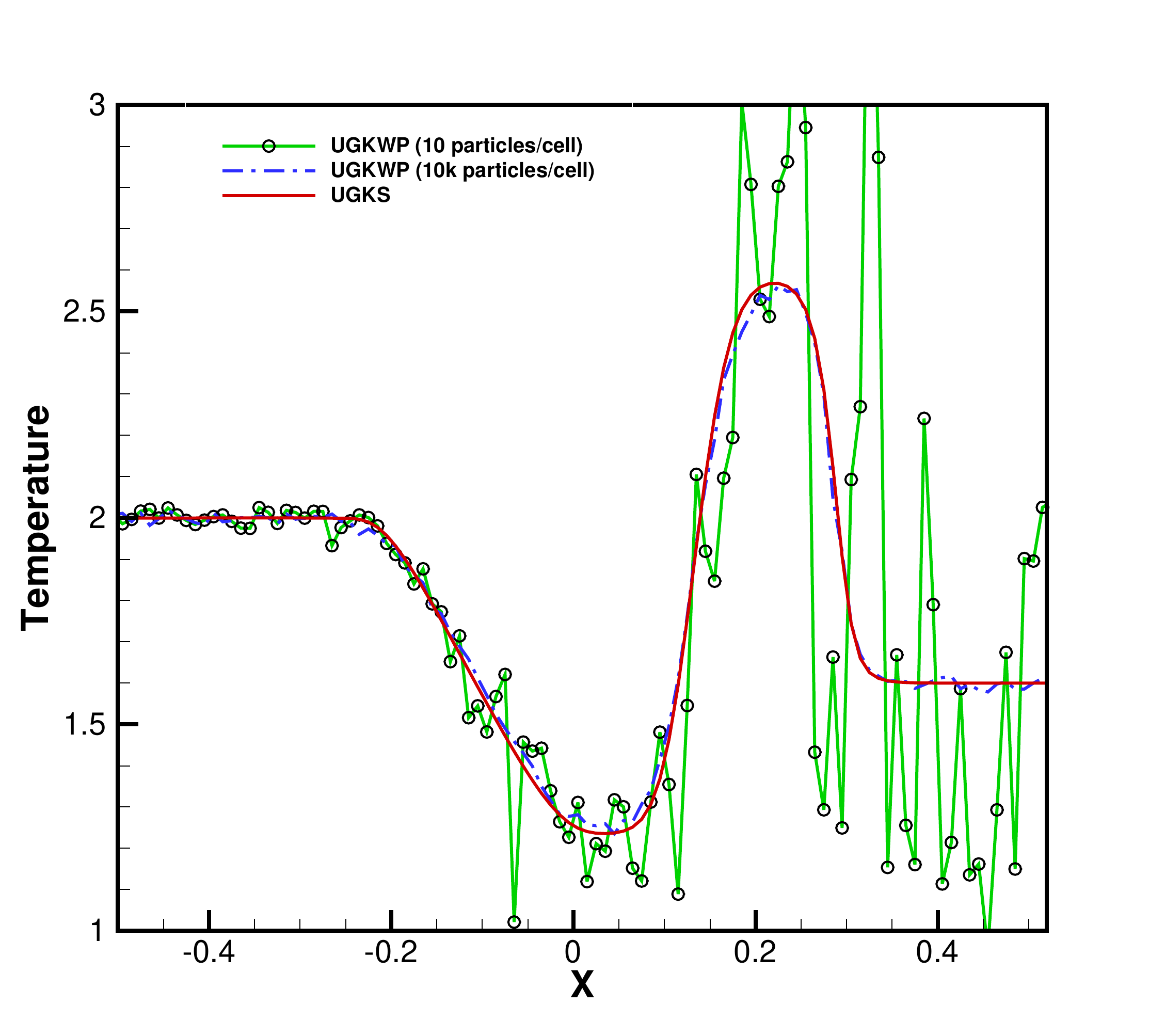}{c}
\includegraphics[width=0.48\textwidth]{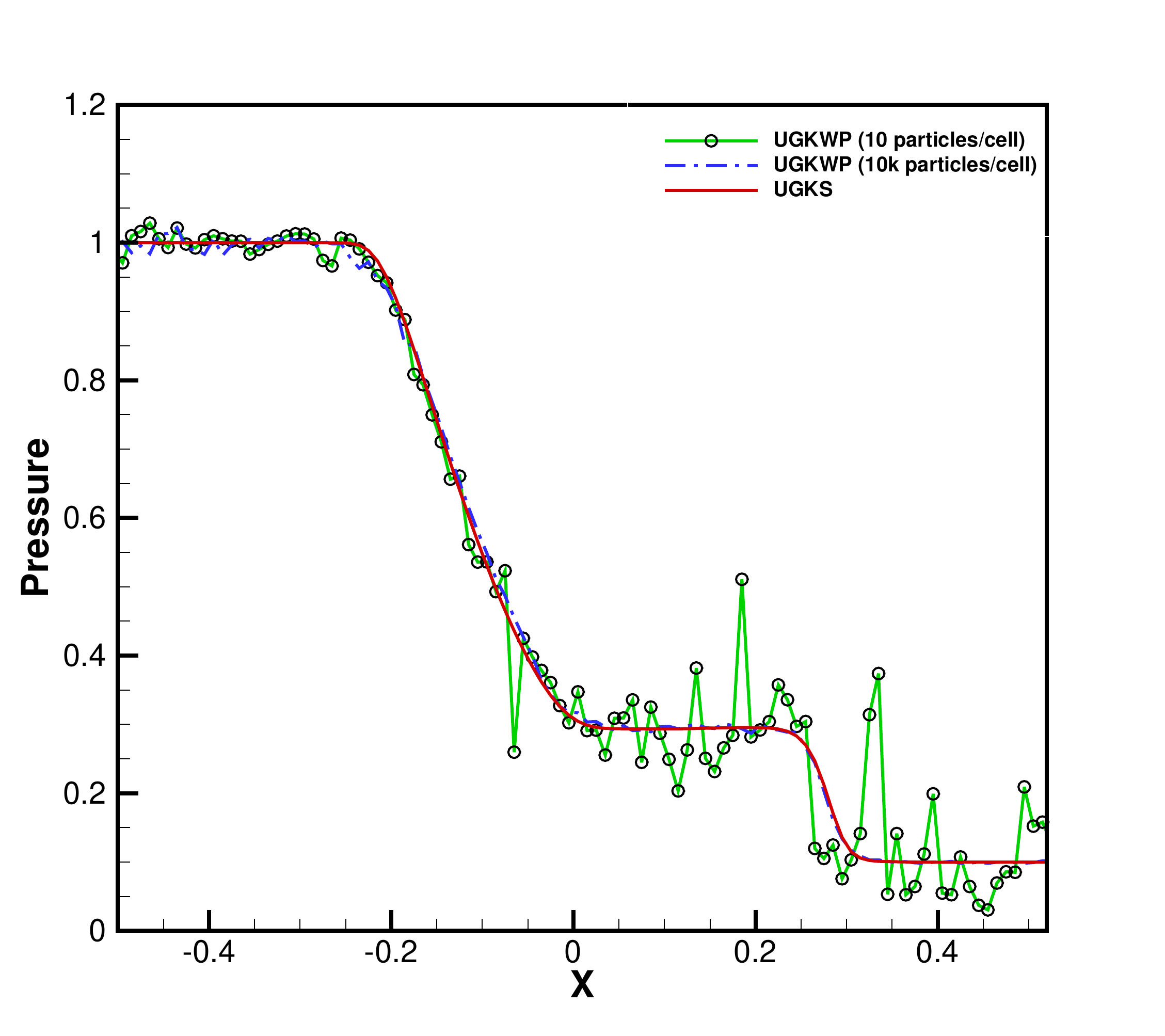}{d}
\caption{(a) Density, (b) velocity, (b) temperature, and (d) pressure profiles of Sod test at $t=0.15$ with Knudsen number $\mathrm{Kn}=10^{-3}$. Symbol and green line are the UGKWP method result with 10 particles per cell; blue dashed line is the UGKWP method result with $10^{4}$ particles per cell; and red solid line is the UGKS solution.}
\label{sod2}
\end{figure}

\begin{figure}
\centering
\includegraphics[width=0.48\textwidth]{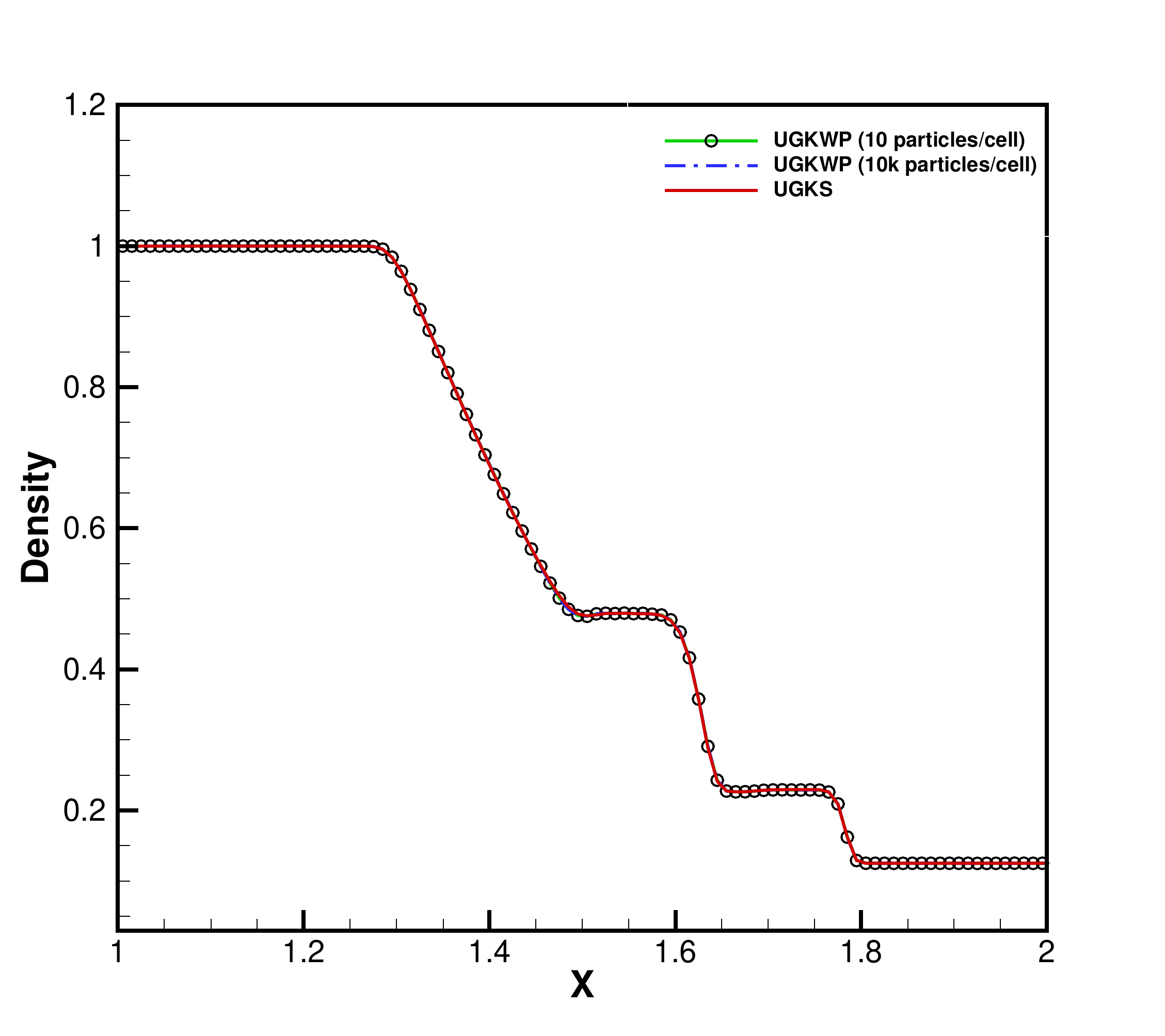}{a}
\includegraphics[width=0.48\textwidth]{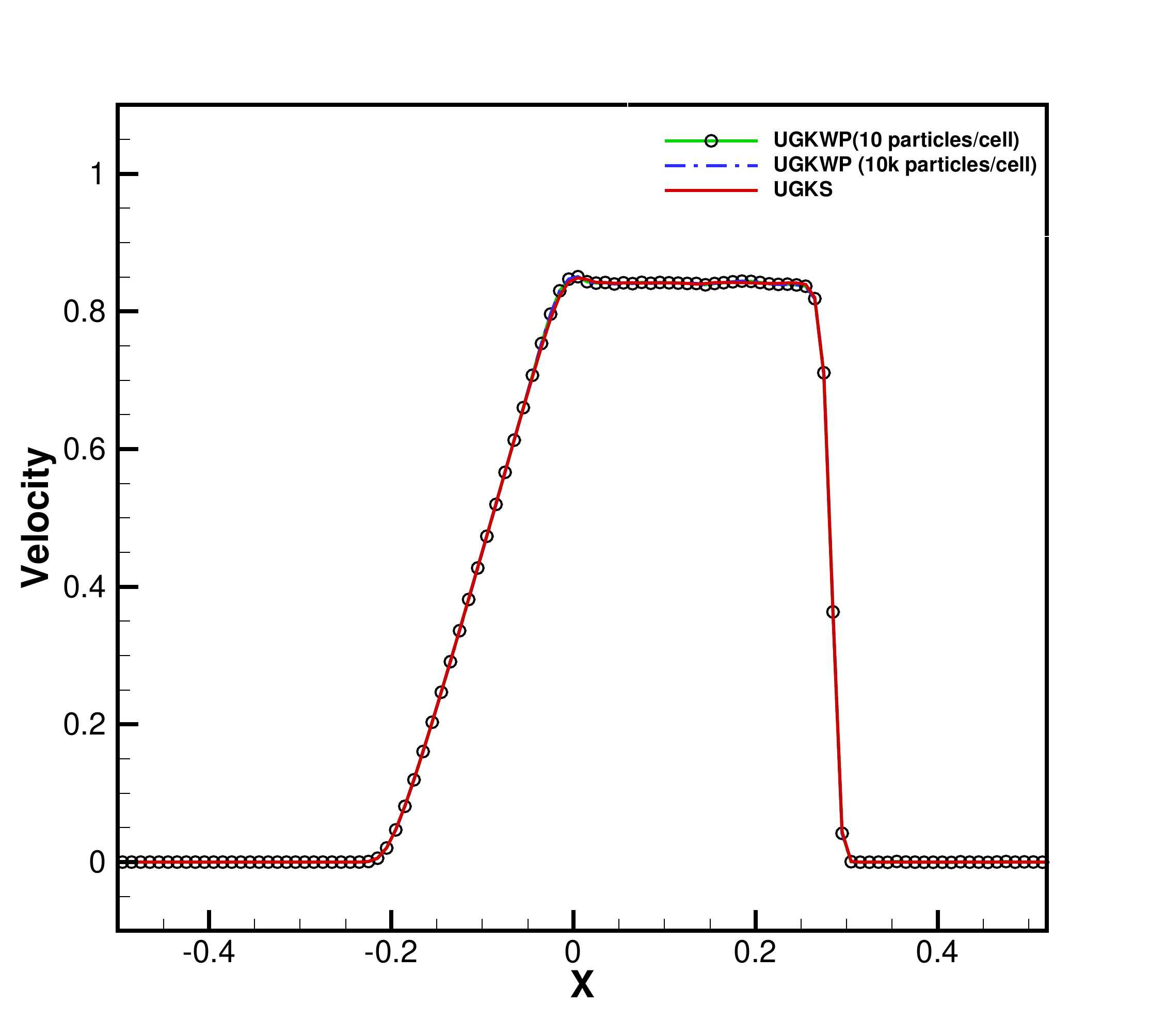}{b}
\includegraphics[width=0.48\textwidth]{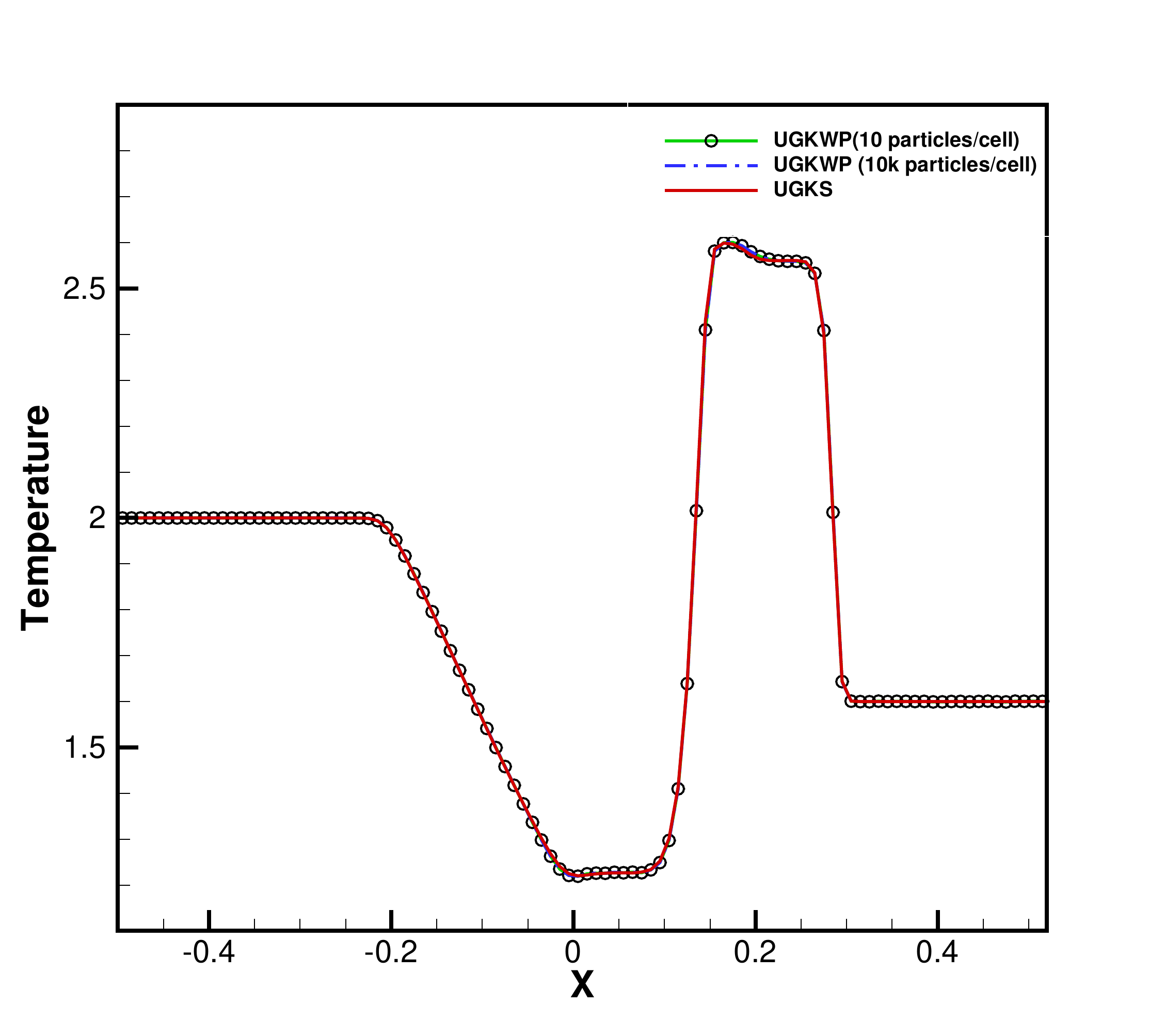}{c}
\includegraphics[width=0.48\textwidth]{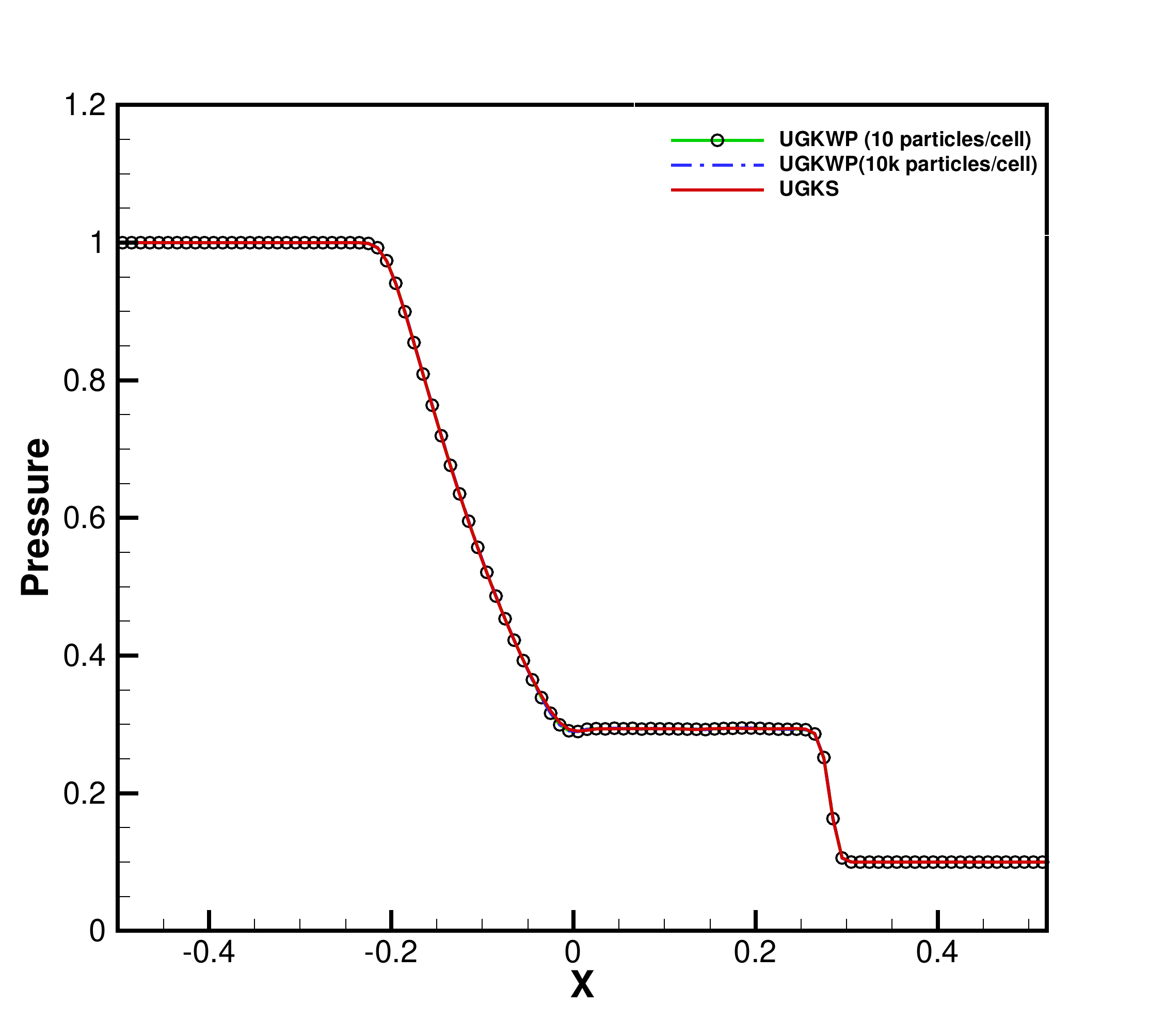}{d}
\caption{(a) Density, (b) velocity, (b) temperature, and (d) pressure profiles of Sod test at $t=0.15$ with Knudsen number $\mathrm{Kn}=10^{-5}$. Symbol and green line are the UGKWP method result with 10 particles per cell; blue dashed line is the UGKWP method result with $10^{4}$ particles per cell; and red solid line is the UGKS solution.}
\label{sod3}
\end{figure}

The comparison between the UGKWP method and UGKS solution at $t=0.15$ is shown in Fig. \ref{sod1}-\ref{sod3}.
The solution in rarefied regime with $\mathrm{Kn}=0.1$ is shown in Fig. \ref{sod1}, the solution of UGKWP method with $10^{4}$ number particle per cell agrees with the UGKS solution, and extreme large statistical noise can be observed when we reduce the number of particle to 10 per cell. The results show that the UGKWP method can capture the rarefied solution.
However, a sufficient number of simulation particles is required or averaging process need to be carried out in order to get sufficiently accurate solution.
Next, we reduce the Knudsen number to $10^{-3}$. As shown in Fig. \ref{sod2}, the solution with large number of simulation particles ($10^{4}$ per cell) agree well with UGKS solution, and the noise of the solution with small number of simulation particles (10 particles per cell) is smaller than the rarefied case with $\mathrm{Kn}=0.1$, especially in the upstream.
When we move to $\mathrm{Kn}=10^{-5}$, even with 10 simulation particles per cell, the UGKWP method well agrees with the UGKS solution. The Sod shows that the UGKWP method can capture the flow physical in different flow regimes, and in the near continuum regime the statistical noise of UGKWP method is much lower than conventional particle methods.

\subsection{Normal shock}
\begin{figure}
\centering
\includegraphics[width=0.7\textwidth]{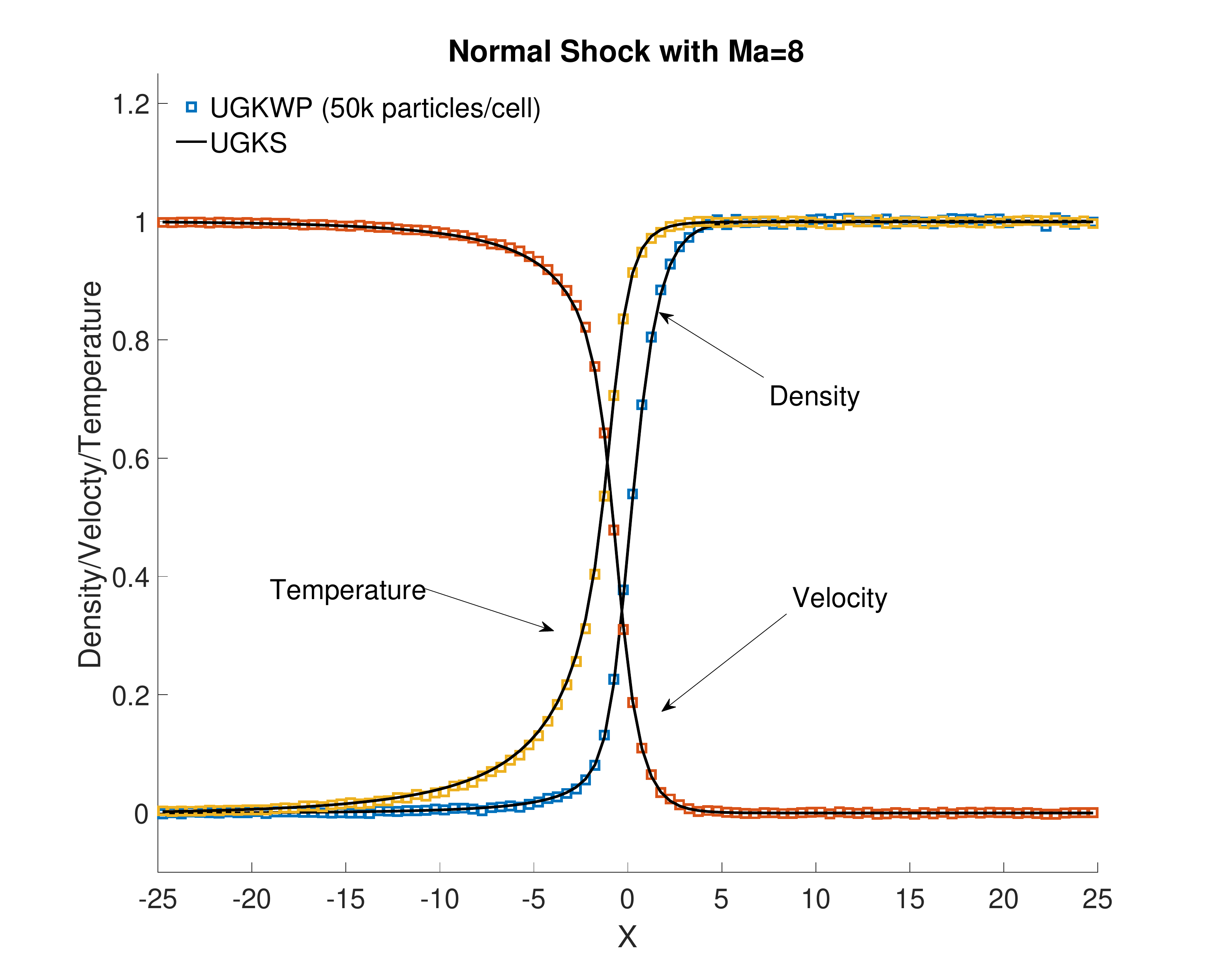}
\includegraphics[width=0.7\textwidth]{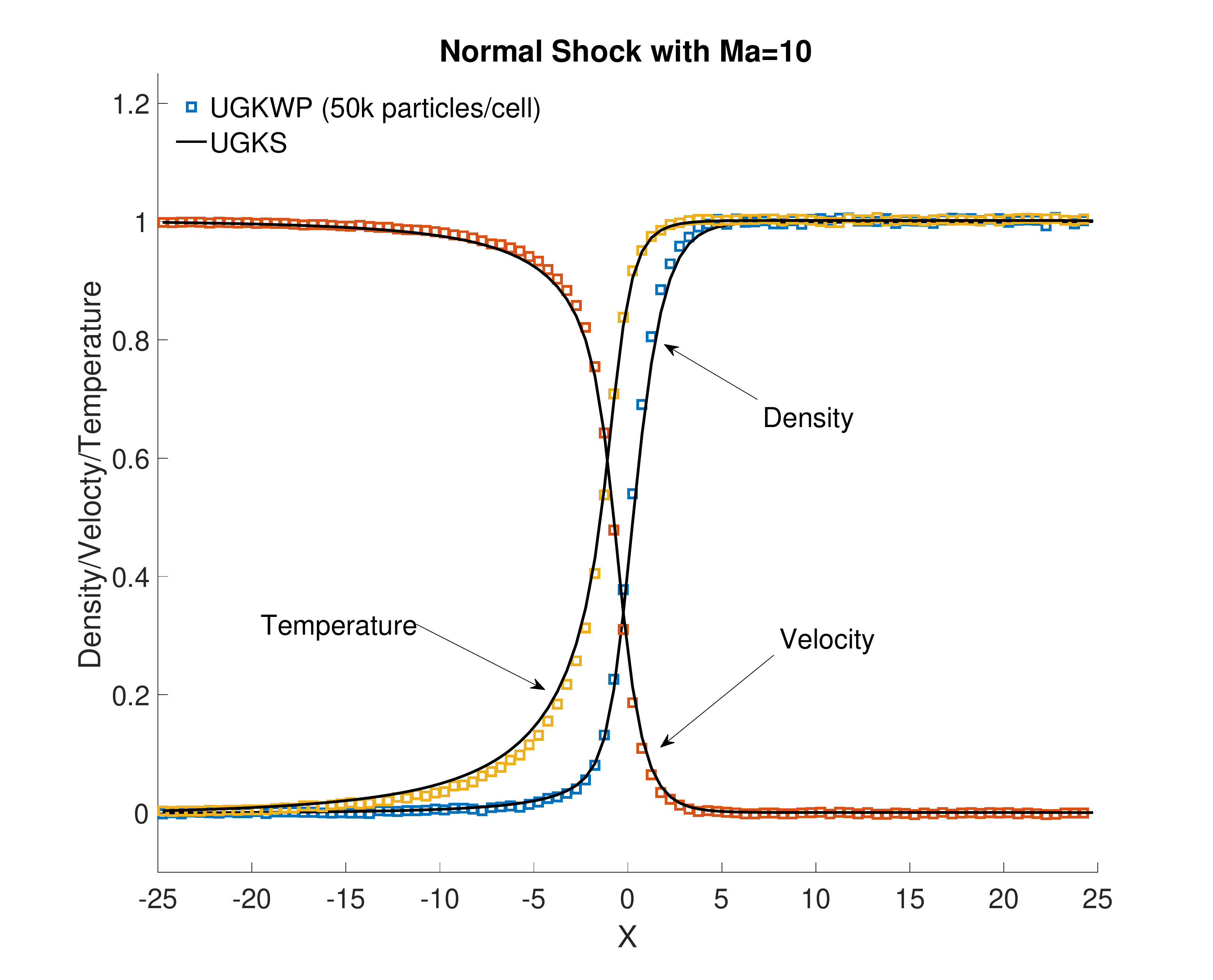}
\caption{Normalized density, velocity, and temperature profile of normal shock wave at $\mathrm{M}=8$ (top) and $\mathrm{M}=10$ (bottom). The UGKWP method solution is shown in symbol, and the UGKS solution is shown in line.}
\label{shock}
\end{figure}

To demonstrate the capability of UGKWP method in capturing the highly non-equilibrium flow, the one dimensional shock wave is studied. The gaseous medium is argon, the viscous coefficient of which follows Eq.\eqref{vhs-vs1}-\eqref{vhs-vs2} with the temperature dependency index $\omega=0.81$.
In this calculation, the reference length is the upstream mean free path, and the computational domain is [-25,25] with 100 cells.
The upstream ($x\le0$) and downstream ($x>0$) is connected by the Rankine–-Hugoniot condition.
We calculate two test cases with upstream Mach number $\mathrm{M}=8$ and $\mathrm{M}=10$.
In order to reduce the statistical noise, $5\times10^{4}$ number of simulation particles are used in each cell.
As shown in Fig. \ref{shock}, the normalized solution of UGKWP method agrees well with the UGKS solution, which proves the ability of UGKWP method in capturing the non-equilibrium flow.

\subsection{Cylinder flow}
In this section, we calculate the supersonic argon gas flow passing over a circular cylinder at different Mach and Knudsen numbers.
For argon gas, the molecular mass is $m_0=6.63\times 10^{-26}$ kg; the molecular diameter $d=4.17\times10^{-10}$ m; and the specific heat ratio $\gamma=5/3$. The variable hard sphere (VHS) model is used to model the molecular interaction, and the viscosity follows
\begin{equation}
  \mu=\frac{15\sqrt{\pi m k T_\infty}}{2\pi d^2(5-2\omega)(7-2\omega)}\left(\frac{T}{T_\infty}\right)^\omega,
\end{equation}
with the temperature dependency index $\omega=0.81$.
As shown in table \ref{cylinder-initial}, the incoming flow Mach number is chosen to be $5, 20, 30$, and the Knudsen number with respect to the cylinder radius is set as $1.0, 0.1, 10^{-4}$. The dimensionless quantities are used with respect to the reference length as the cylinder radius $L_{ref}=R$, the reference velocity $U_{ref}=\sqrt{2RT_\infty}$, the reference time $t_{ref}=L_{ref}/U_{ref}$, the reference density $\rho_{ref}= \rho_\infty$, and the reference temperature $T_{ref}=T_\infty$.

\begin{table}[ht]
\caption{Incoming flow condition for cylinder flow}
\centering
\begin{tabular}{|c|c|c|c|c|c|c|}
\hline
$\mathrm{Kn}_{\infty}$ & $n_{\infty}\mathrm{[particles/m^3]}$ & $\rho_{\infty}\mathrm{[Kg/m^3]}$ &
$T_{\infty}\mathrm{[K]}$ & $T_{w}\mathrm{[K]}$ & $U_{\infty}\mathrm{[m/s]}$ & $\mathrm{M}_{\infty}$ \\ \hline
\multirow{3}{*}{$10^{-4}$} & \multirow{3}{*}{$1.294\times10^{24}$} & \multirow{3}{*}{$8.585\times10^{-2}$} & \multirow{3}{*}{273} & \multirow{3}{*}{273} & 1538.794 & 5 \\ \cline{6-7}
 &  &  &  &  & 6155.17 & 20 \\ \cline{6-7}
 &  &  &  &  & 9232.76 & 30 \\ \hline
\multirow{3}{*}{0.1} & \multirow{3}{*}{$1.294\times10^{21}$} & \multirow{3}{*}{$8.585\times10^{-5}$} & \multirow{3}{*}{273} & \multirow{3}{*}{273} & 1538.794 & 5 \\ \cline{6-7}
 &  &  &  &  & 6155.17 & 20 \\ \cline{6-7}
 &  &  &  &  & 9232.76 & 30 \\ \hline
\multirow{3}{*}{1} & \multirow{3}{*}{$1.294\times10^{20}$} & \multirow{3}{*}{$8.585\times10^{-6}$} &
\multirow{3}{*}{273} & \multirow{3}{*}{273} & 1538.794 & 5 \\ \cline{6-7}
 &  &  &  &  & 6155.17 & 20 \\ \cline{6-7}
 &  &  &  &  & 9232.76 & 30 \\ \hline
\end{tabular}
\label{cylinder-initial}
\end{table}

\begin{figure}
\centering
\includegraphics[width=0.48\textwidth]{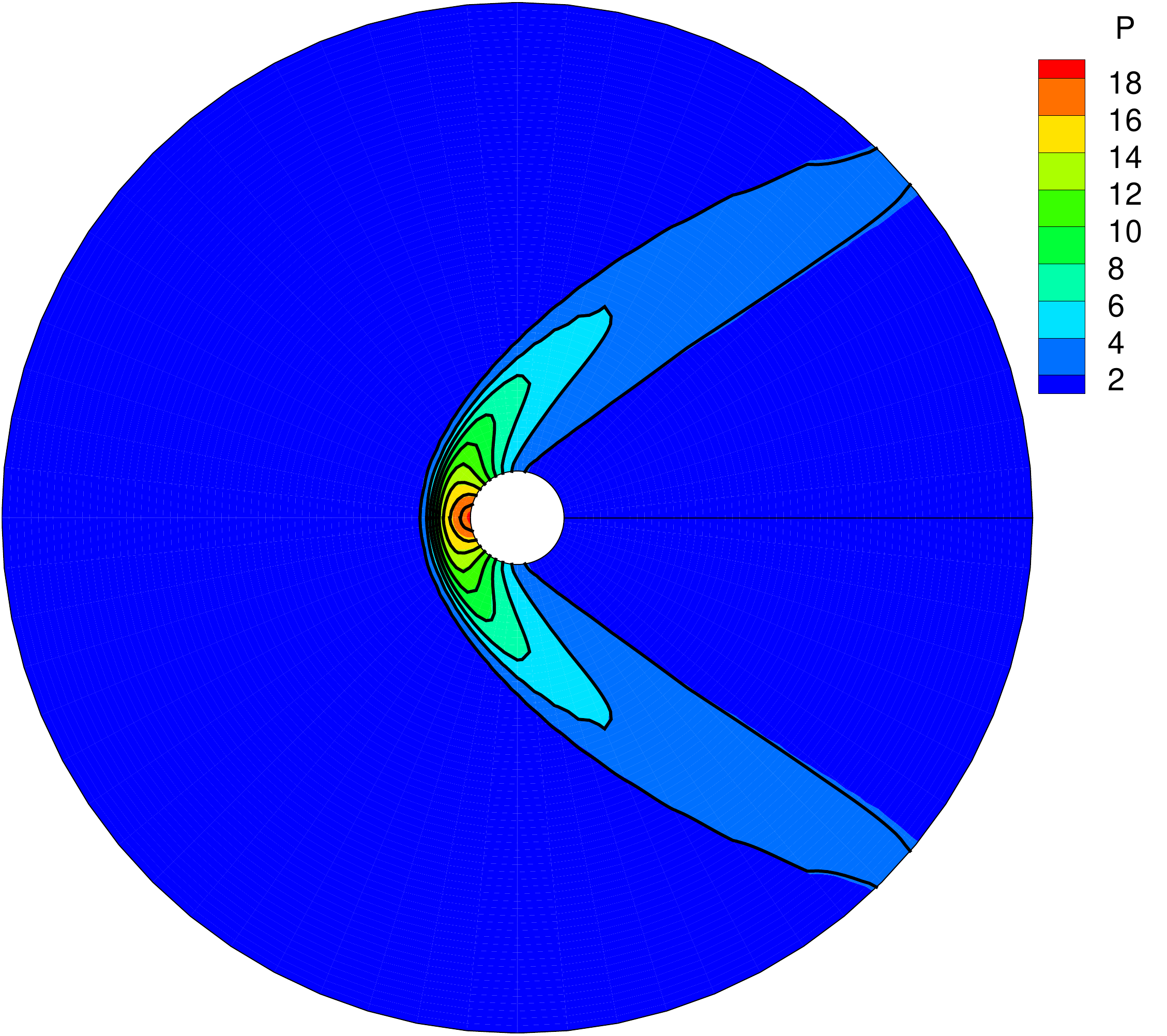}{a}
\includegraphics[width=0.48\textwidth]{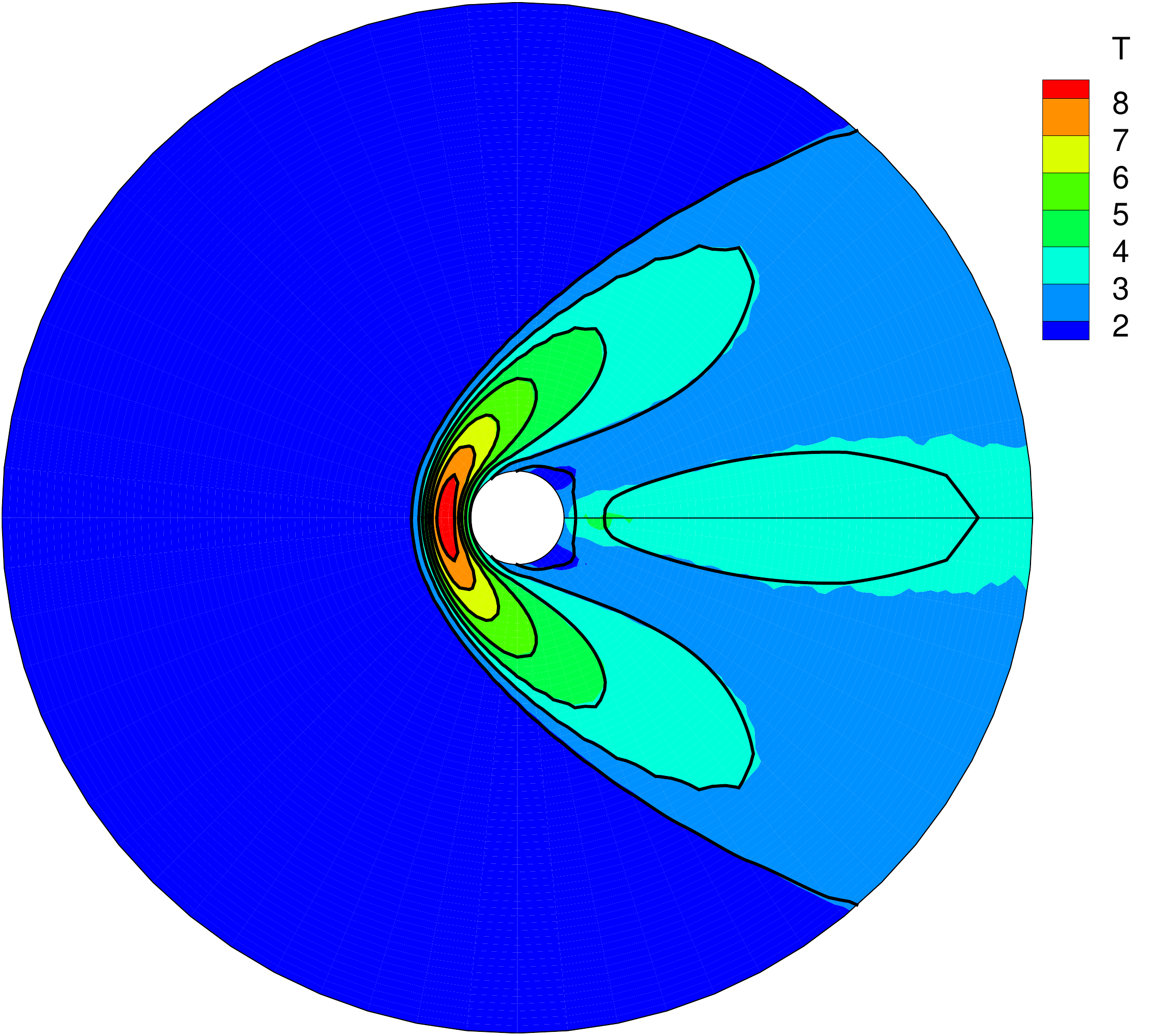}{b}
\includegraphics[width=0.48\textwidth]{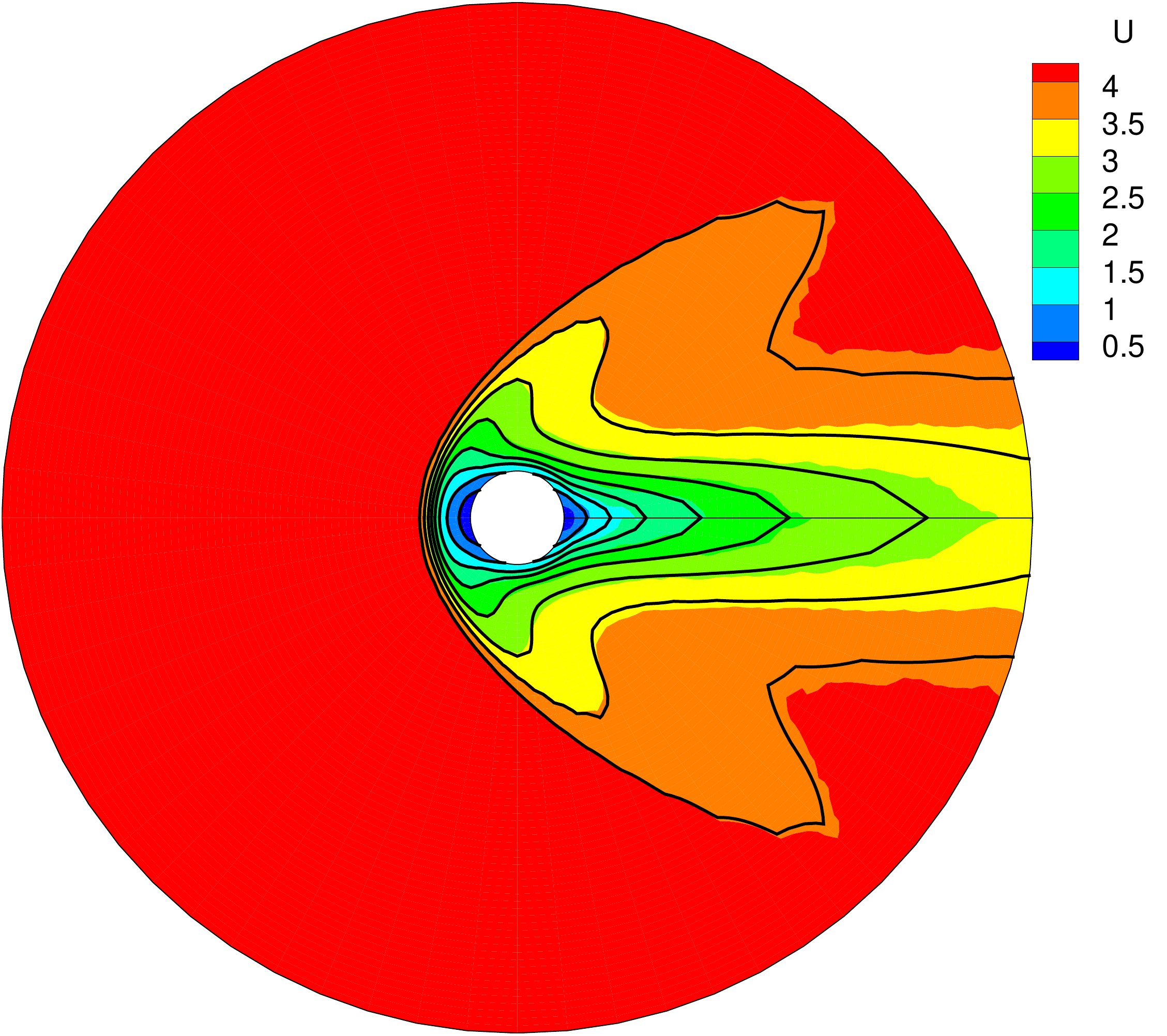}{c}
\includegraphics[width=0.48\textwidth]{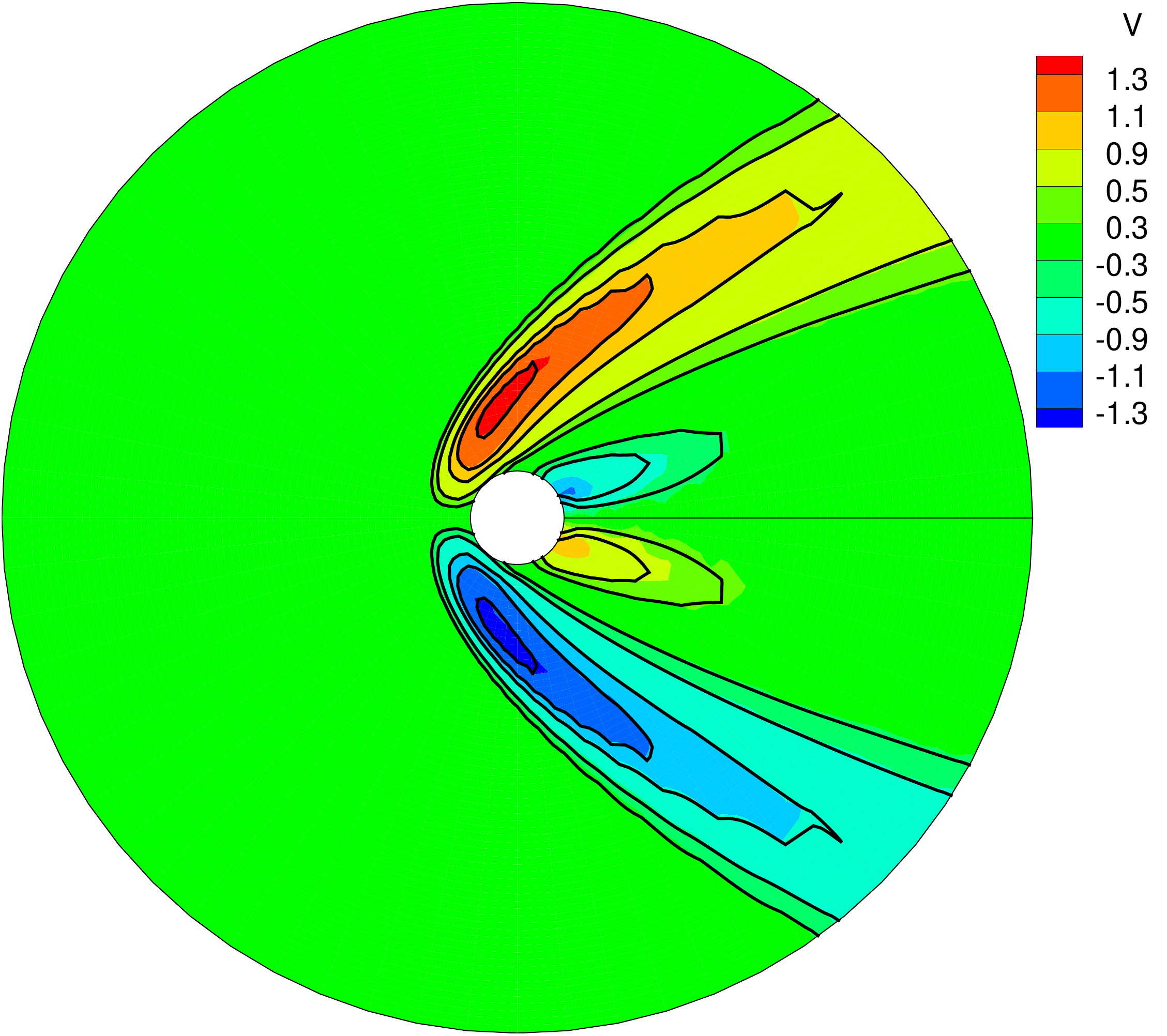}{d}
\caption{(a) Pressure, (b) temperature, (c) x directional velocity, and (d) y directional velocity contour for $\mathrm{M}=5$ and $\mathrm{Kn}=0.1$. The UGKWP method solution is shown in flood, and the UGKS solution is shown in contour line.}
\label{cylinder11}
\end{figure}

\begin{figure}
\centering
\includegraphics[width=0.48\textwidth]{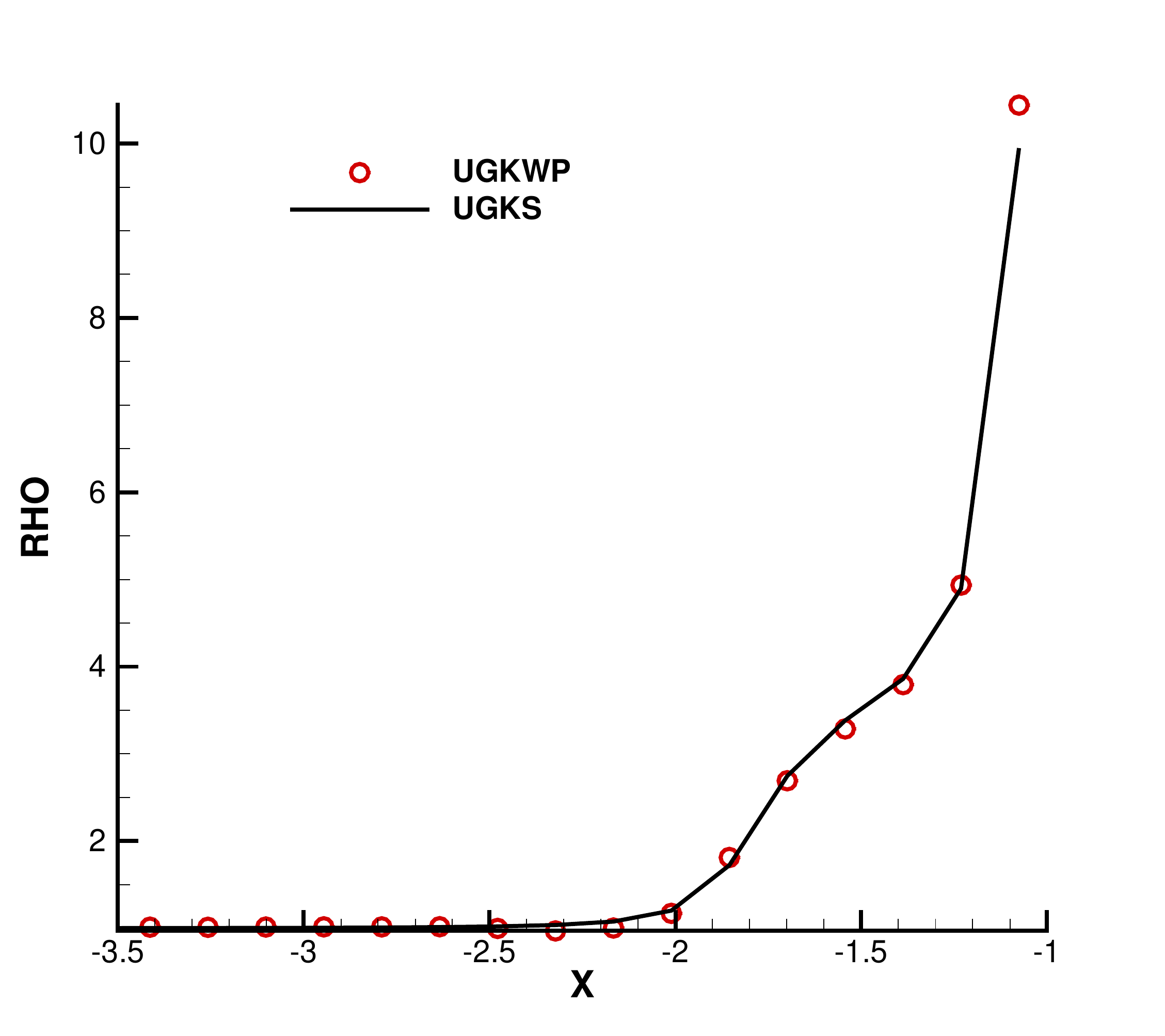}{a}
\includegraphics[width=0.48\textwidth]{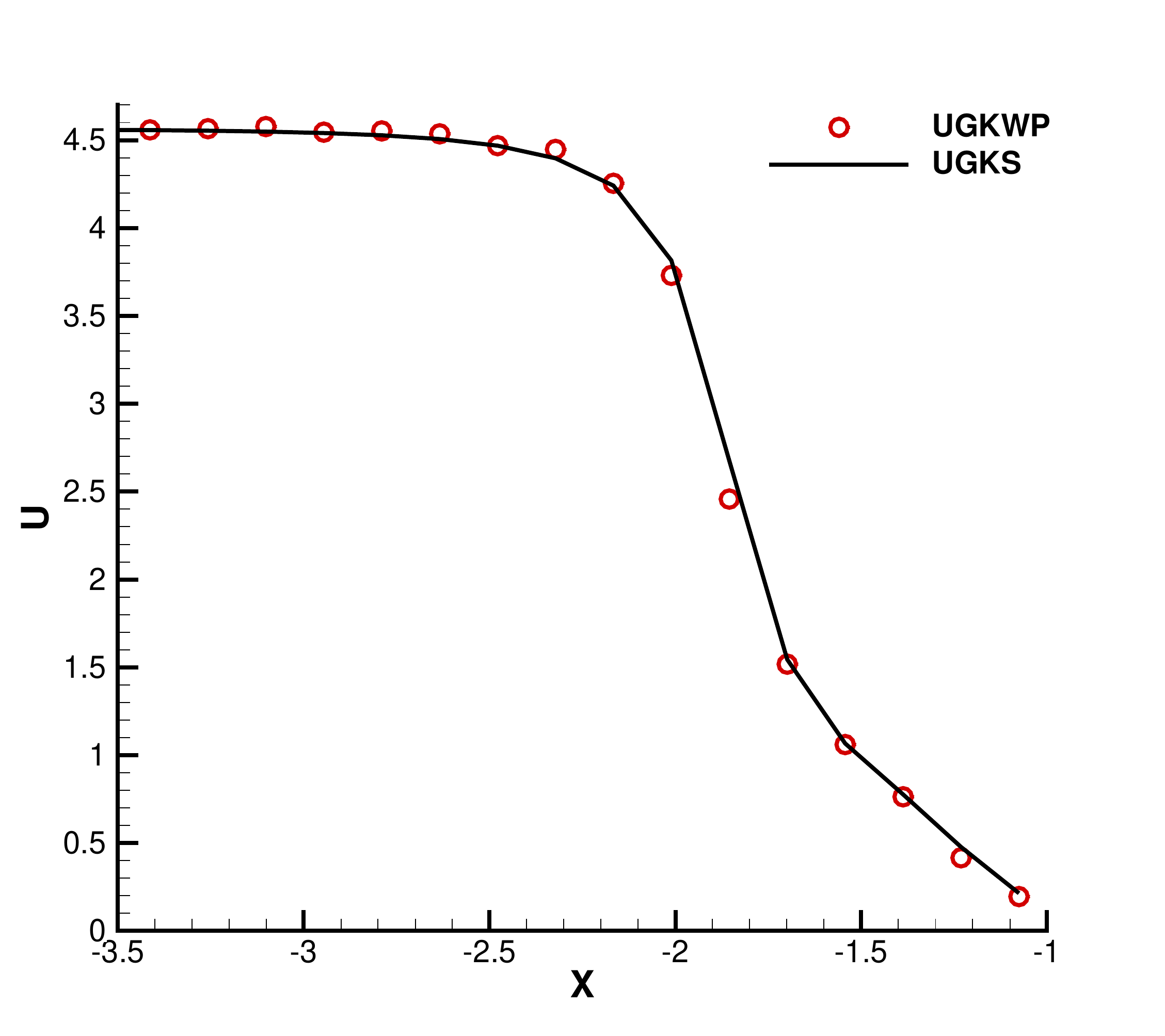}{b}
\includegraphics[width=0.48\textwidth]{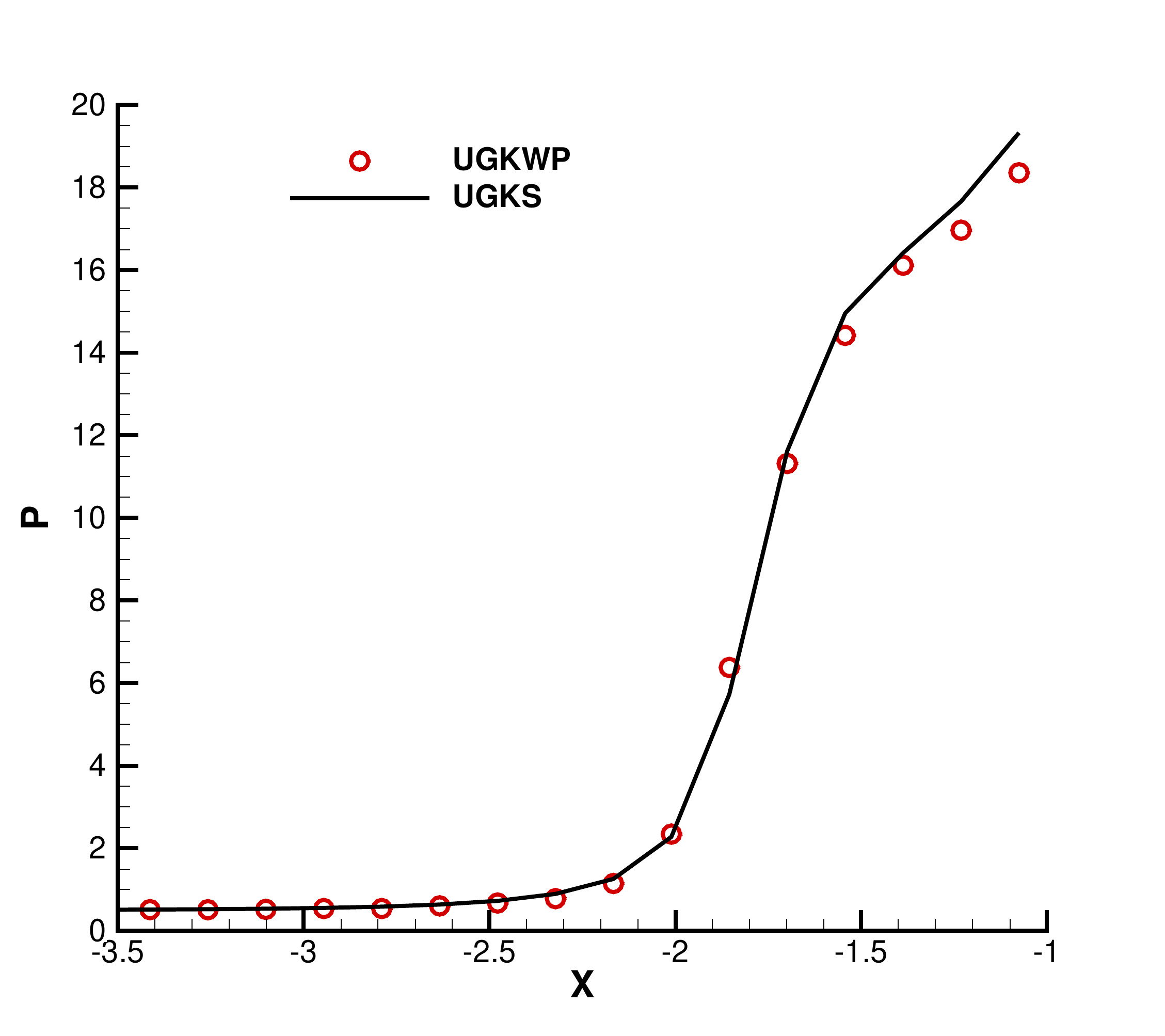}{c}
\includegraphics[width=0.48\textwidth]{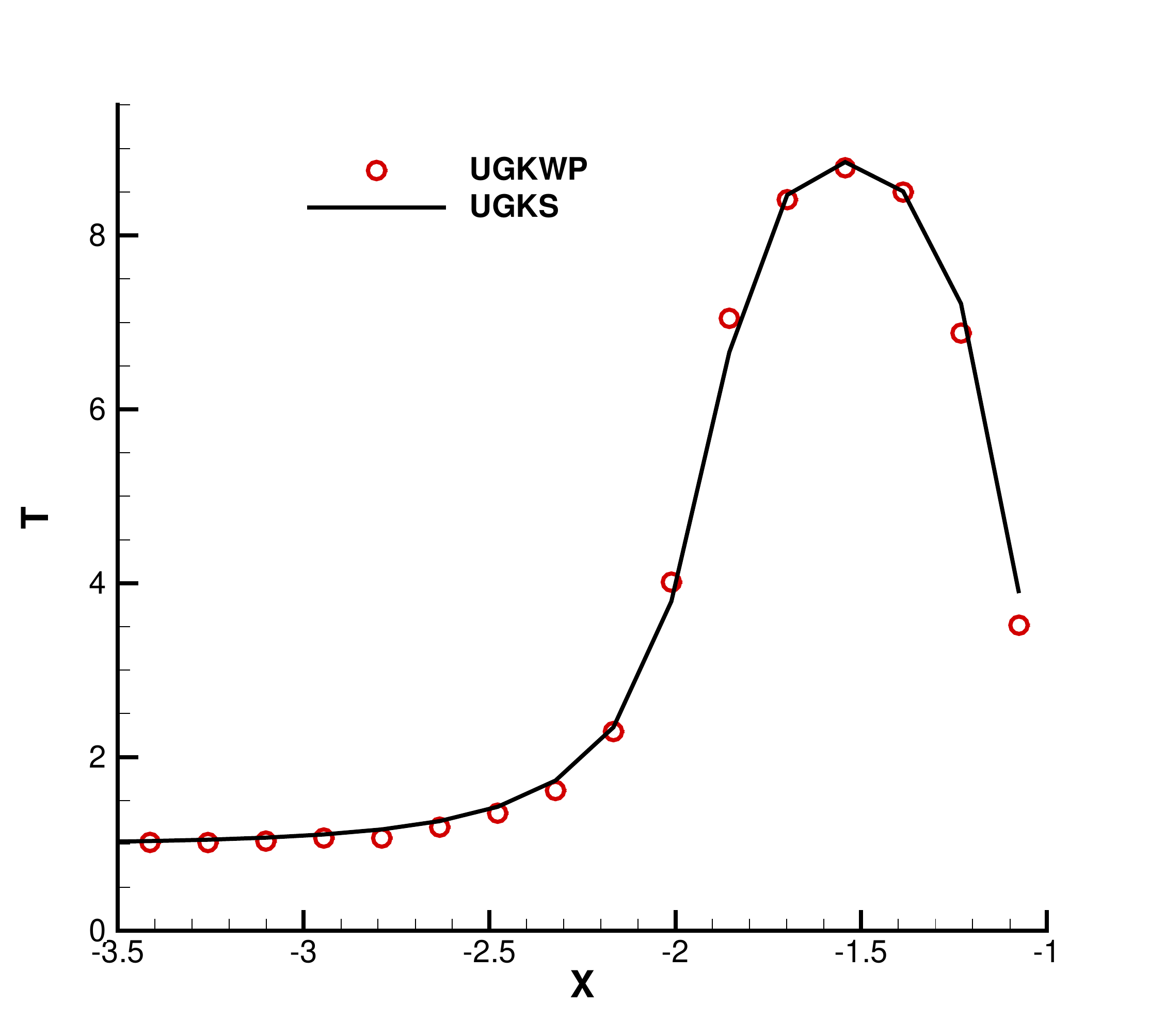}{d}
\caption{(a) Density, (b) x direction velocity, (c) pressure, (d) temperature profile along stagnation line for $\mathrm{M}=5$ and $\mathrm{Kn}=0.1$. The UGKWP method solution is shown in symbol, and the UGKS solution is shown in line.}
\label{cylinder12}
\end{figure}

\begin{figure}
\centering
\includegraphics[width=0.48\textwidth]{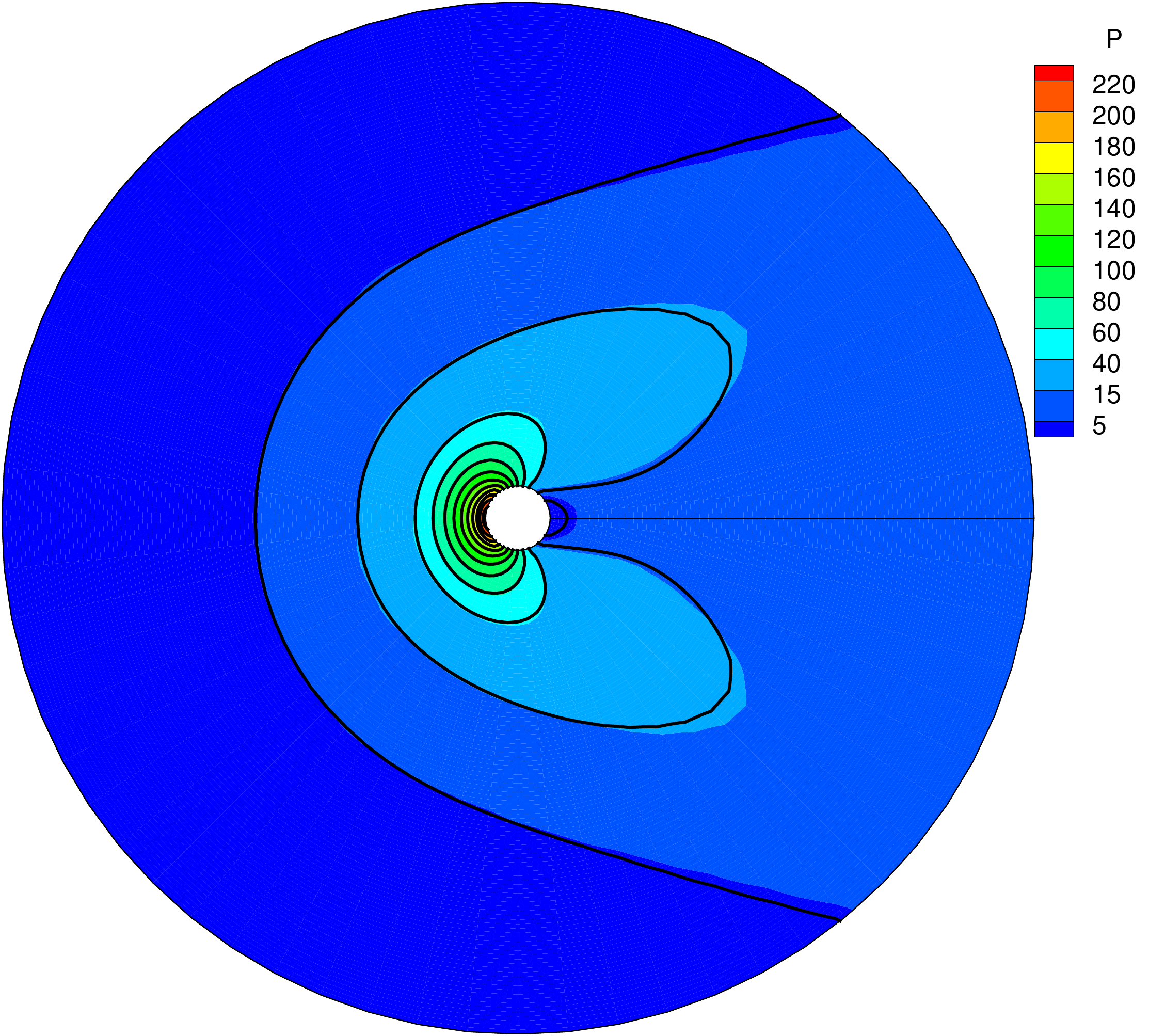}{a}
\includegraphics[width=0.48\textwidth]{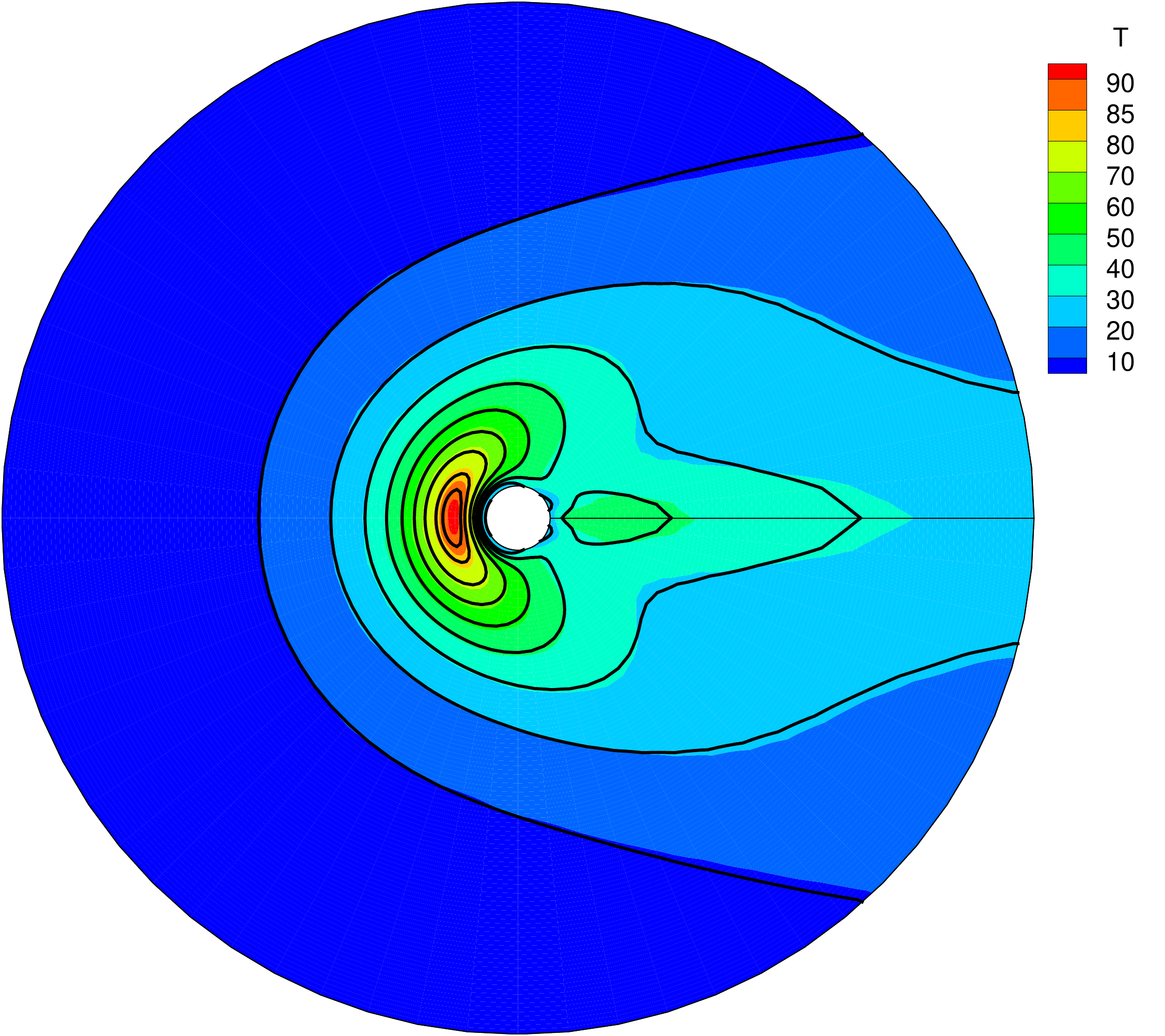}{b}
\includegraphics[width=0.48\textwidth]{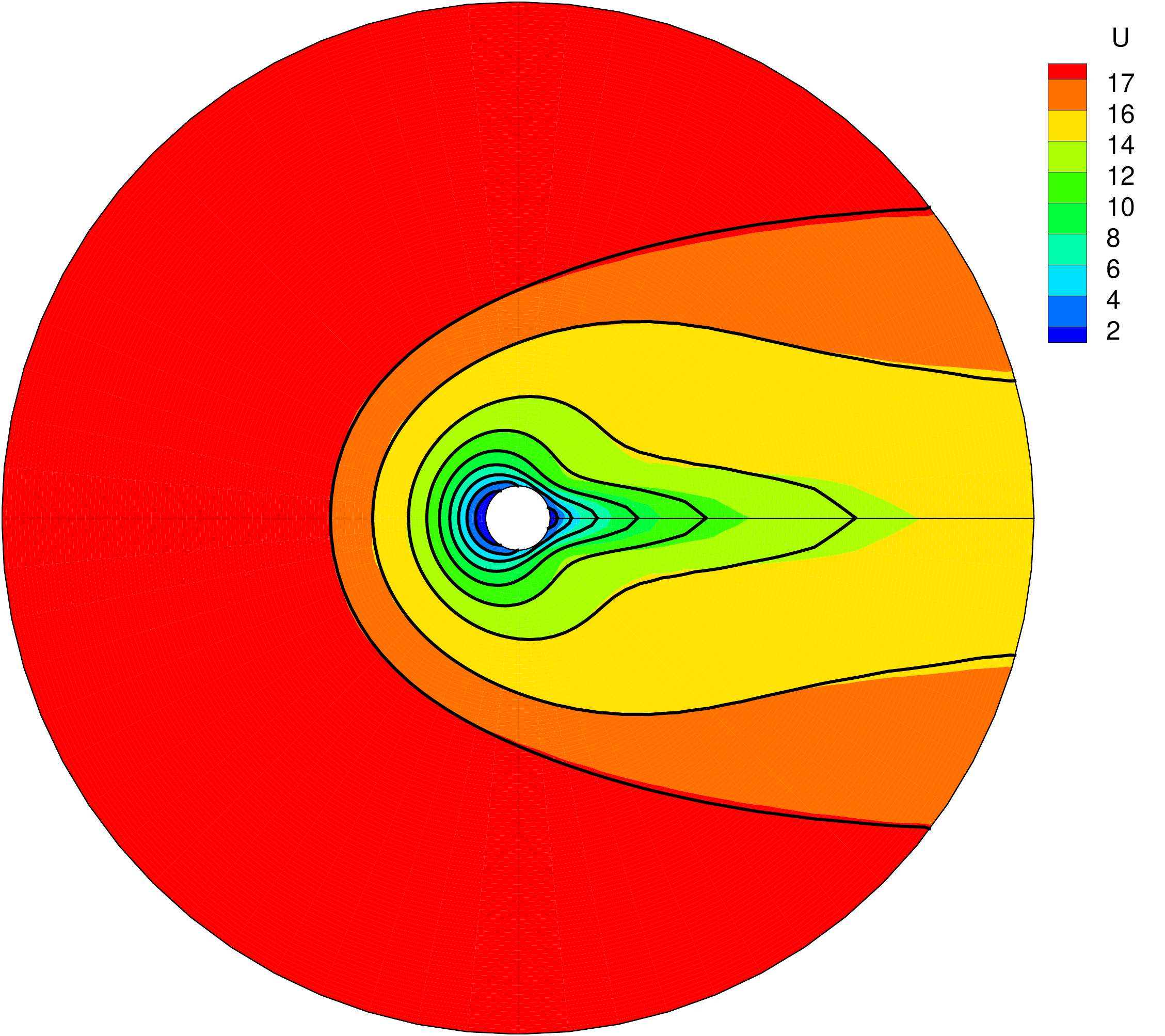}{c}
\includegraphics[width=0.48\textwidth]{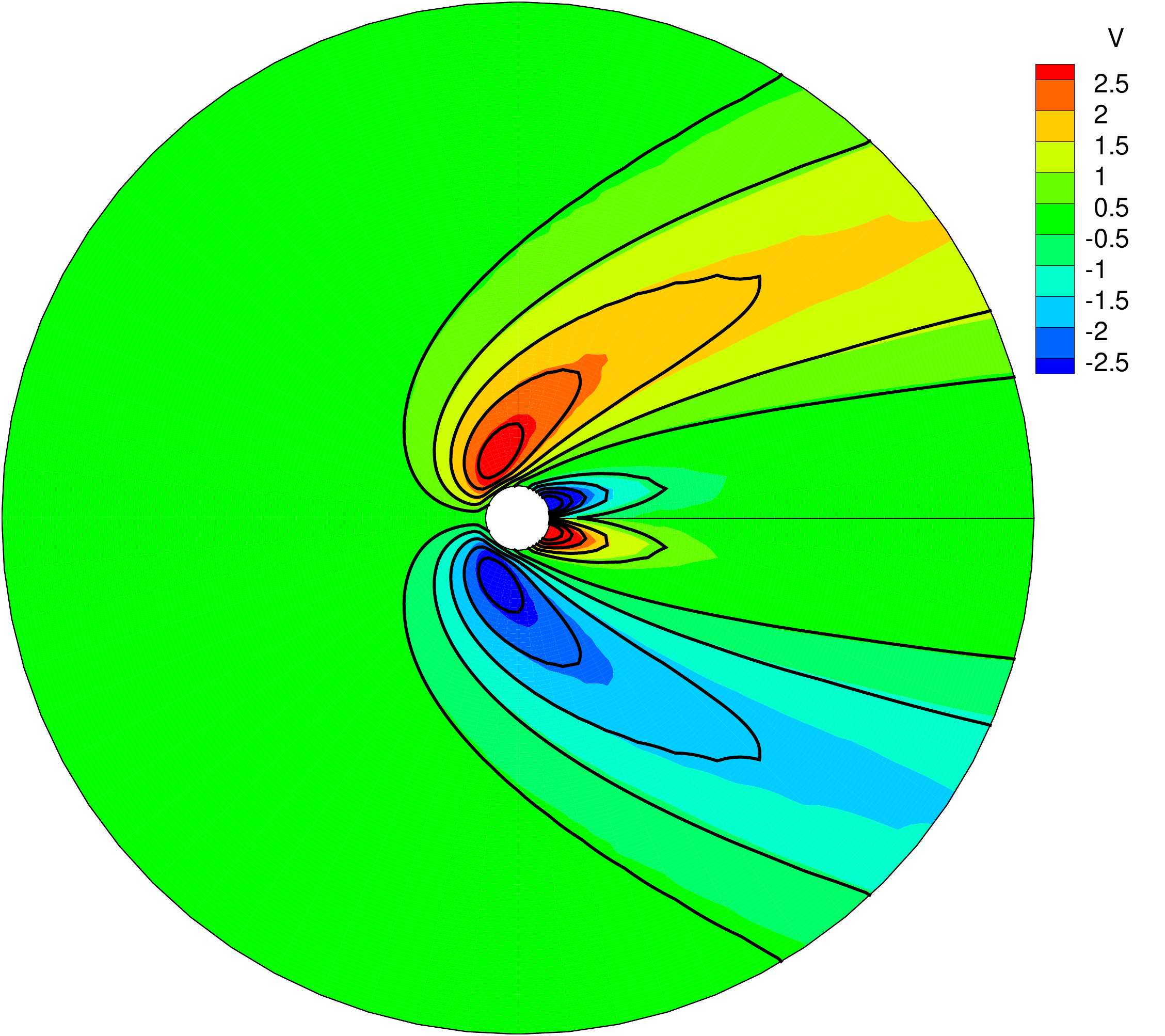}{d}
\caption{(a) Pressure, (b) temperature, (c) x directional velocity, and (d) y directional velocity contour for $\mathrm{M}=20$ and $\mathrm{Kn}=1$. The UGKWP method solution is shown in flood, and the UGKS solution is shown in contour line.}
\label{cylinder21}
\end{figure}

\begin{figure}
\centering
\includegraphics[width=0.48\textwidth]{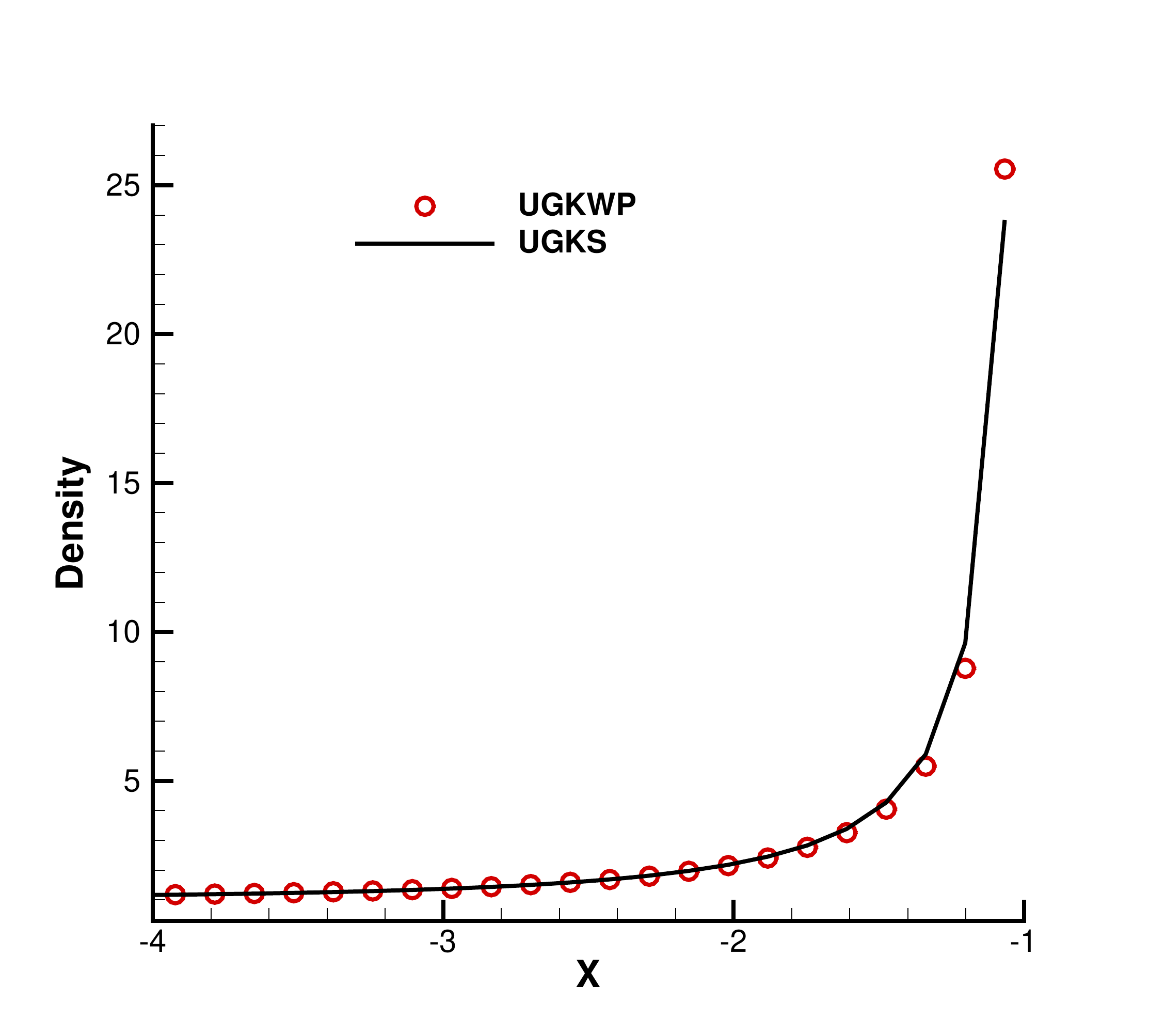}{a}
\includegraphics[width=0.48\textwidth]{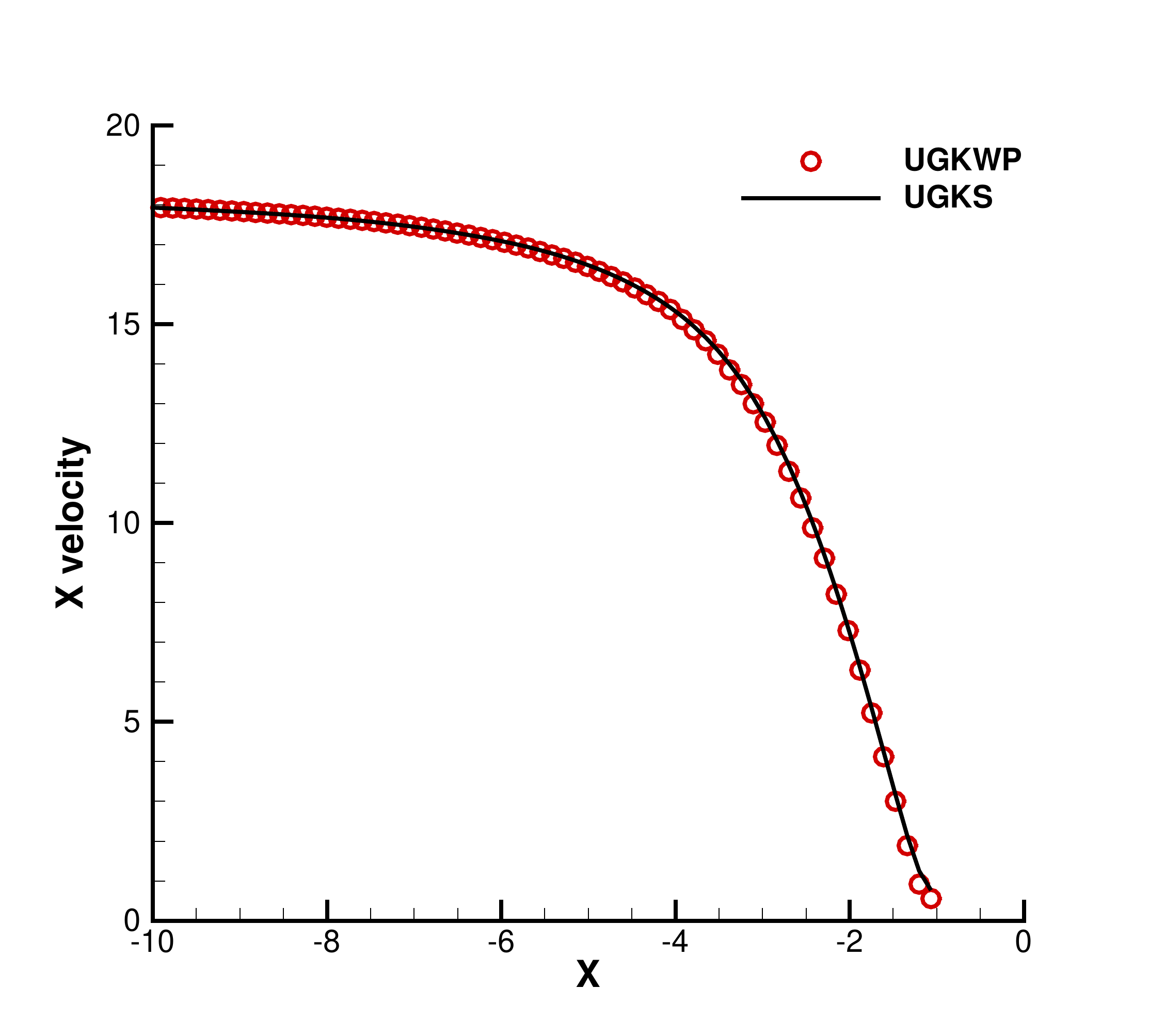}{b}
\includegraphics[width=0.48\textwidth]{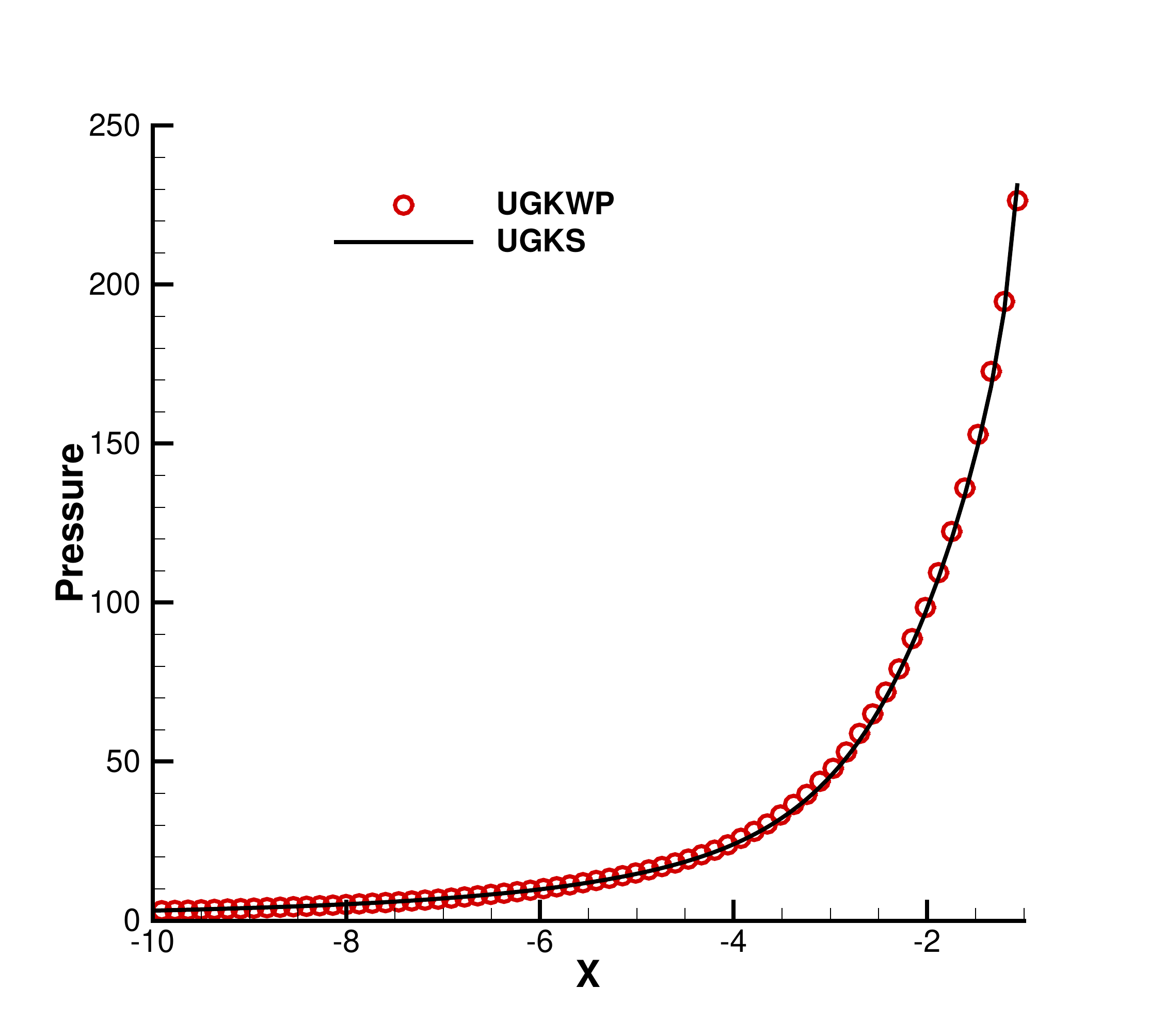}{c}
\includegraphics[width=0.48\textwidth]{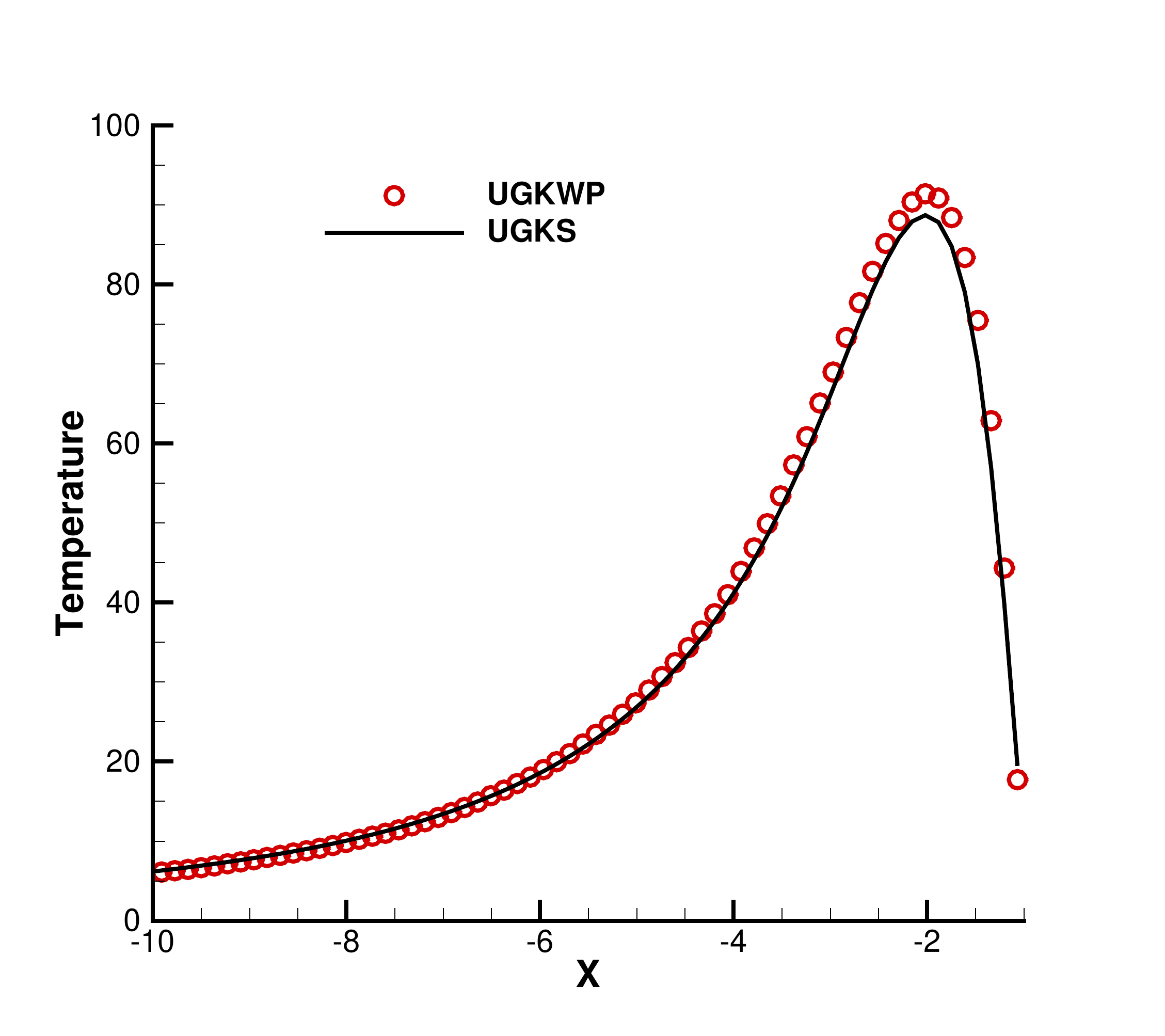}{d}
\caption{(a) Density, (b) x direction velocity, (c) pressure, (d) temperature profile along stagnation line for $\mathrm{M}=20$ and $\mathrm{Kn}=1$. The UGKWP method solution is shown in symbol, and the UGKS solution is shown in line.}
\label{cylinder22}
\end{figure}

\begin{figure}
\centering
\includegraphics[width=0.48\textwidth]{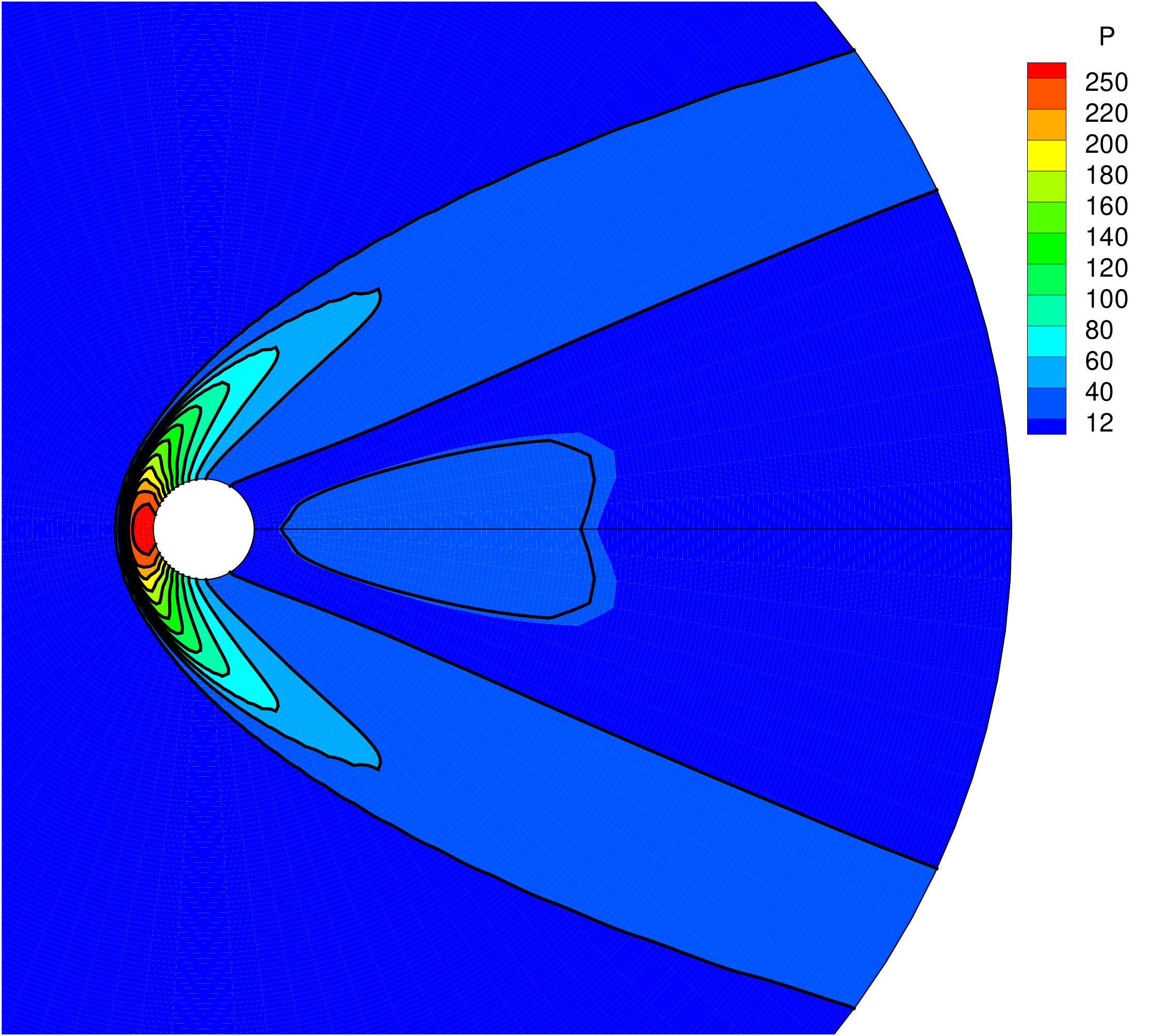}{a}
\includegraphics[width=0.48\textwidth]{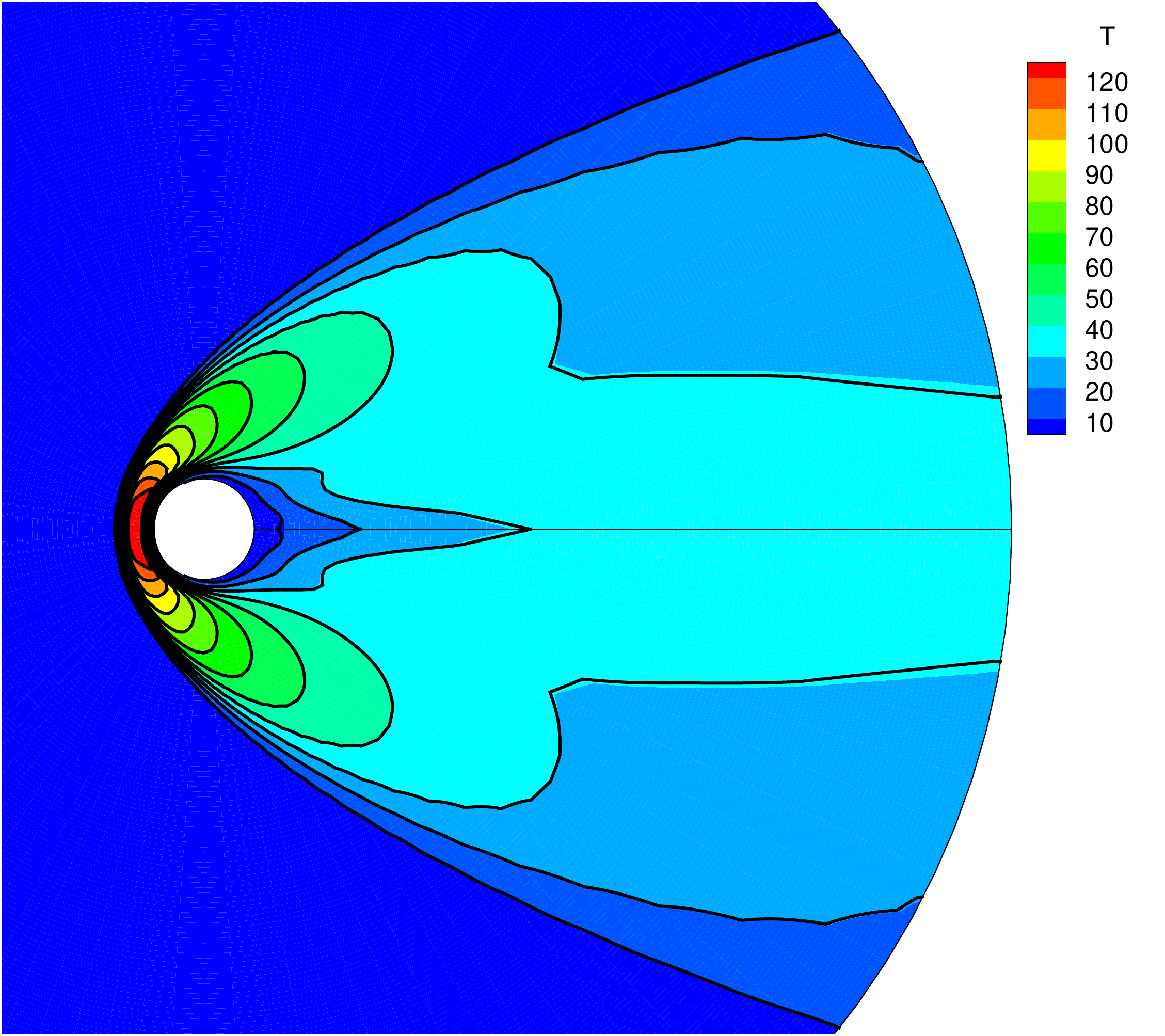}{b}
\includegraphics[width=0.48\textwidth]{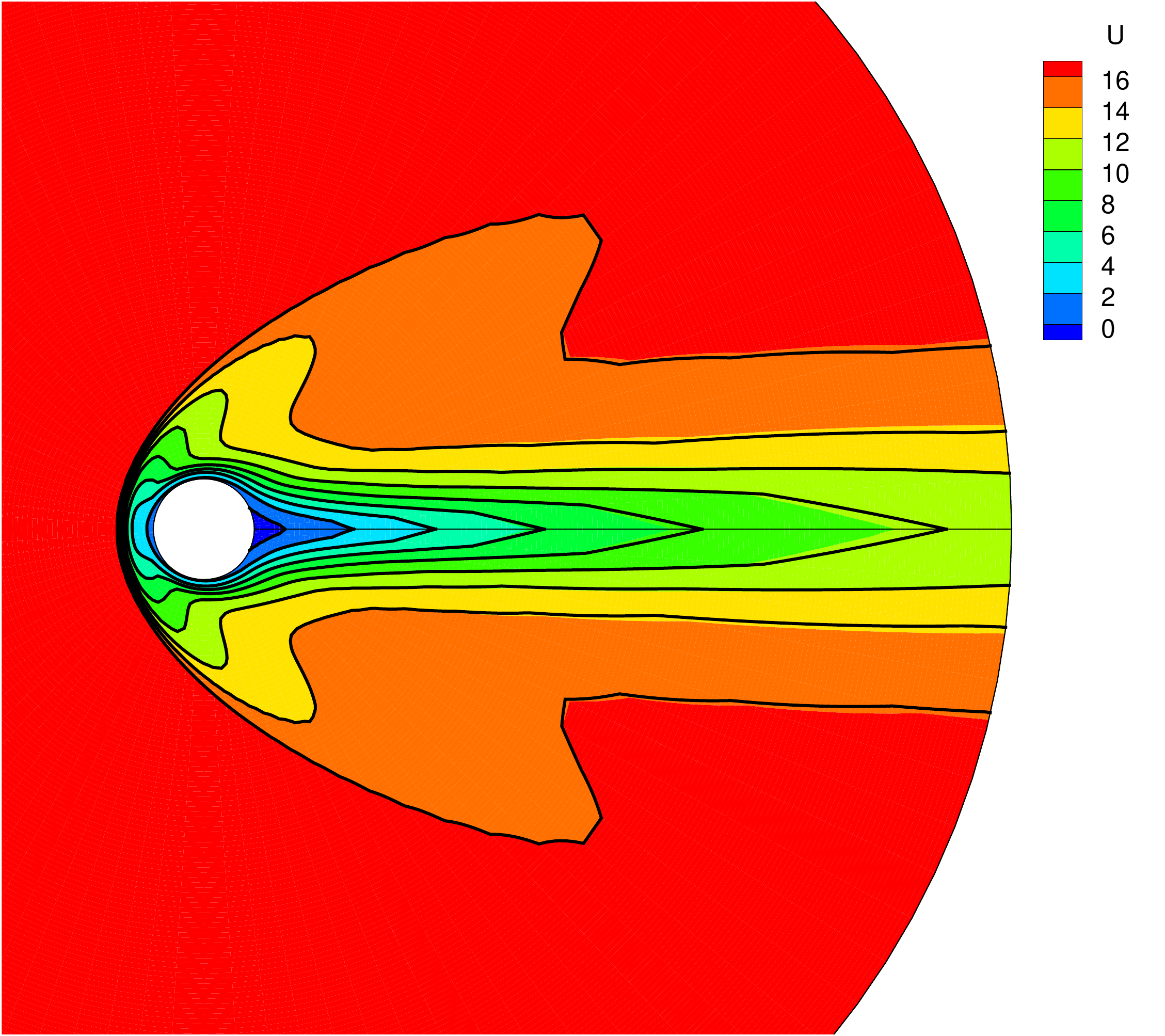}{c}
\includegraphics[width=0.48\textwidth]{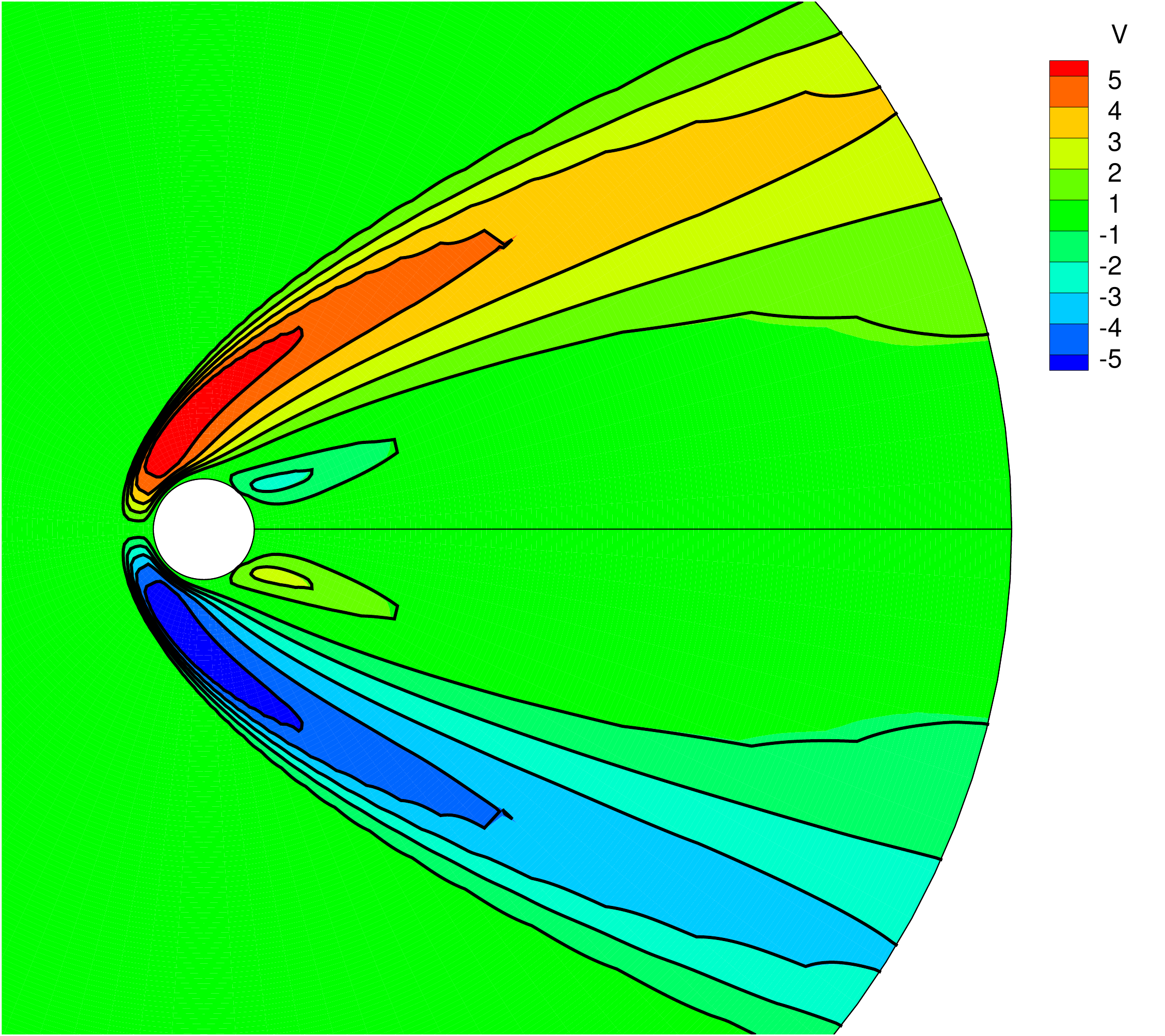}{d}
\caption{(a) Pressure, (b) temperature, (c) x directional velocity, and (d) y directional velocity contour for $\mathrm{M}=20$ and $\mathrm{Kn}=10^{-4}$. The UGKWP method solution is shown in flood, and the GKS solution is shown in contour line.}
\label{cylinder31}
\end{figure}

\begin{figure}
\centering
\includegraphics[width=0.48\textwidth]{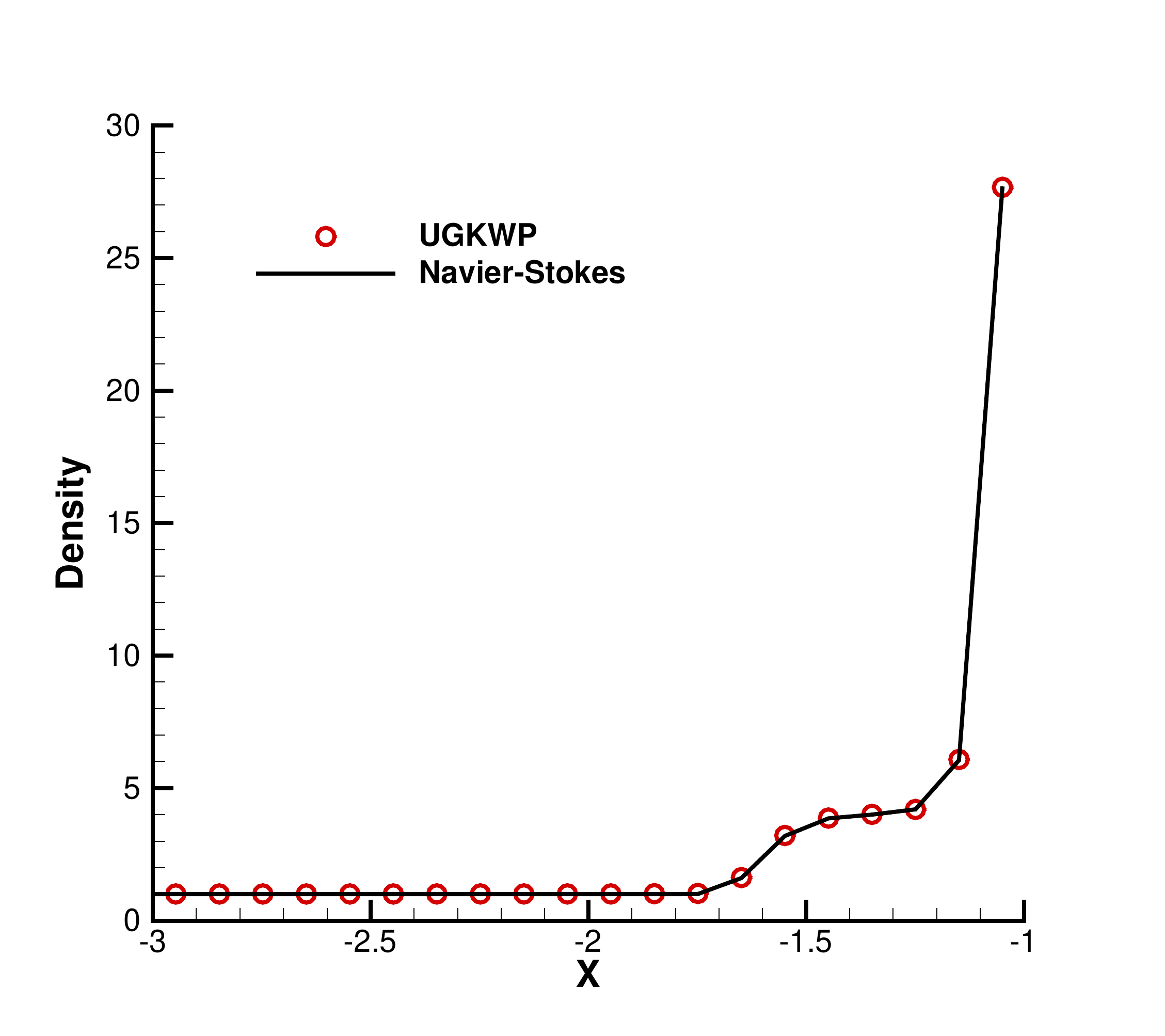}{a}
\includegraphics[width=0.48\textwidth]{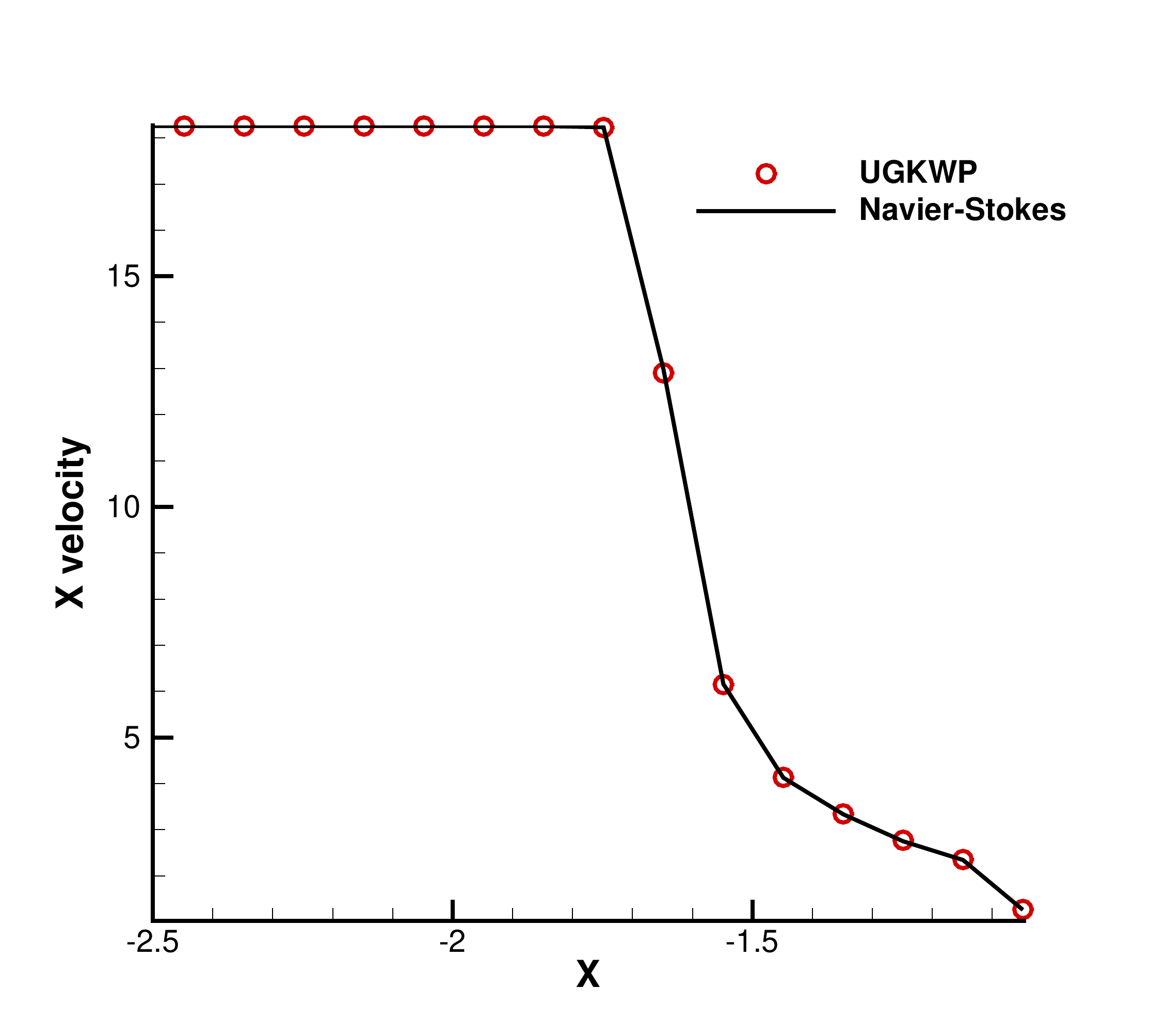}{b}
\includegraphics[width=0.48\textwidth]{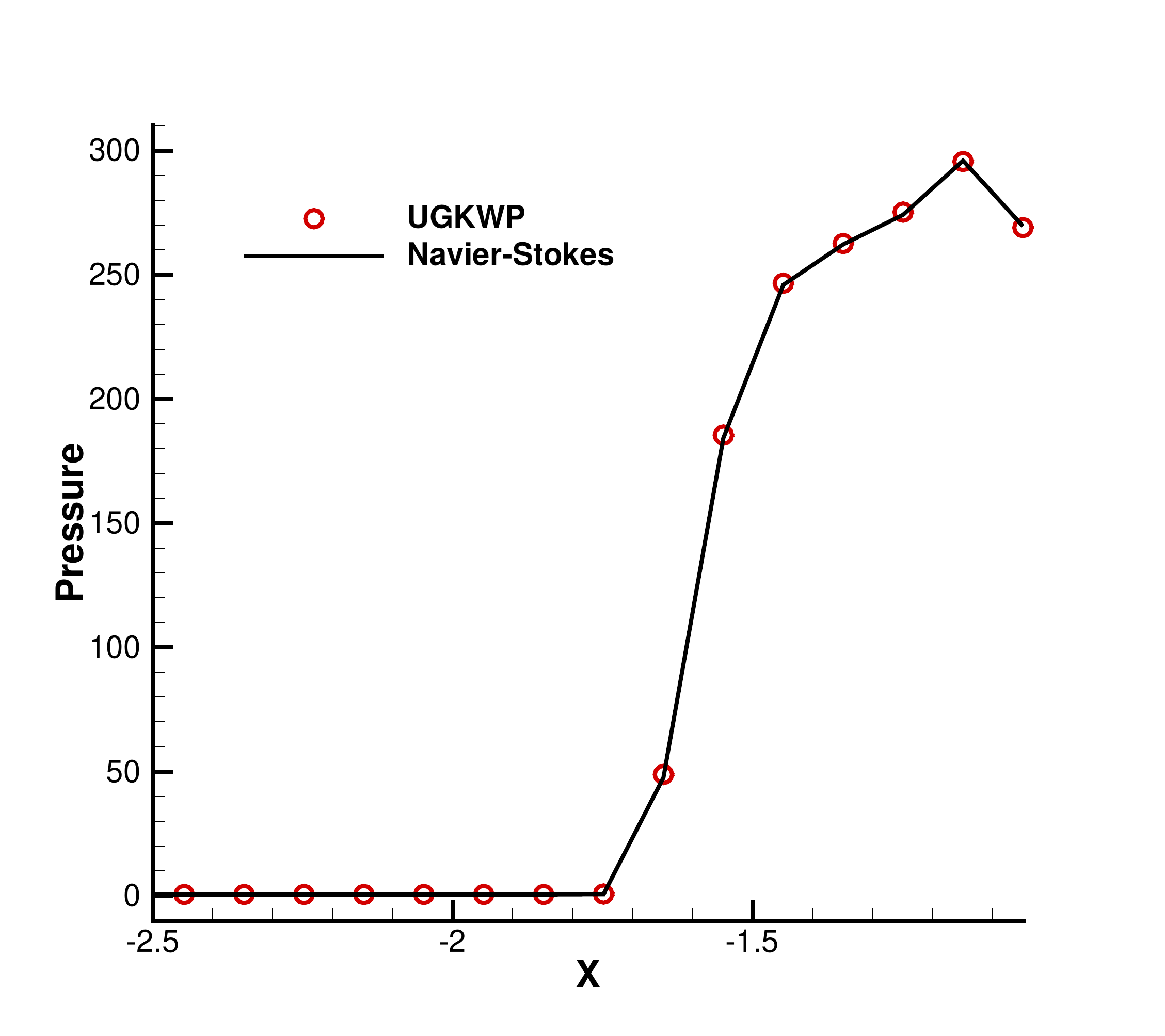}{c}
\includegraphics[width=0.48\textwidth]{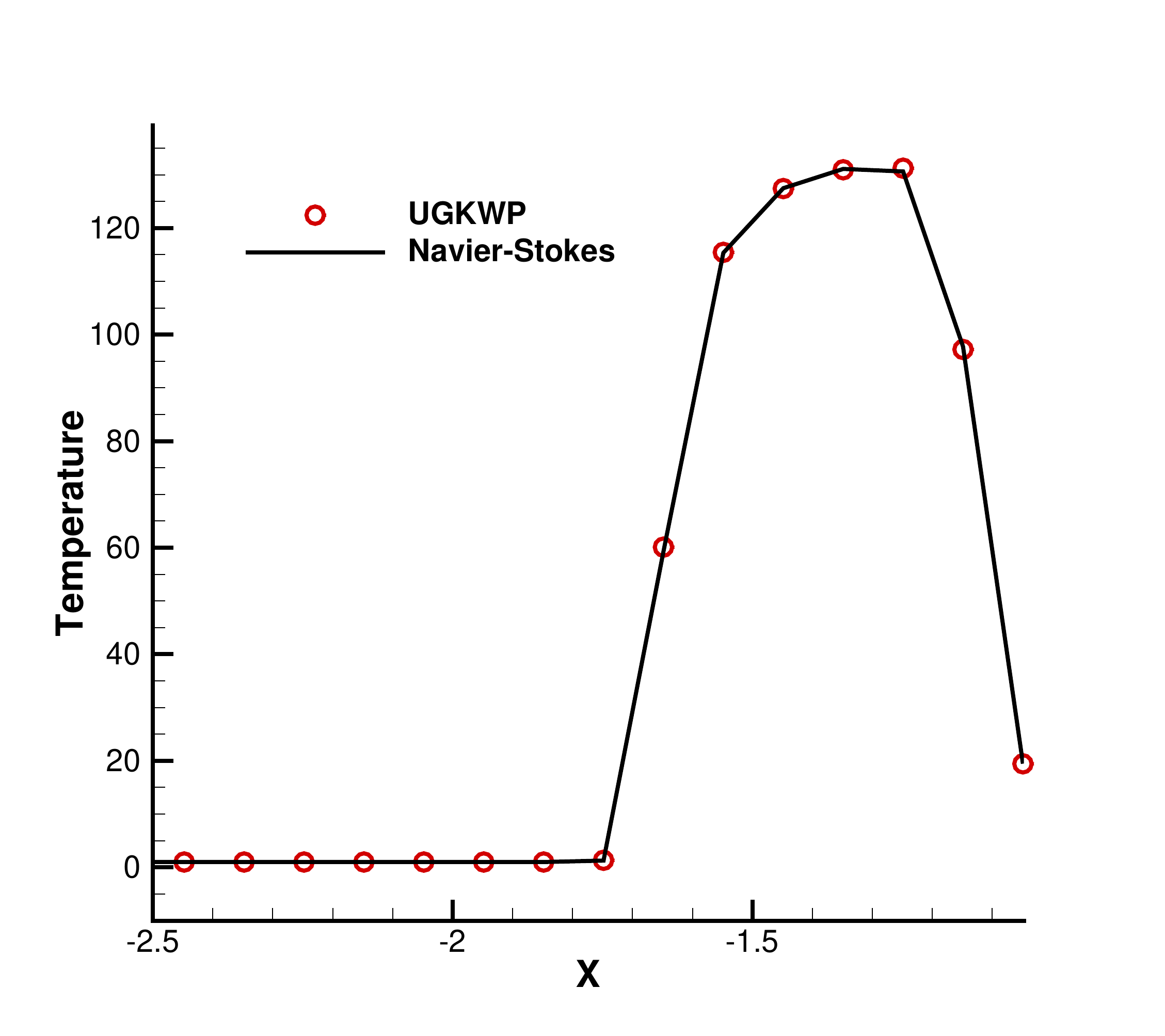}{d}
\caption{(a) Density, (b) x direction velocity, (c) pressure, (d) temperature profile along stagnation line for $\mathrm{M}=20$ and $\mathrm{Kn}=10^{-4}$. The UGKWP method solution is shown in symbol, and the GKS solution is shown in line.}
\label{cylinder32}
\end{figure}

\begin{figure}
\centering
\includegraphics[width=0.48\textwidth]{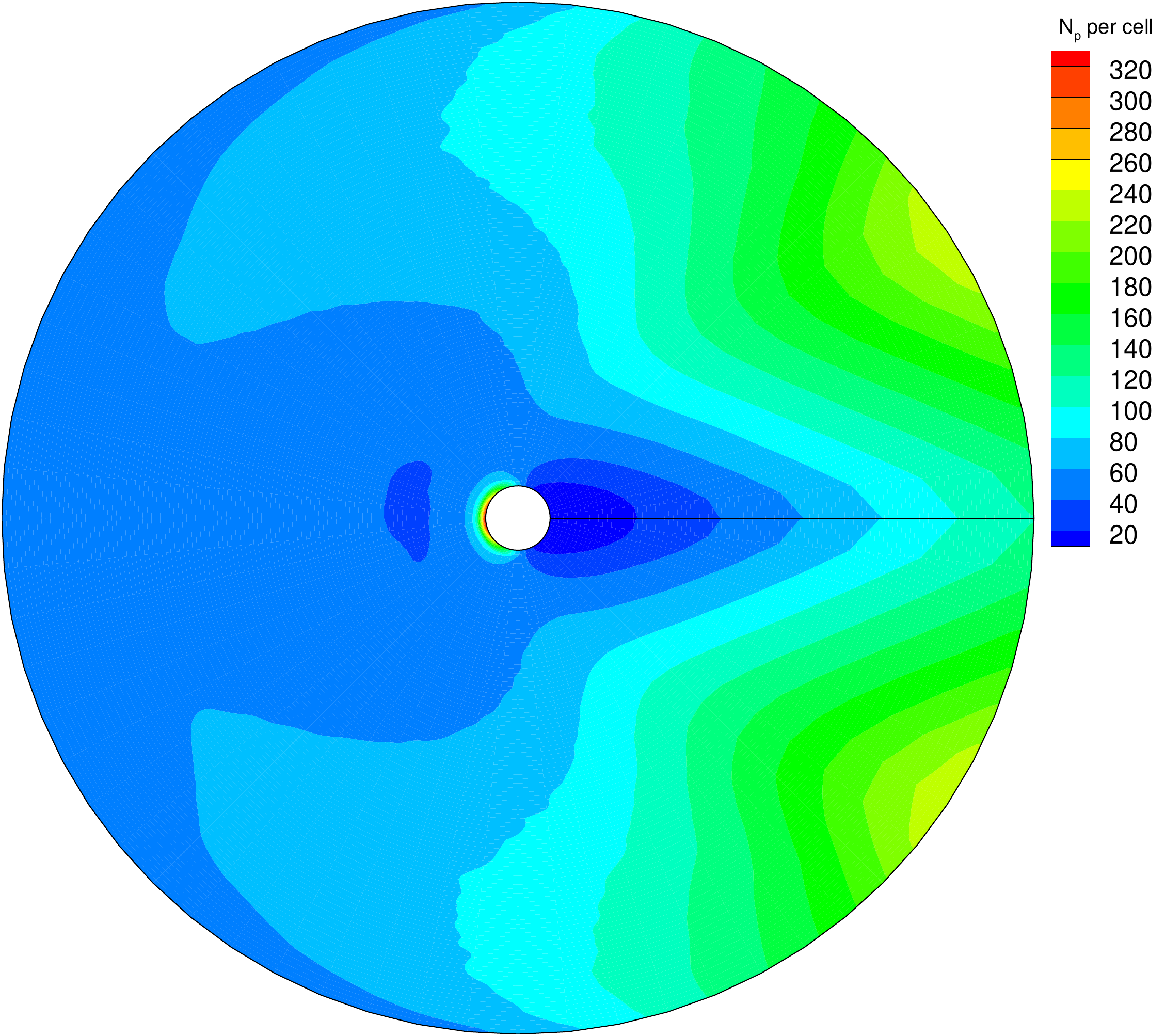}{a}
\includegraphics[width=0.48\textwidth]{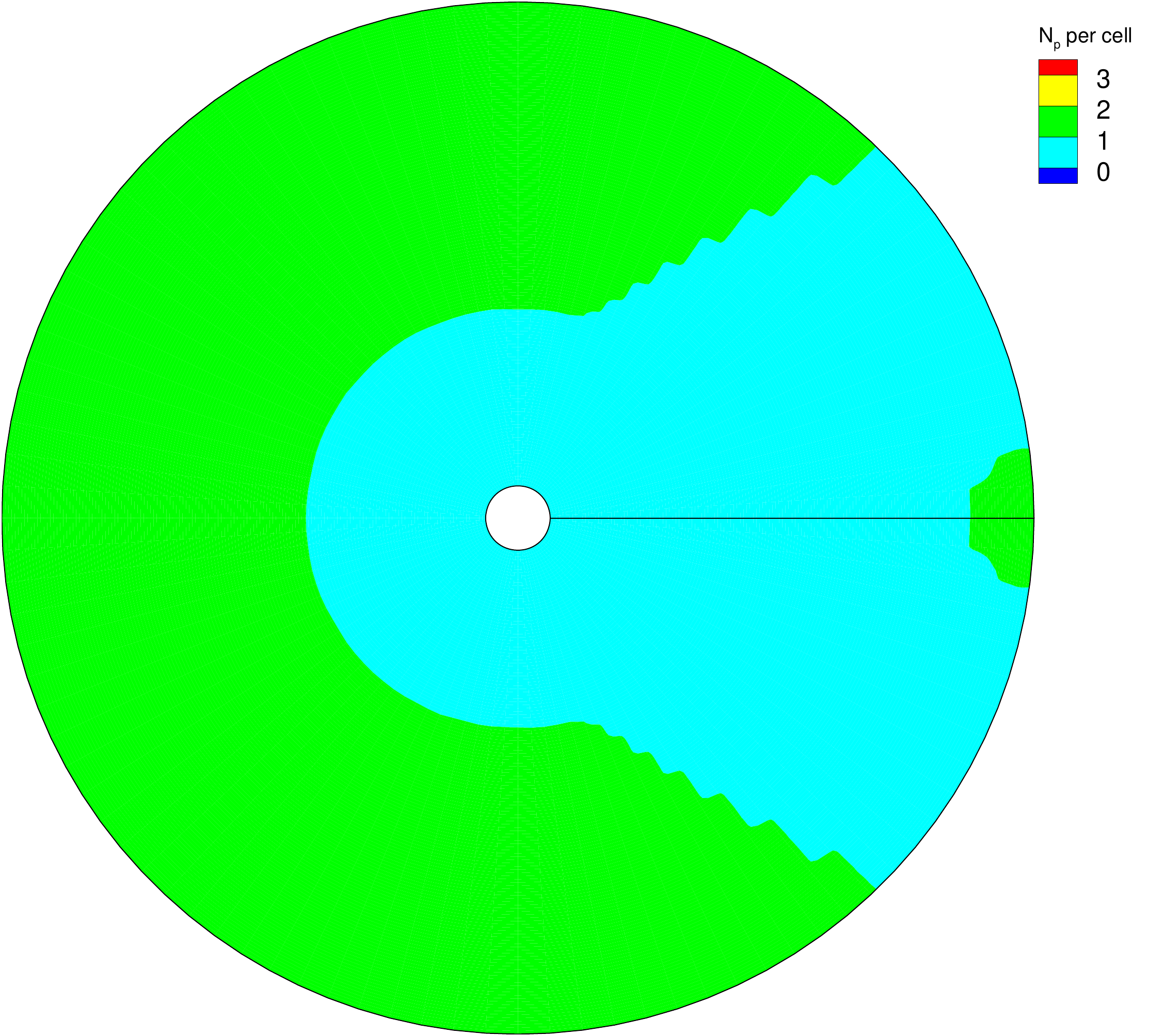}{b}
\caption{Number of simulation particles per cell for the cylinder flow with Mach number 20: (a)Kn=1.0; (b)Kn=$1.0\times10^{-4}$.}
\label{cylinder-number}
\end{figure}

\begin{figure}
\centering
\includegraphics[width=0.48\textwidth]{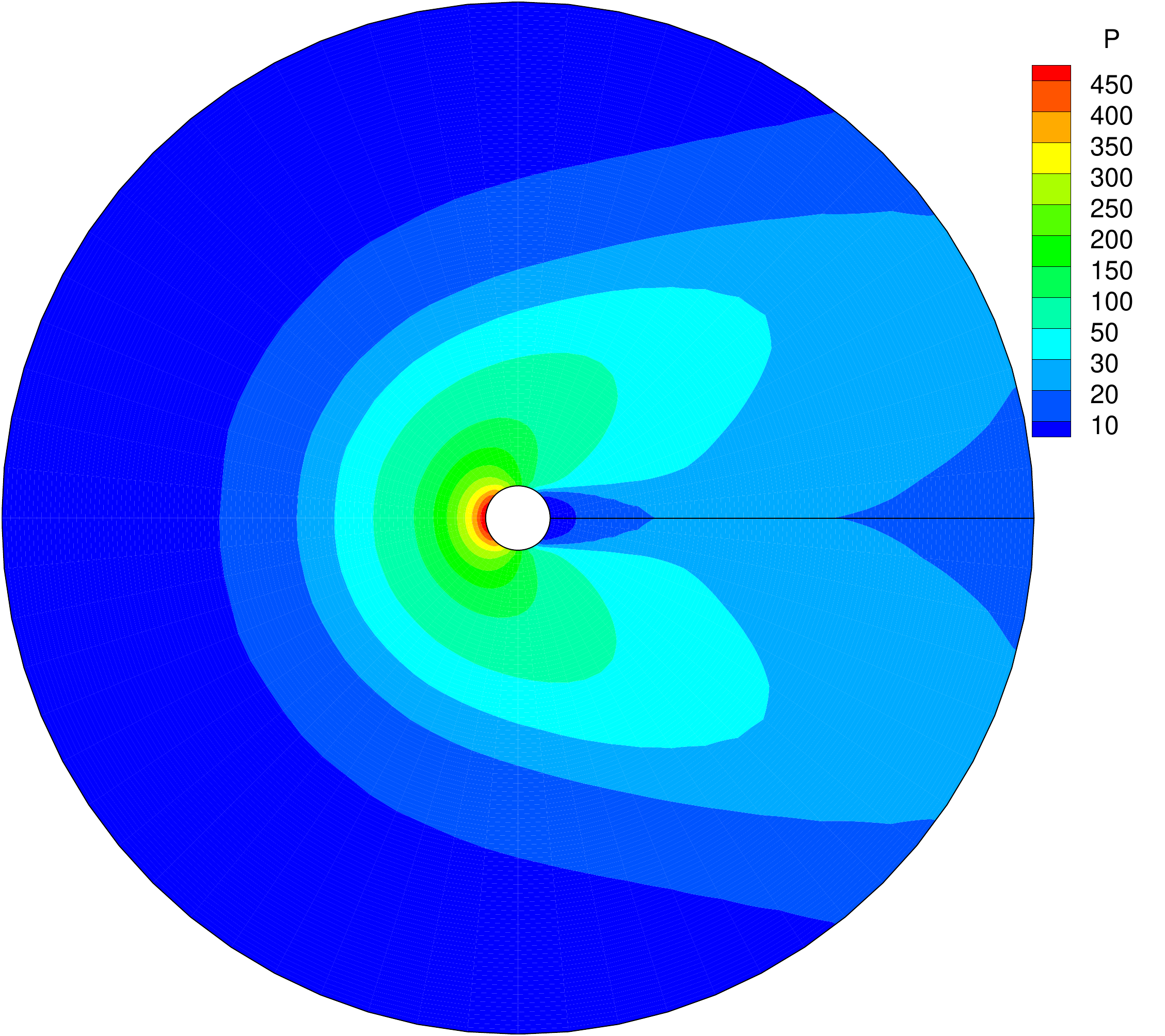}{a}
\includegraphics[width=0.48\textwidth]{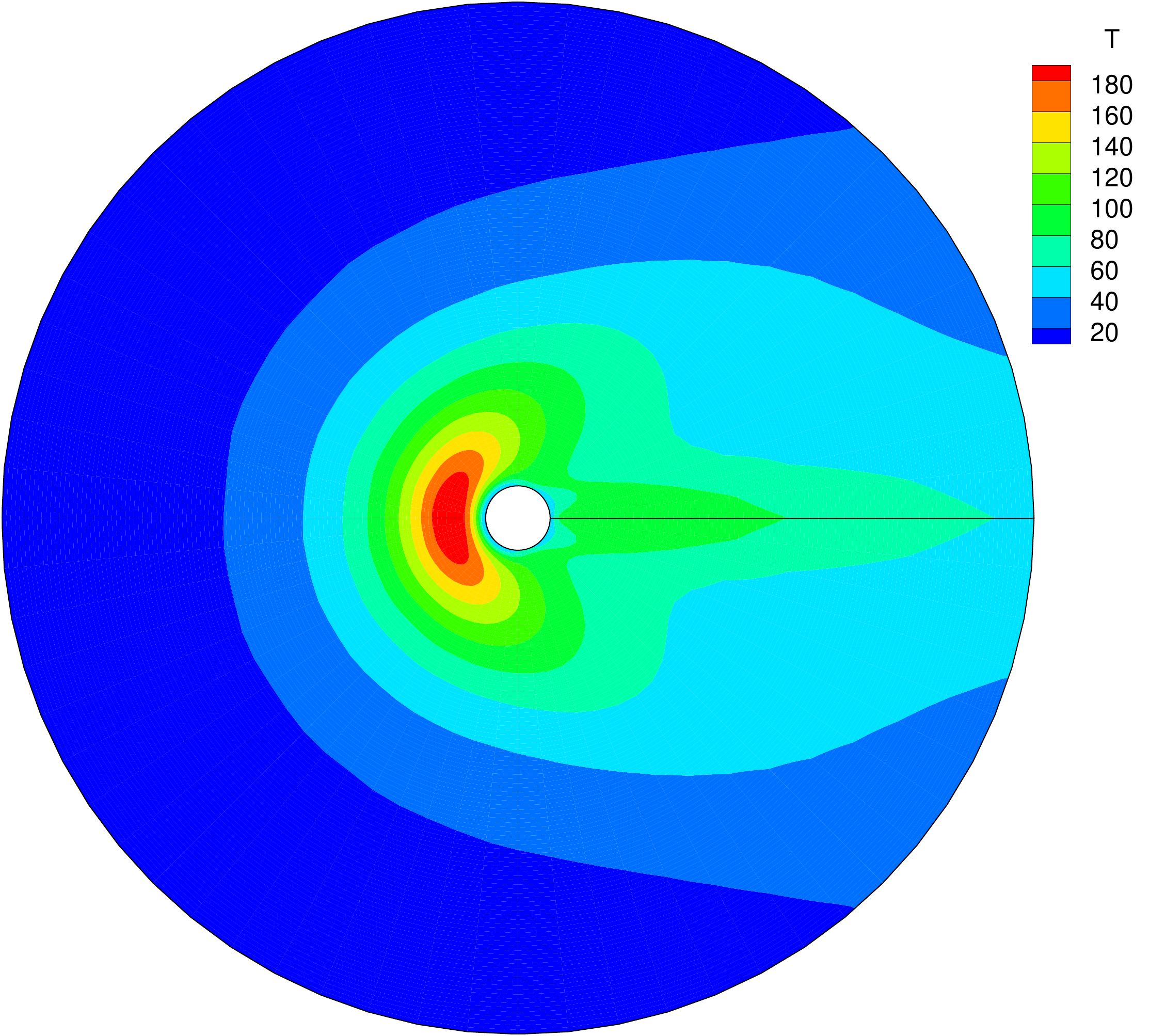}{b}
\includegraphics[width=0.48\textwidth]{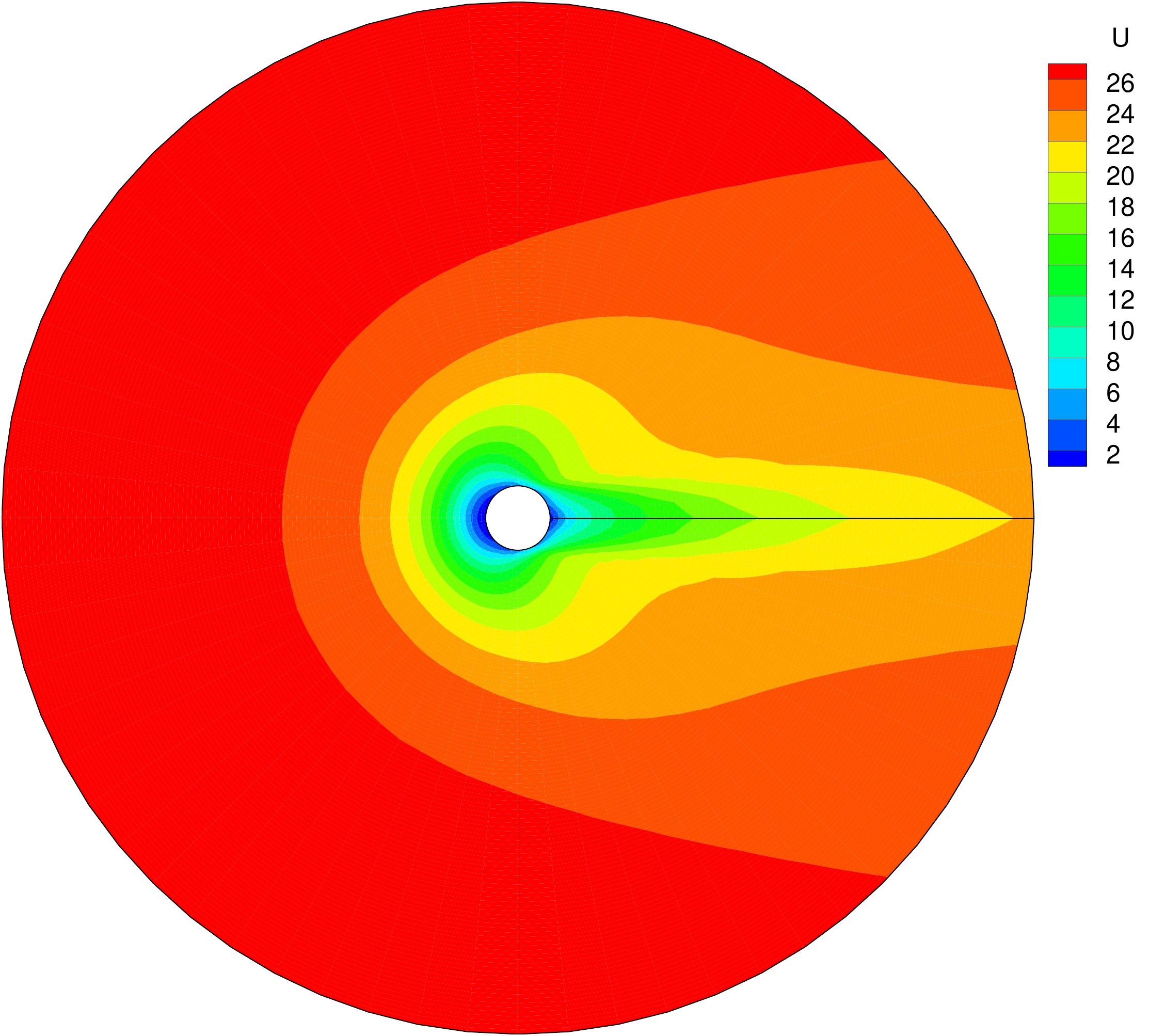}{c}
\includegraphics[width=0.48\textwidth]{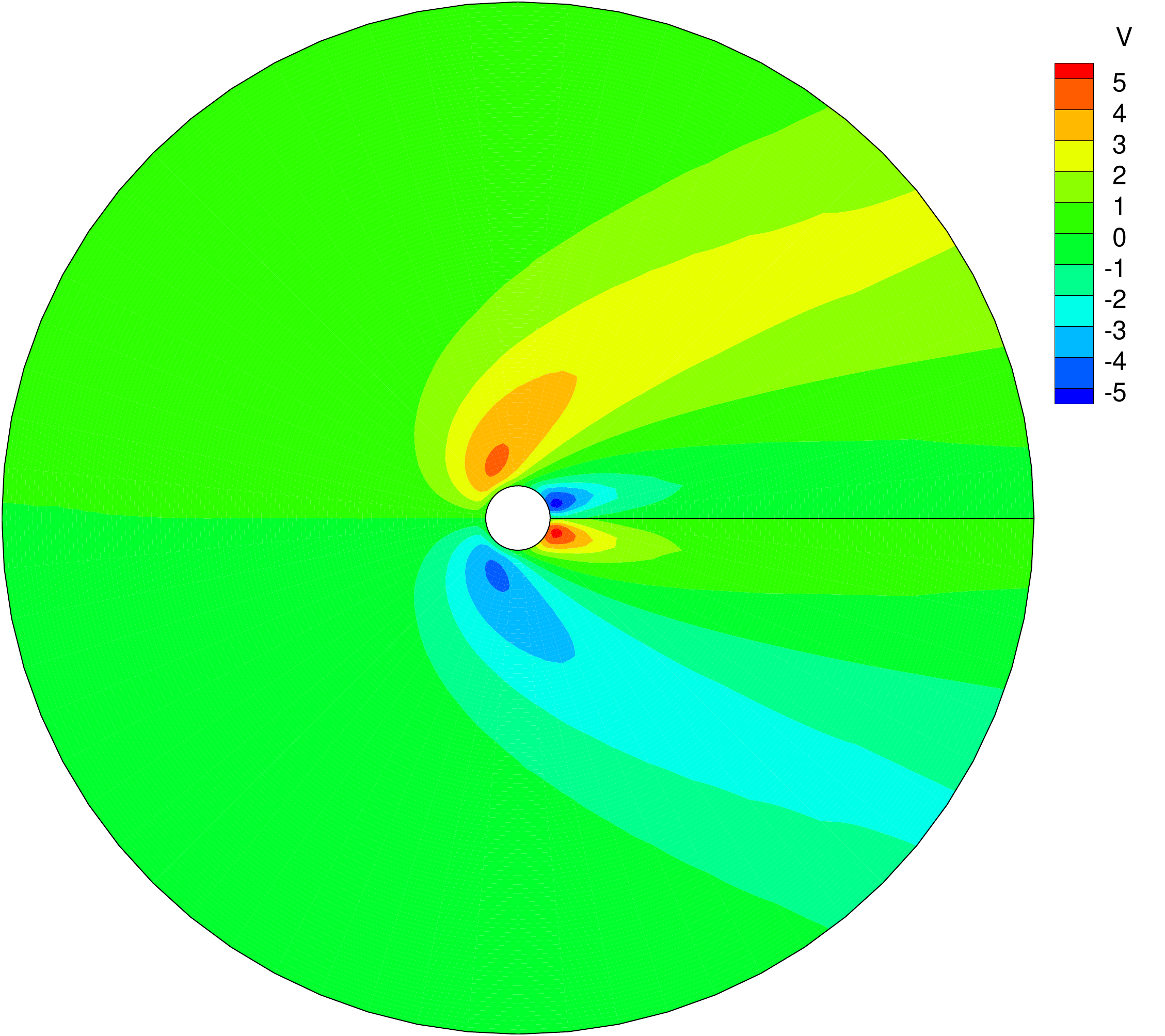}{d}
\caption{(a) Pressure, (b) temperature, (c) x directional velocity, and (d) y directional velocity contour of UGKWP method for $\mathrm{M}=30$ and $\mathrm{Kn}=1$.}
\label{cylinder4}
\end{figure}

\begin{figure}
\centering
\includegraphics[width=0.48\textwidth]{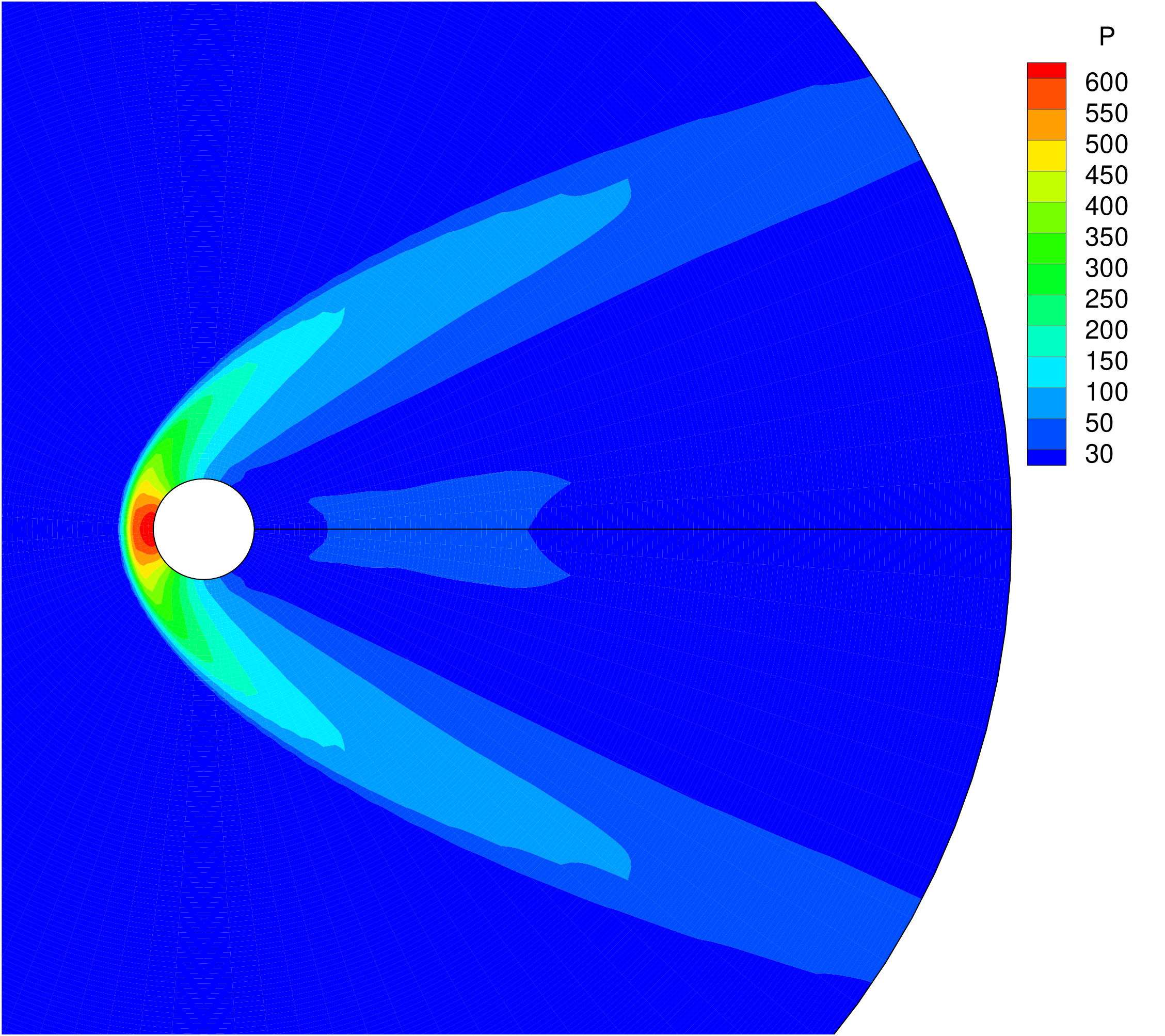}{a}
\includegraphics[width=0.48\textwidth]{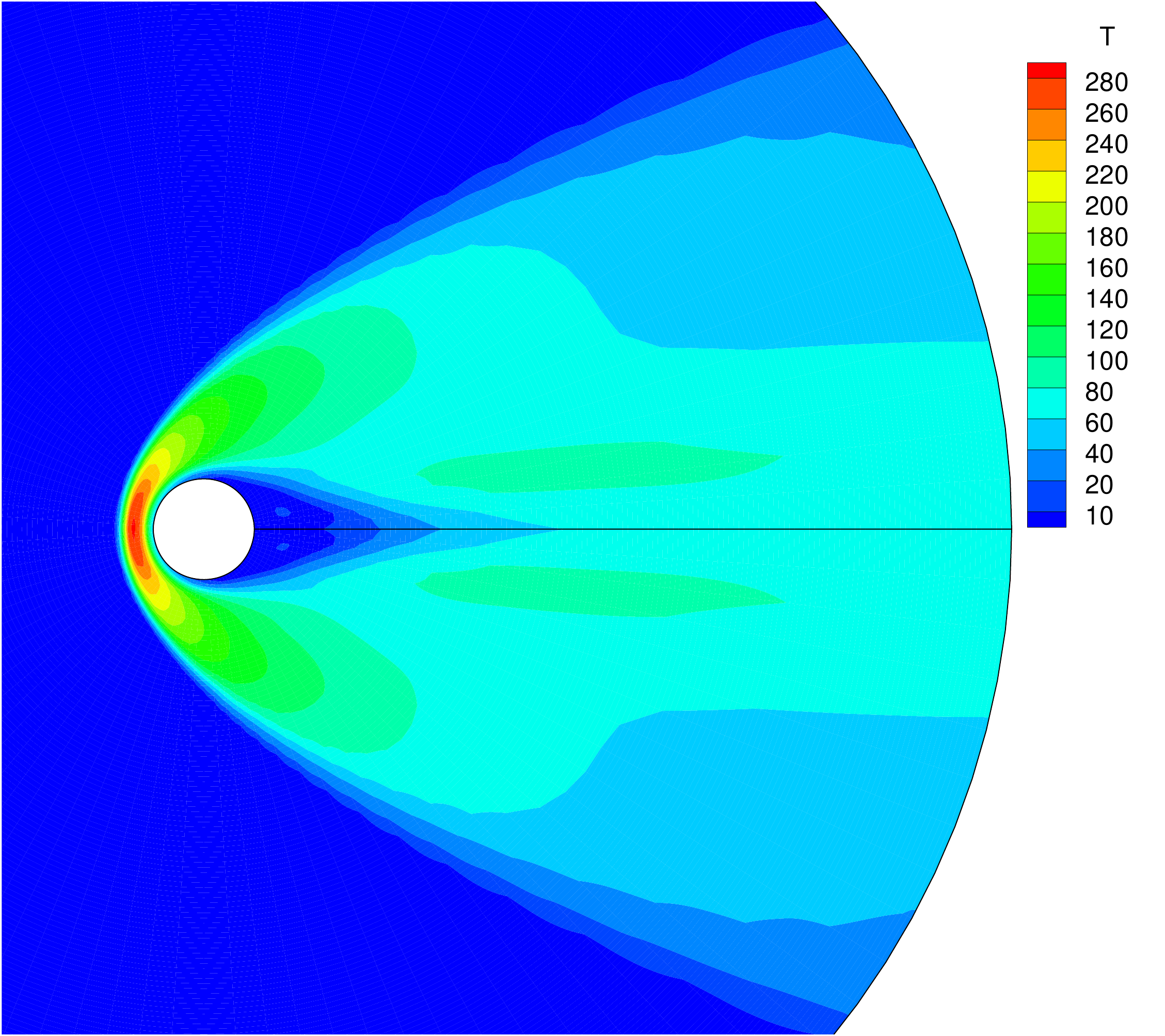}{b}
\includegraphics[width=0.48\textwidth]{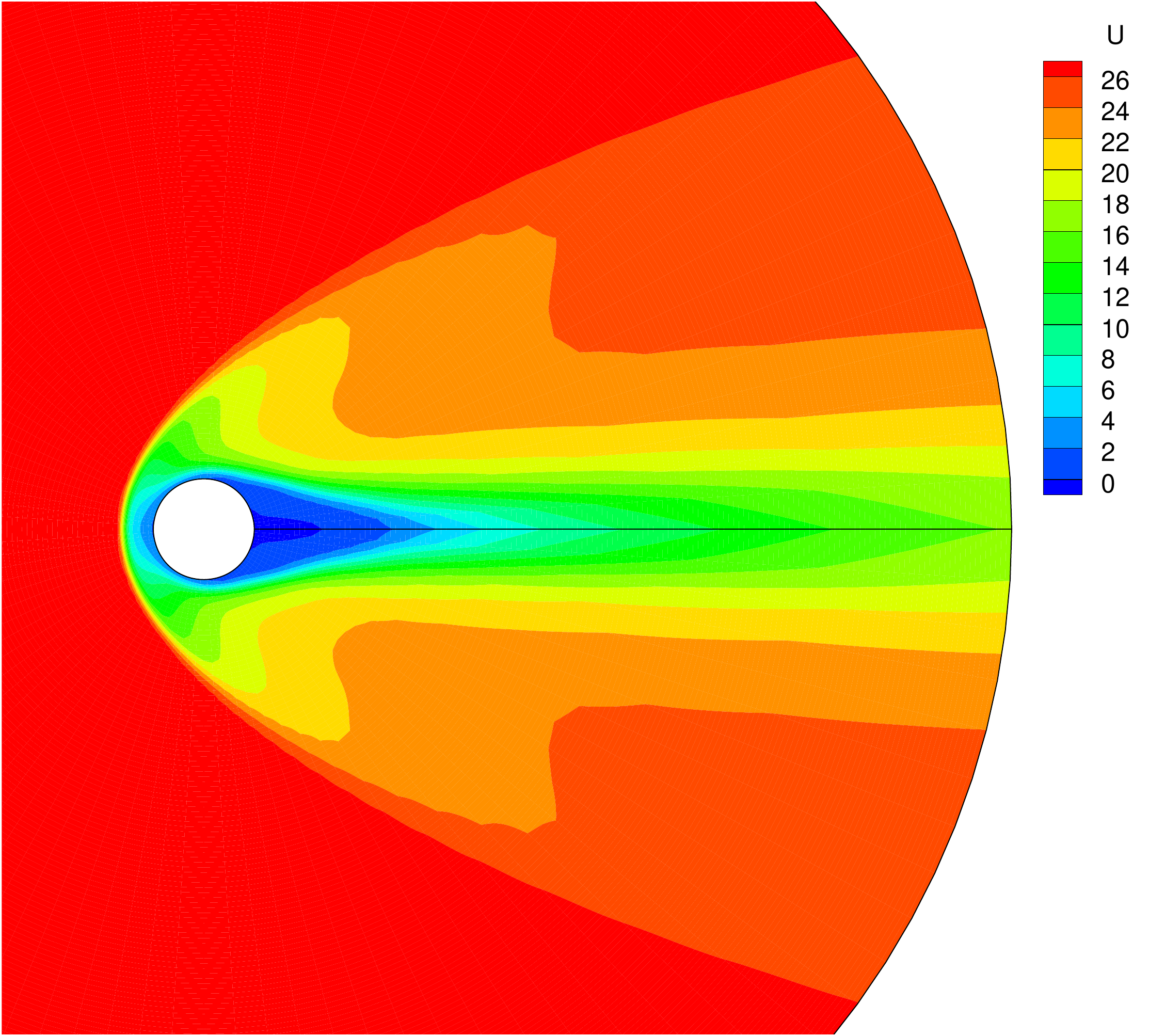}{c}
\includegraphics[width=0.48\textwidth]{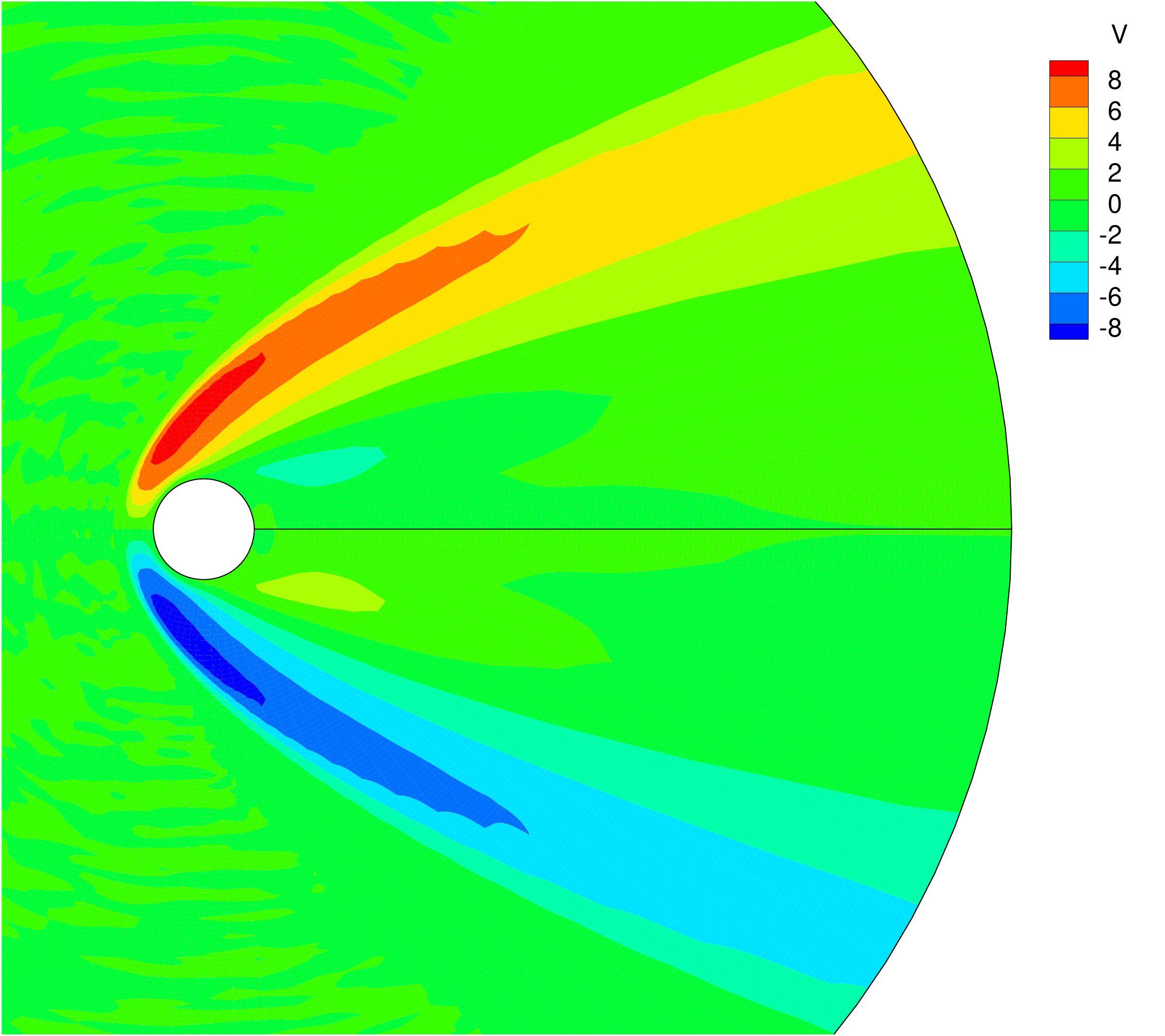}{d}
\caption{(a) Pressure, (b) temperature, (c) x directional velocity, and (d) y directional velocity contour of UGKWP method for $\mathrm{M}=30$ and $\mathrm{Kn}=10^{-4}$.}
\label{cylinder5}
\end{figure}

We compare the results and computational efficiency of UGKWP method to UGKS/GKS.
For the case of $\mathrm{M}=5$ and $\mathrm{Kn}=0.1$, the physical domain of the UGKWP method is discretized by a mesh with $64\times64$ cells, and the mass of simulation particle is set $m_p=1.52\times10^{-3}$.
For UGKS, the same discretization in physical space is used and the velocity space is $[-10,10]\times[-10,10]$ discretized by $100\times100$ velocity points. The CFL condition for both UGKS and UGKWP method is set to be 0.9. Fig. \ref{cylinder11} shows the contours of steady state pressure, temperature, and velocity, where the flood is the UGKWP method solution and lines are the UGKS solution. The density, velocity, pressure, temperature profiles along the stagnation line are shown in Fig. \ref{cylinder12}, where the solution of UGKWP method is denoted by symbols and the solution of UGKS is shown in lines.
To get the steady state solution, the computational time for UGKS is 10.9 hours and the UGKWP method is 21.5 minutes (including the averaging procedure).
The memory cost for UGKS is 3.4 GB and for UGKWP method is 75 MB.
The computational time of UGKWP method is about 30 times faster than explicit UGKS,
and the memory cost is 46 times less than UGKS.
For the case of $\mathrm{M}=20$ and $\mathrm{Kn}=1.0$, the discretization of the physcial domain is 64 cells along azimuth direction, and 110 cells along radial direction. For UGKWP method, the mass of simulation particle is set $m_p=1.0\times 10^{-3}$, and the velocity space for UGKS is $[-50,50]\times[-50,50]$ discretized by $200\times200$ velocity points. Fig. \ref{cylinder21} shows the contours of steady state pressure, temperature, and velocity, where the flood is the UGKWP method solution and lines are the UGKS solution. The density, velocity, pressure, temperature profiles along the stagnation line are shown in Fig. \ref{cylinder22}, where the solution of UGKWP method is denoted by symbols and the solution of UGKS is shown in lines.
The averaged number of simulation particle per cell is shown in Fig. \ref{cylinder-number}-(a).
To get the steady state solution, the computational time for UGKS is about 429 hours and the UGKWP method is 36.1 minutes (including the averaging procedure).
The memory cost for UGKS is 22.3 GB and for UGKWP method is less than 100 MB.
The computational time of UGKWP method is about 713 times faster than explicit UGKS,
and the memory cost is 228 times less than UGKS.
For the case of $\mathrm{M}=20$ and $\mathrm{Kn}=10^{-4}$, a discretization of $100\times150$ cells is used in physical space and the simulation particle mass is set $m_p=4.7\times 10^{-3}$. The solution of the UGKWP method is compared with the Navier-Stokes solution calculated by the gas-kinetic scheme, and the solution contour and profile along the stagnation line are shown in Fig. \ref{cylinder31} and \ref{cylinder32} respectively.
The averaged number of simulation particle per cell is shown in Fig. \ref{cylinder-number}-(b).
The simulation time for UGKWP method is 17.2 minutes which is comparable to GKS.
To demonstrate the capability of UGKWP method in the simulation of multiscale high Mach number flow, we calculate the cylinder flow at Mach number 30 and Knudsen number $\mathrm{Kn}=1.0,10^{-4}$. The solution contours are shown in Fig. \ref{cylinder4}-\ref{cylinder5}. The computational time for the rarefied case is around 40 minutes and for $\mathrm{Kn}=10^{-4}$ is around 20 minutes, both including the averaging procedure. The comparison between the computational cost between UGKWP method and UGKS is shown in table \ref{efficiency1}-\ref{efficiency2}.

\begin{table}[ht] \footnotesize
\centering
\begin{tabular}{ |c | c | c| c|c|c|}
  \hline			
  $\mathrm{M}$ and $\mathrm{Kn}$ of cylinder flow & Time of UGKS &Time of UGKWP & Time ratio $\frac{\text{UGKS}}{\text{UGKWP}}$ \\
  \hline	
  $\mathrm{M}=5$,\ $\mathrm{Kn}=0.1$& 10.9 hours & 21.5 minutes& 30  \\
  \hline	
  $\mathrm{M}=20$,\ $\mathrm{Kn}=1$& 429 hours & 36.1 minutes & 713\\
  \hline
\end{tabular}
\caption{Comparison between the computational time between UGKWP method and UGKS}
\label{efficiency1}
\end{table}
\begin{table}[ht]\footnotesize
\centering
\begin{tabular}{ |c | c | c| c|c|c|}
  \hline			
  $\mathrm{M}$ and $\mathrm{Kn}$ of cylinder flow & Memory of UGKS &Memory of UGKWP & Memory ratio $\frac{\text{UGKS}}{\text{UGKWP}}$ \\
  \hline	
  $\mathrm{M}=5$,\ $\mathrm{Kn}=0.1$& 3.4 GB & 75 MB & 46 \\
  \hline	
  $\mathrm{M}=20$,\ $\mathrm{Kn}=1$& 22.3GB & 100 MB&228\\
  \hline
\end{tabular}
\caption{Comparison between the memory cost between UGKWP method and UGKS}
\label{efficiency2}
\end{table}

\subsection{Cavity flow}
To show the capability of UGKWP method in simulating the multiscale low speed flow, we calculate the lid-driven cavity flow at different Knudsen number $\mathrm{Kn}=1.0, 0.75, 1.42\times10^{-4}$.
The gaseous medium is assumed to consist of monatomic argon gas, which is modeled by the VHS model. The particle parameters of argon and the formulation of viscosity coefficient are the same as in the calculation of cylinder flow.
The wall temperature is kept $T_w=273$ K and the top lid is moving with a fixed velocity of $U_w=50$ m/s.
The dimensionless quantities are used with respect to the reference length as the cavity, the reference temperature as the initial gas temperature $T_{ref}=T_0$, the reference velocity $U_{ref}=\sqrt{2RT_0}$, the reference time $t_{ref}=L_{ref}/U_{ref}$, the reference density as the gas initial density $\rho_{ref}= \rho_0$.

\begin{figure}
\centering
\includegraphics[width=0.48\textwidth]{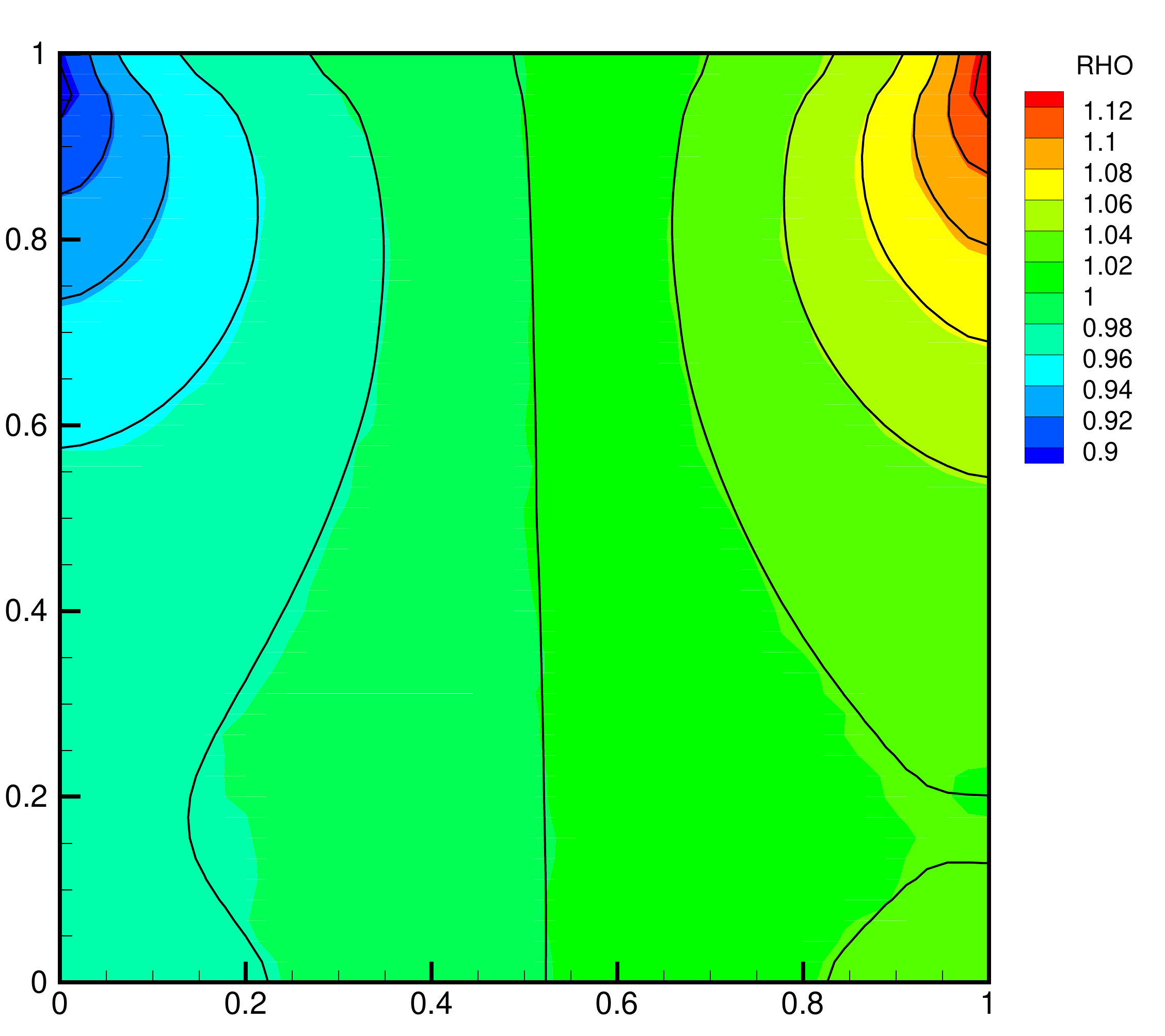}{a}
\includegraphics[width=0.48\textwidth]{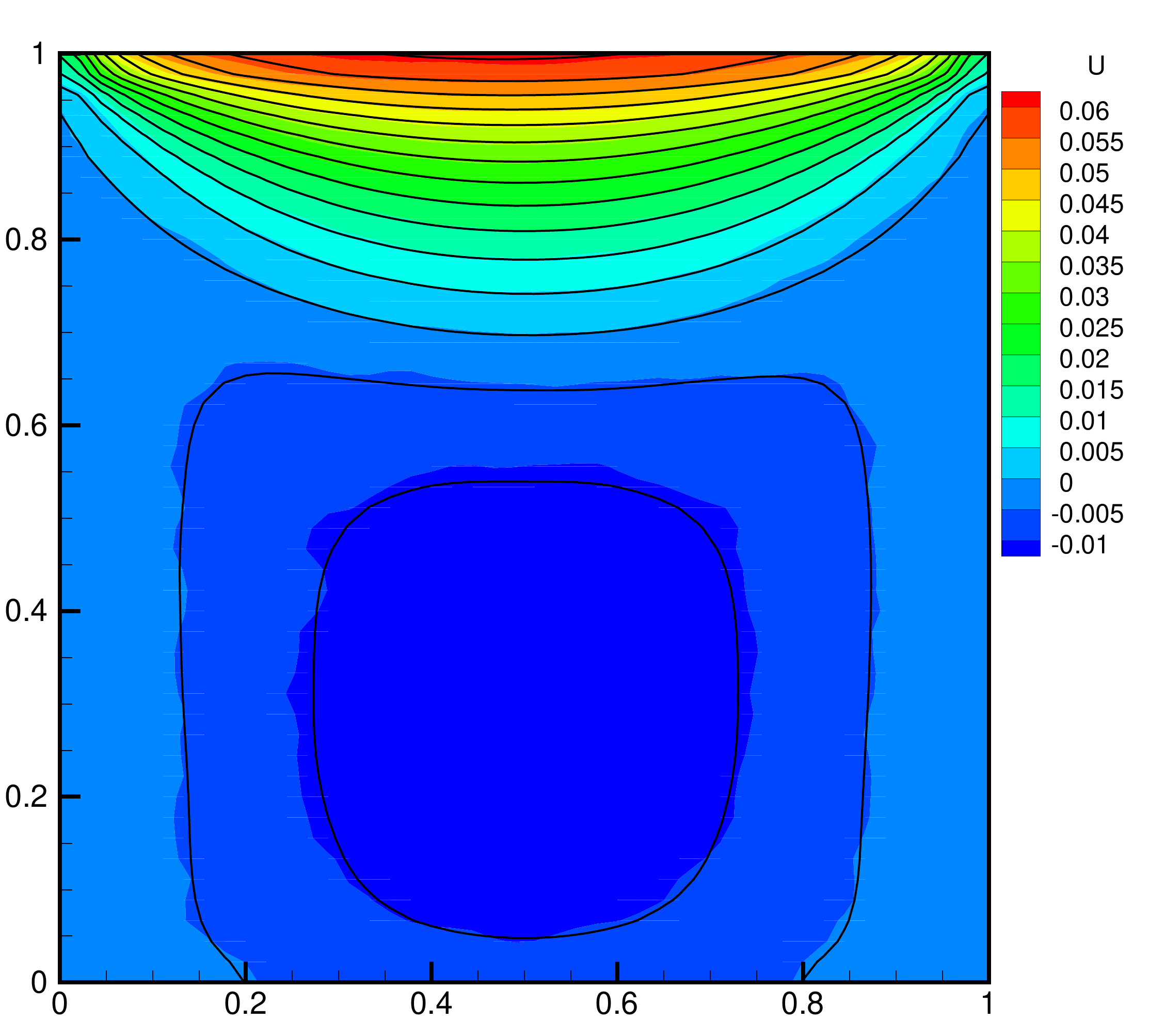}{b}
\includegraphics[width=0.48\textwidth]{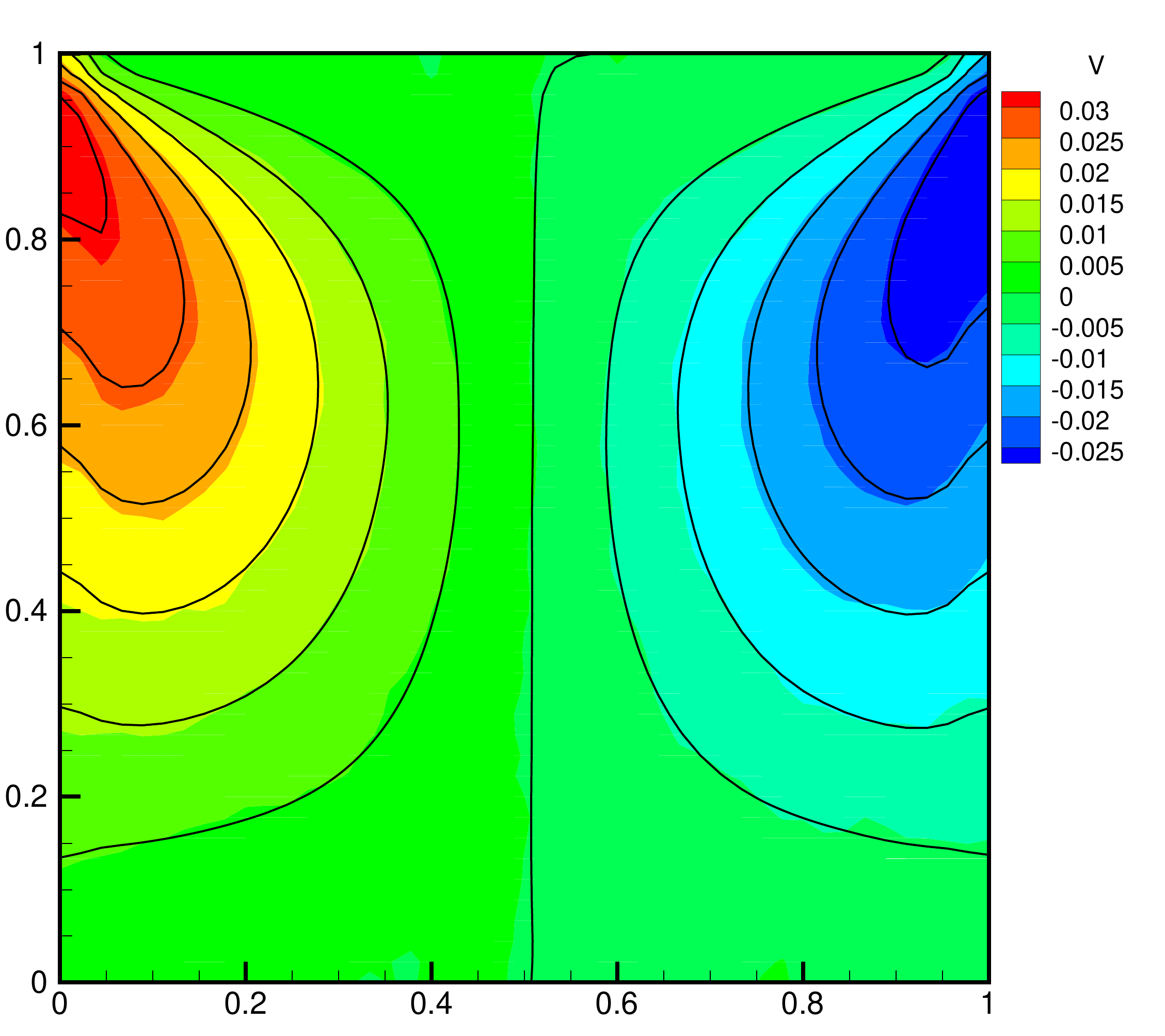}{c}
\includegraphics[width=0.48\textwidth]{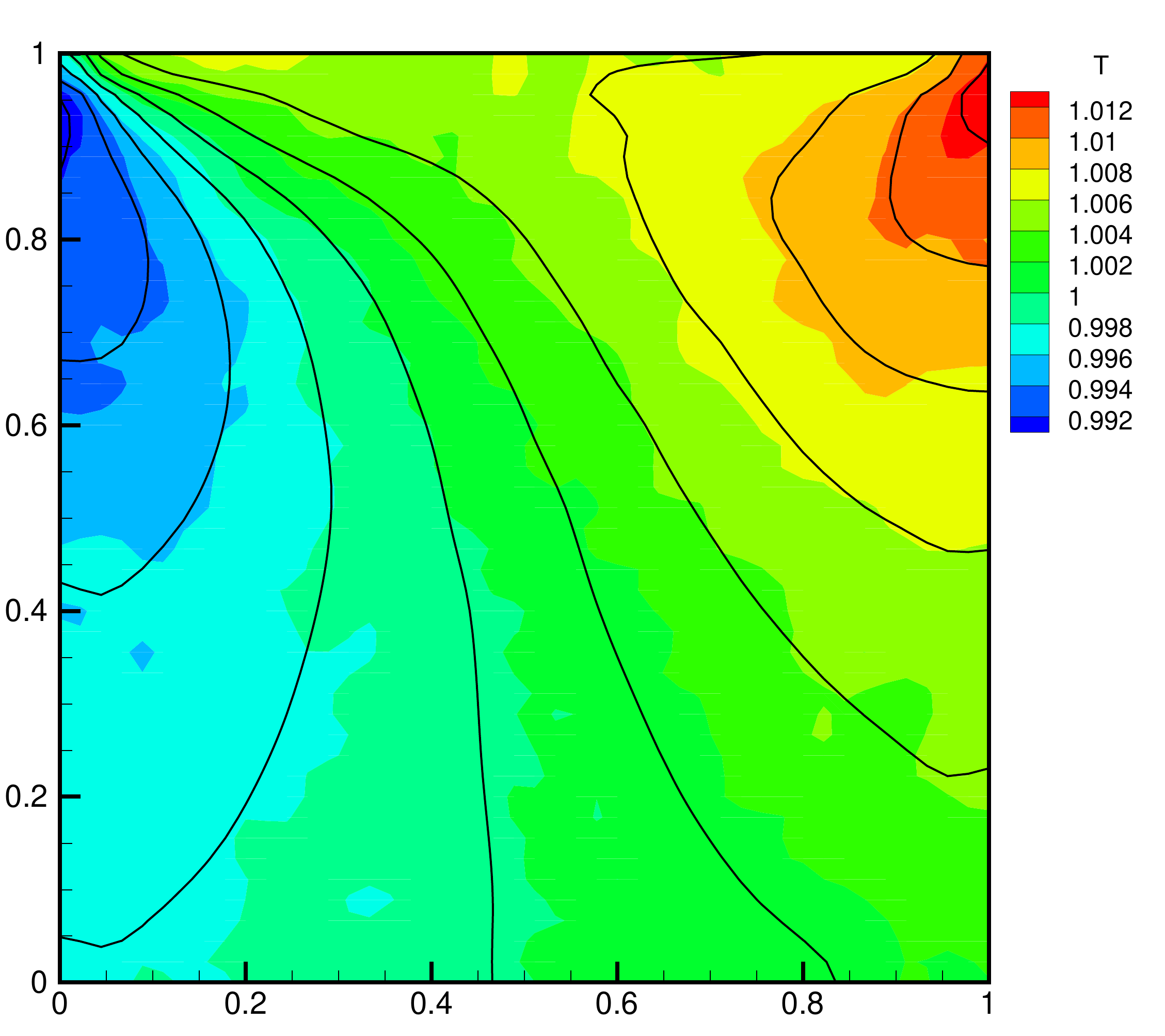}{d}
\caption{(a) Density, (b) x directional velocity, (c) y directional velocity, and (d) temperature contour for the lid-driven cavity flow at $\mathrm{Kn}=1$. The UGKWP method solution is shown in flood, and the UGKS solution is shown in contour line.}
\label{cavity11}
\end{figure}

\begin{figure}
\centering
\includegraphics[width=0.48\textwidth]{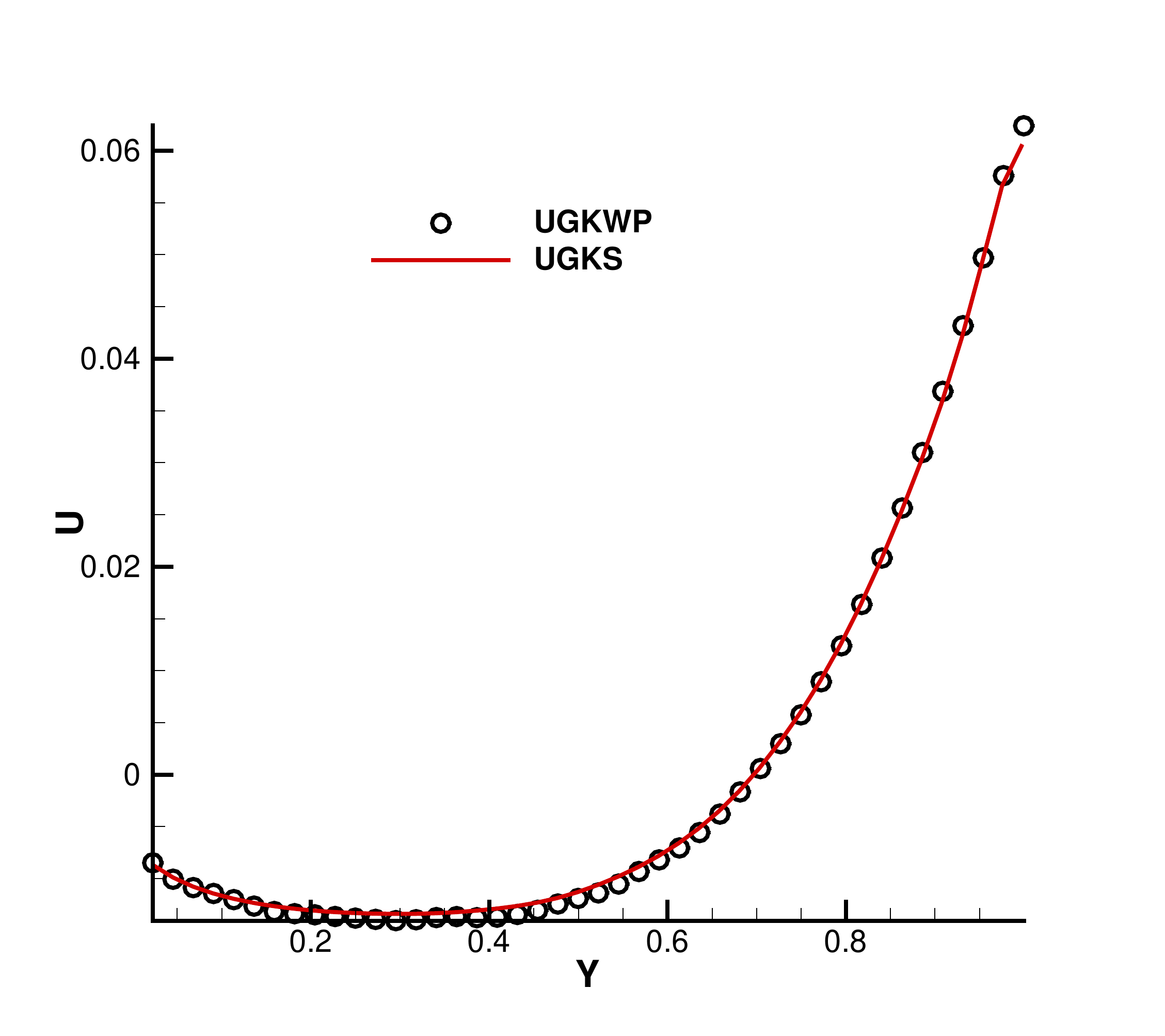}
\includegraphics[width=0.48\textwidth]{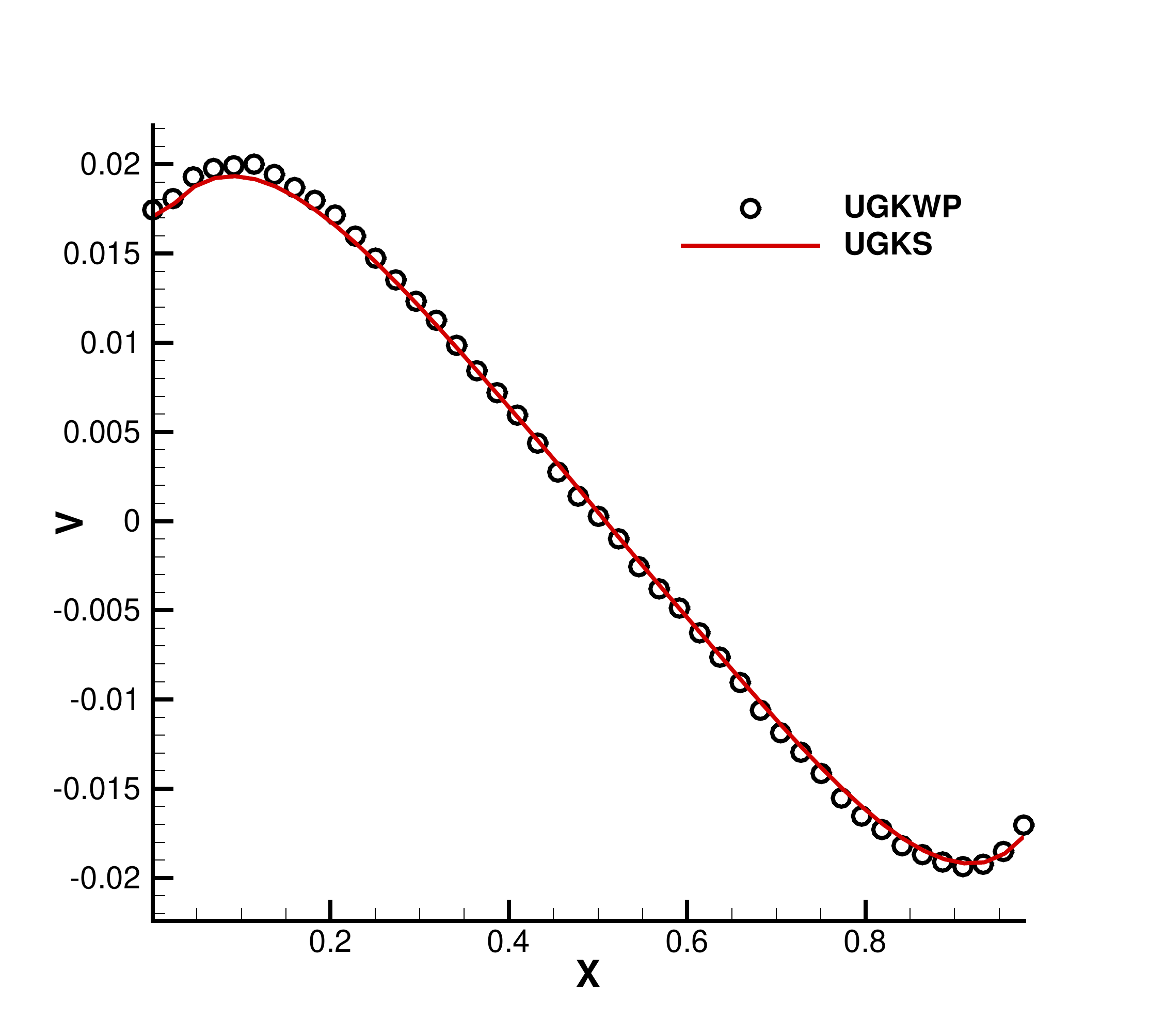}
\caption{Left figure shows x directional velocity along $x=0.5$, and right figure shows y directional velocity along $y=0.5$ for lid-driven cavity flow at $\mathrm{Kn}=1$. Solution of UGKWP method is shown in symbol and UGKS solution is shown in line.}
\label{cavity12}
\end{figure}

\begin{figure}
\centering
\includegraphics[width=0.48\textwidth]{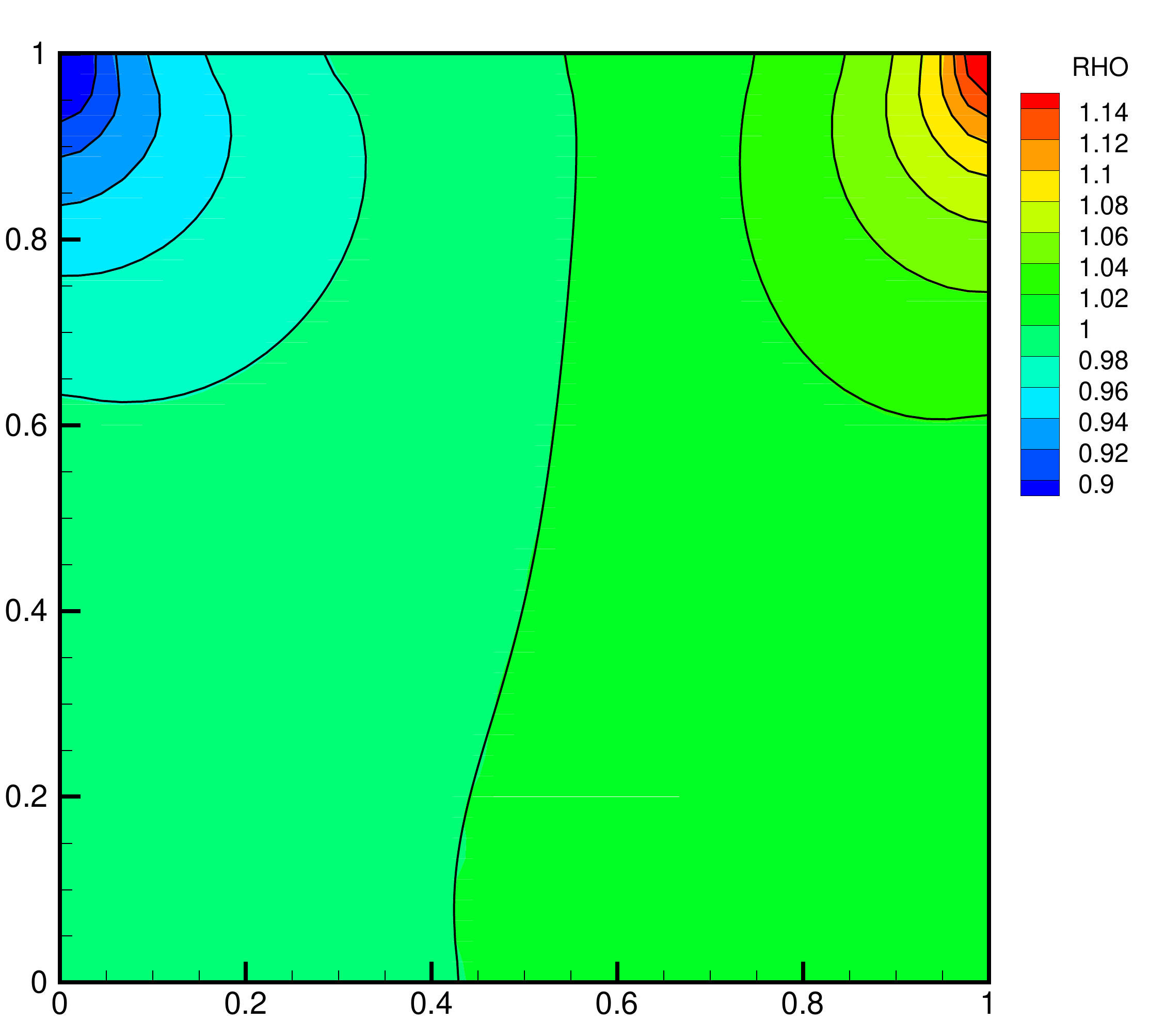}{a}
\includegraphics[width=0.48\textwidth]{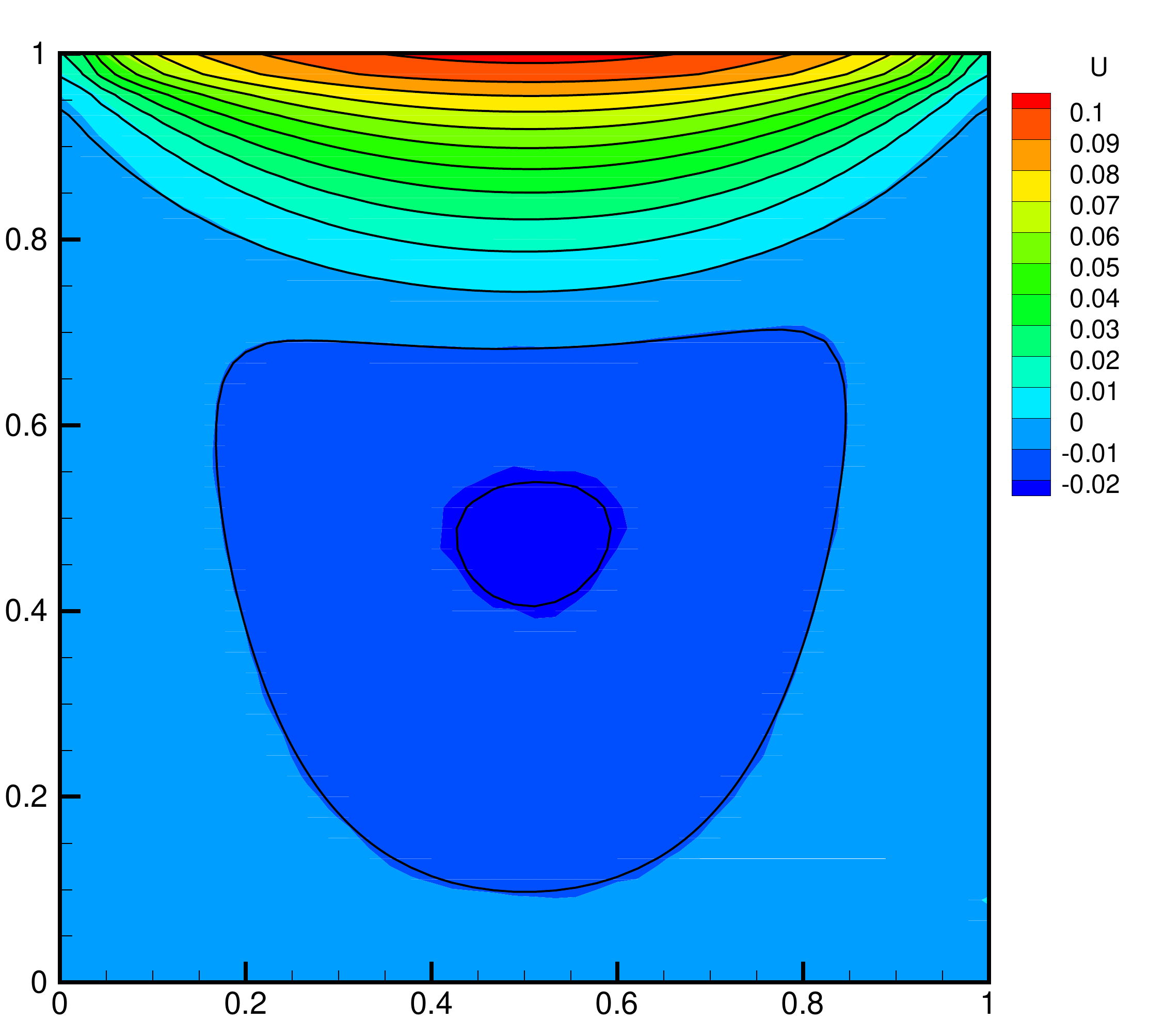}{b}
\includegraphics[width=0.48\textwidth]{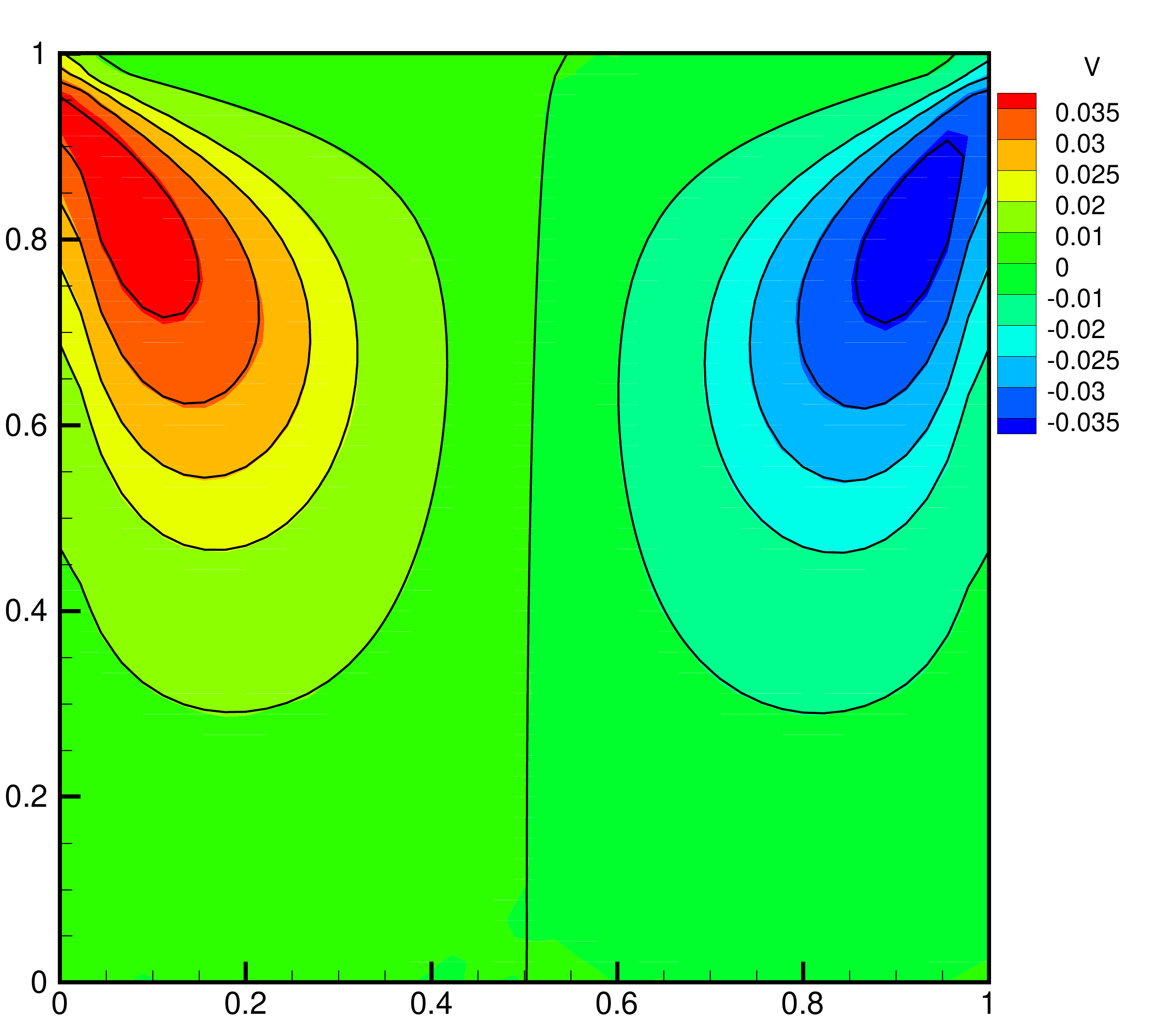}{c}
\includegraphics[width=0.48\textwidth]{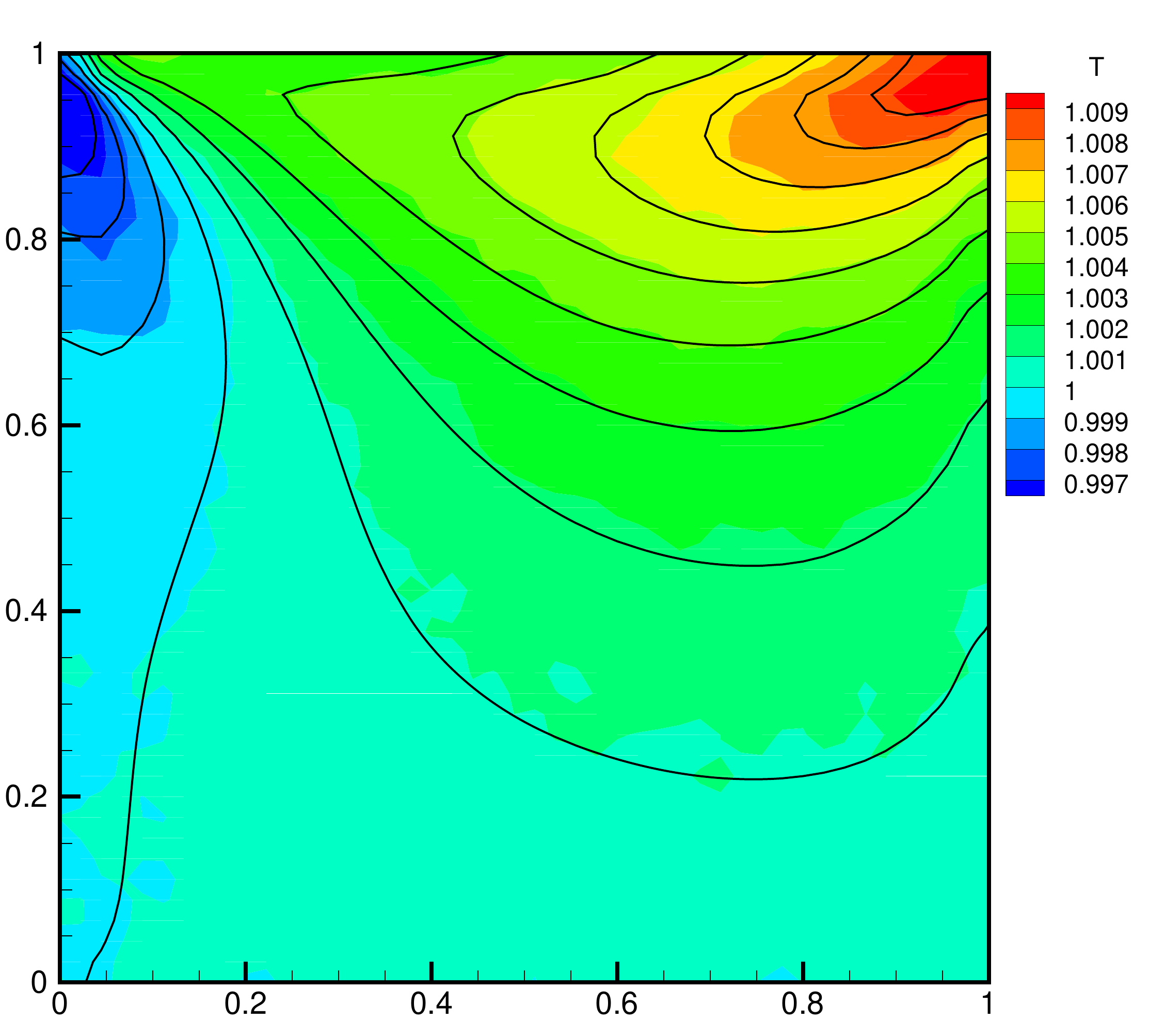}{d}
\caption{(a) Density, (b) x directional velocity, (c) y directional velocity, and (d) temperature contour for the lid-driven cavity flow at $\mathrm{Kn}=0.075$. The UGKWP method solution is shown in flood, and the UGKS solution is shown in contour line.}
\label{cavity21}
\end{figure}

\begin{figure}
\centering
\includegraphics[width=0.48\textwidth]{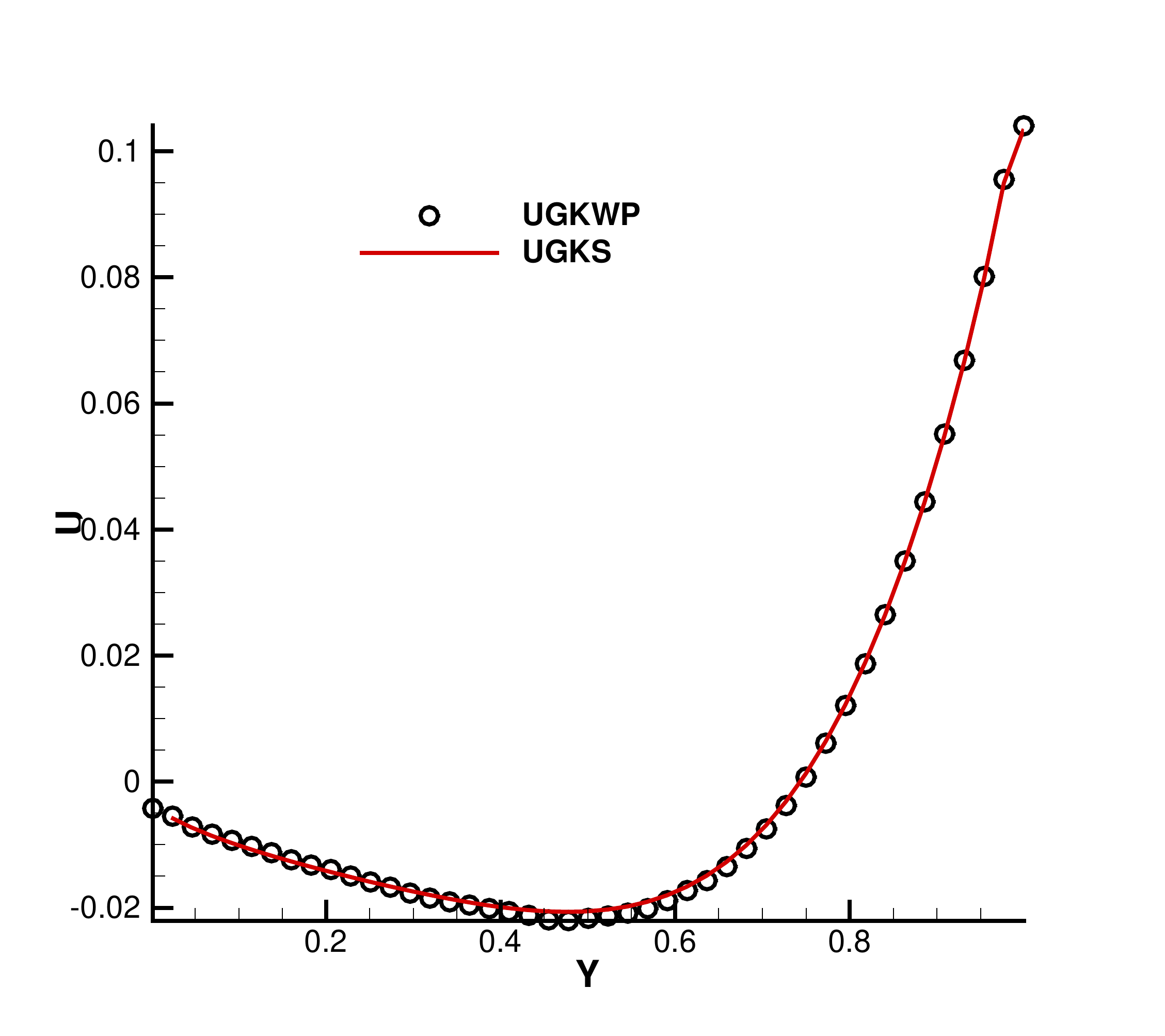}
\includegraphics[width=0.48\textwidth]{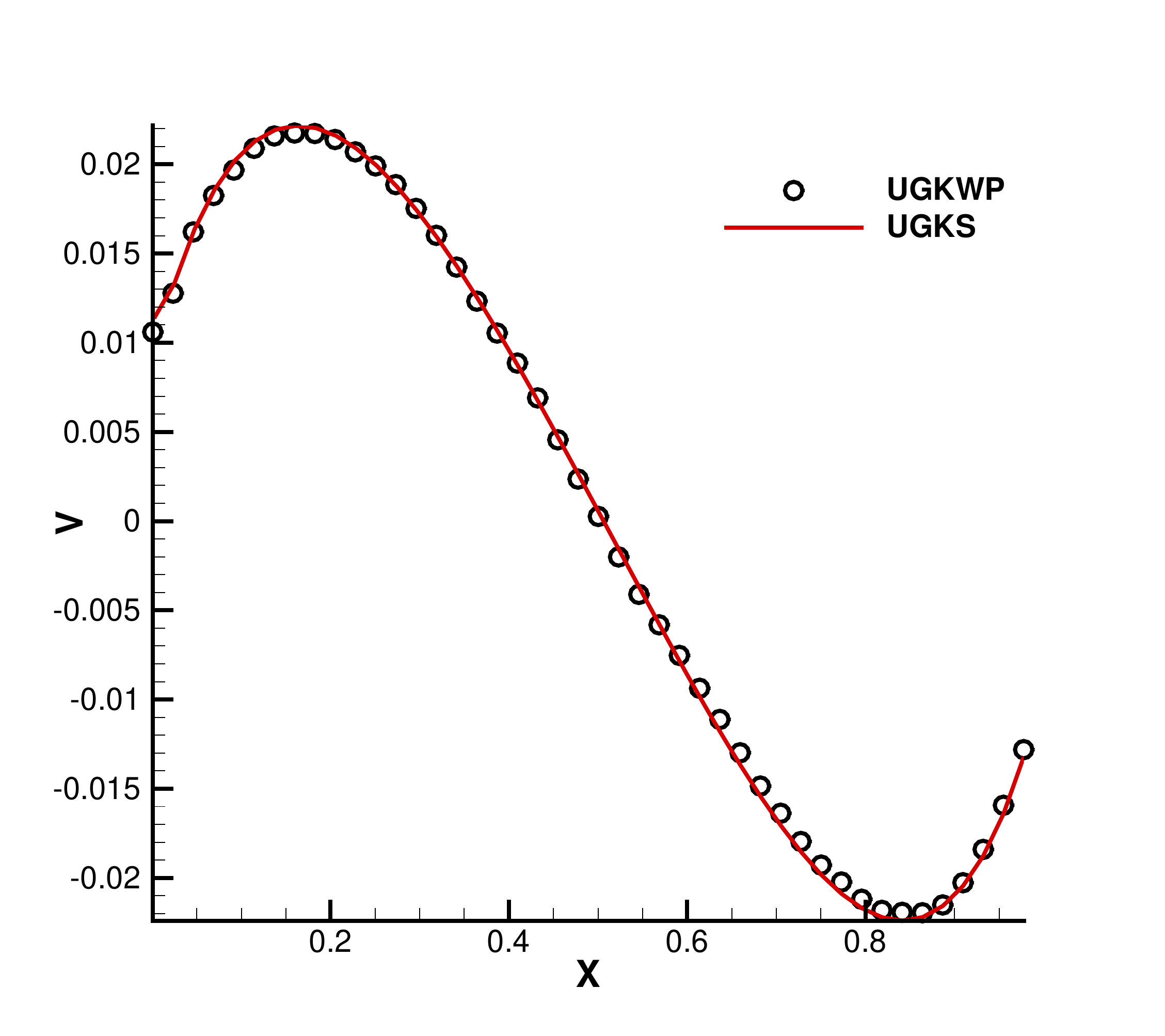}
\caption{Left figure shows x directional velocity along $x=0.5$, and right figure shows y directional velocity along $y=0.5$ for lid-driven cavity flow at $\mathrm{Kn}=0.075$. Solution of UGKWP method is shown in symbol and UGKS solution is shown in line.}
\label{cavity22}
\end{figure}

For $\mathrm{Kn}=1.0$, $5000$ number of simulation particles per cell is used for UGKWP method and UGKS use $50\times50$ discrete velocity points. The computational time for UGKS is 14.4 hours, and for UGKWP method is 20 hours
(including 1800 averaging steps). The memory cost for UGKS is 500 MB and for UGKWP method is 2.5 GB. The computational time of UGKWP method is 1.38 times slower than explicit UGKS, and the memory cost is 5 times larger than UGKS. For $\mathrm{Kn}=0.075$, the same numerical set up is used with a similar computational cost as the case of $\mathrm{Kn}=1.0$.
The solution of the UGKWP method is shown in Fig. \ref{cavity11}-\ref{cavity22}, compared to the UGKS solution. Fig. \ref{cavity11} and \ref{cavity21} show the density, velocity, and temperature contours of UGKWP method (flood) and UGKS (lines). For both $\mathrm{Kn}=1.0$ and $\mathrm{Kn}=0.075$, the density and velocity contours agree well between two methods, while the temperature shows relative large statistical noise for UGKWP method.
Fig. \ref{cavity12} and \ref{cavity21} show the velocity profile along the vertical and horizontal lines, and good agreement can be observed between the UGKWP method and UGKS.

\begin{figure}
\centering
\includegraphics[width=0.48\textwidth]{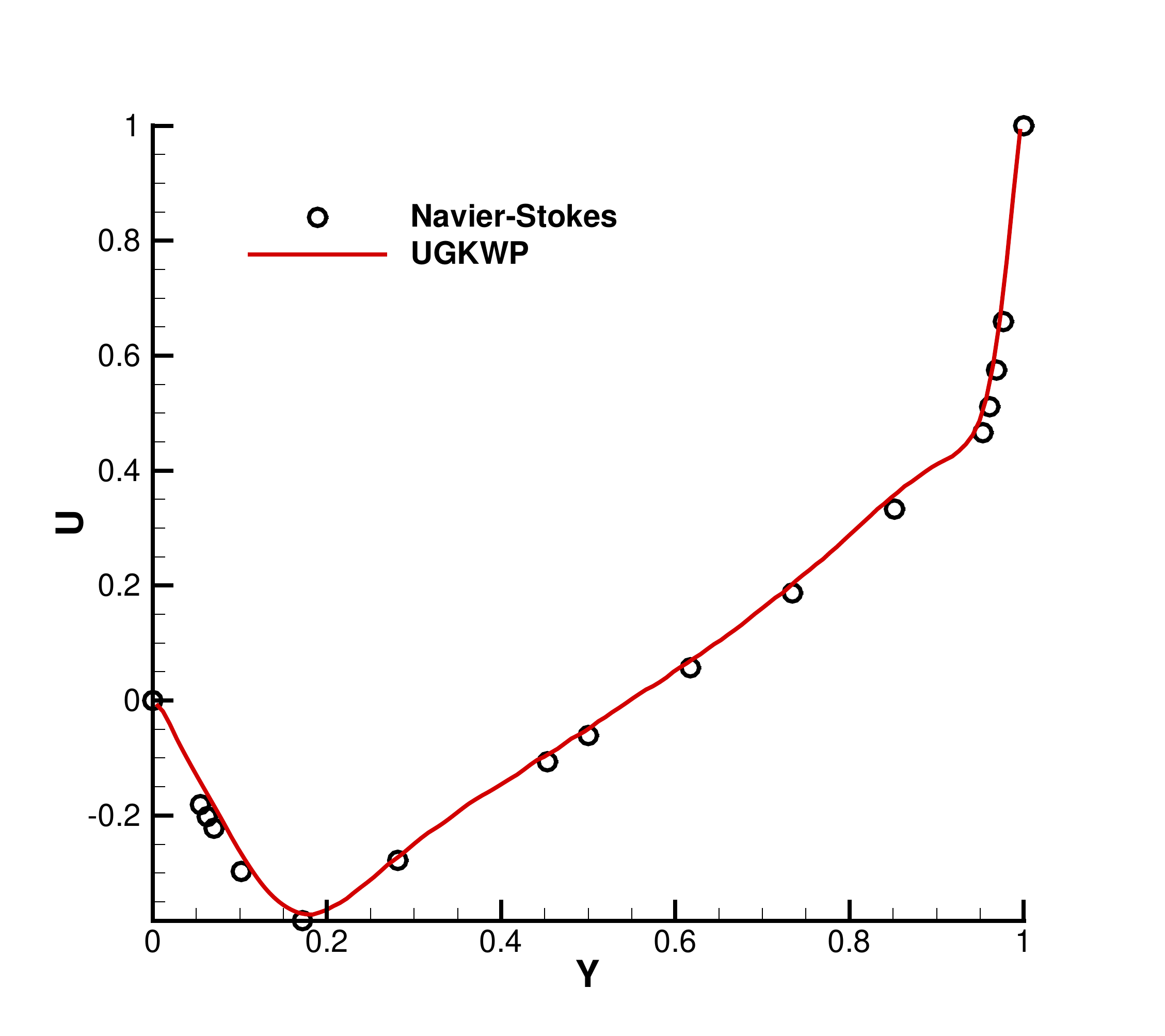}
\includegraphics[width=0.48\textwidth]{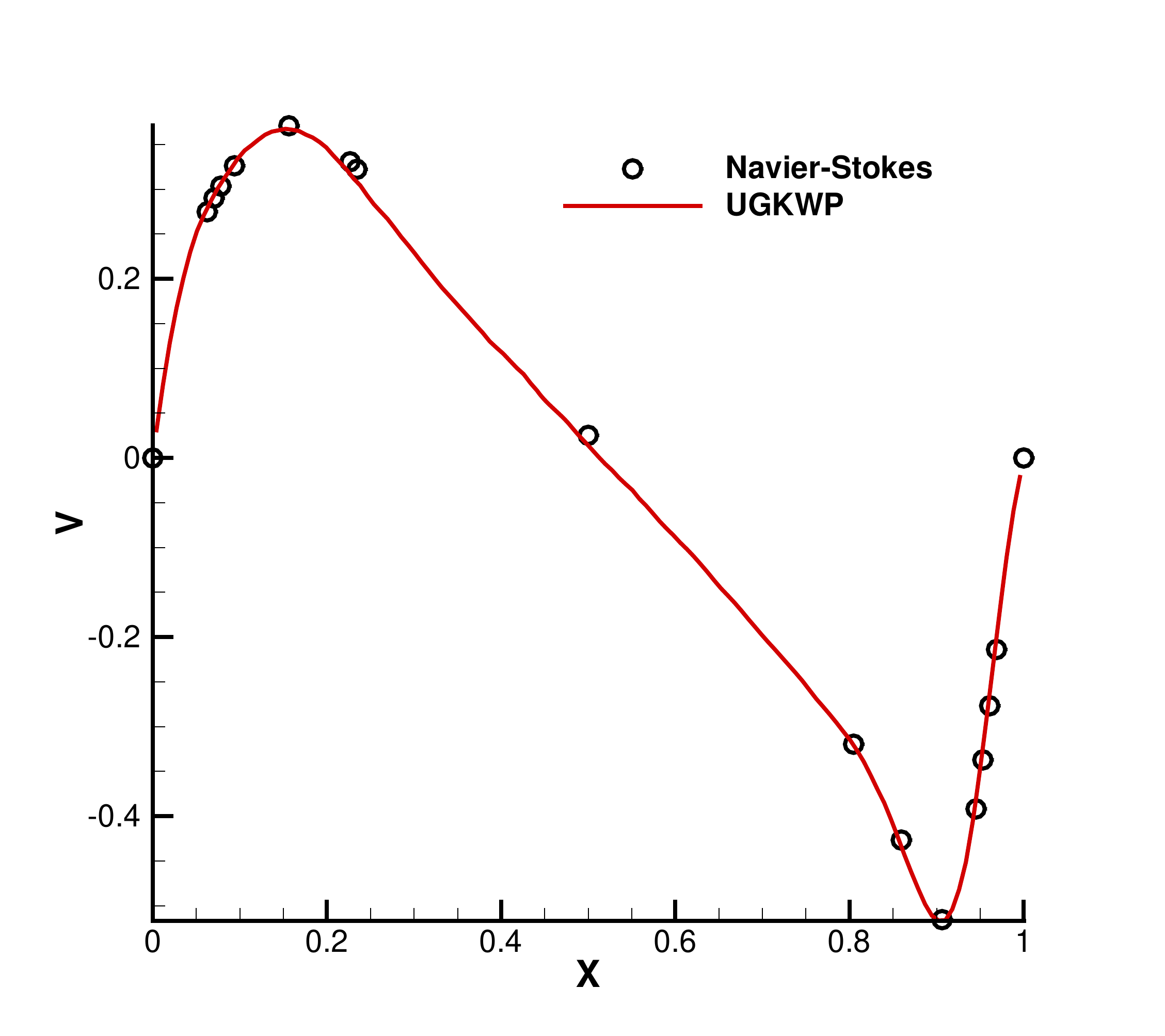}
\caption{Left figure shows x directional velocity along $x=0.5$, and right figure shows y directional velocity along $y=0.5$ for lid-driven cavity flow at $\mathrm{Re}=1000$. Solution of UGKWP method is shown in symbol and Navier-Stokes solution is shown in line.}
\label{cavity32}
\end{figure}

Next, we calculate the cavity flow at $\mathrm{Kn}=1.42 \times 10^{-4}$, i.e., $\mathrm{Re}=1000$.
In this calculation, the number of particle used for UGKWP method is 100 particles per cell,
and UGKS use $28\times28$ discrete velocity points.
The computational time for UGKS is 11 hours, and for UGKWP method is 38.7 minutes.
The memory cost for UGKS is 13 GB and for UGKWP method is 500 MB.
The computational time of UGKWP method is 17 times faster than explicit UGKS, and the memory cost is 26 times less than UGKS. The velocity profiles of UGKWP method along vertical and horizontal lines agrees well with the Navier-Stokes solution as shown in Fig. \ref{cavity32}.

For the low speed flow calculation, in the rarefied flow regime the UGKWP method is more expensive than the UGKS, and in the continuum flow, the UGKWP method is more efficient than the explicit UGKS, with a numerical cost close to GKS for the Navier-Stokes solution.
UGKS can use other acceleration techniques, such as implicit and multigrid \cite{zhu2016implicit,zhu2017implicit}, which improve the efficiency of UGKS by hundreds of times.
But, the memory in UGKS can be hardly reduced due to the discretization of particle velocity space in all flow regimes.

\subsection{Boundary layer}
\begin{figure}
\centering
\includegraphics[width=0.48\textwidth]{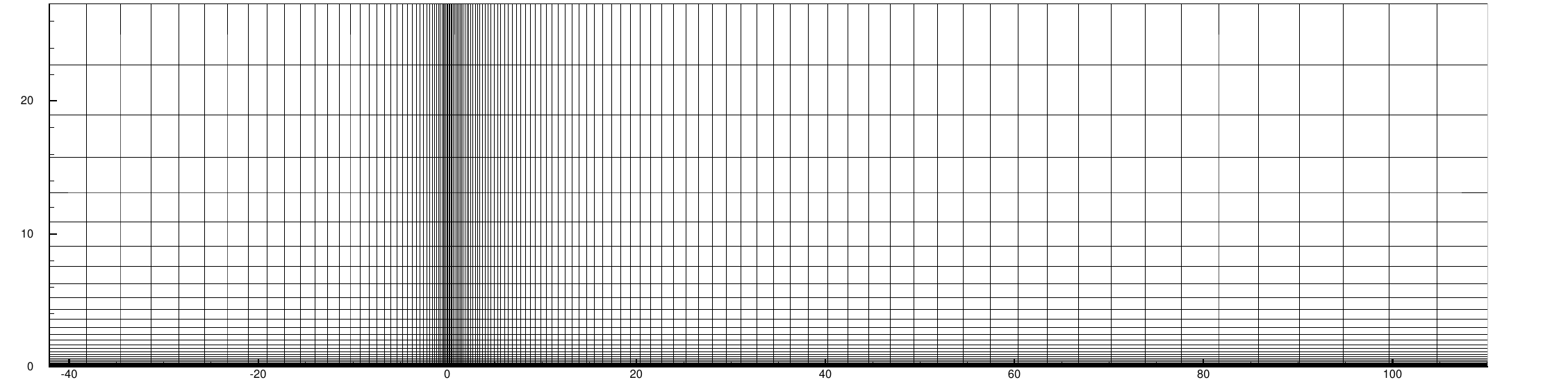}{a}
\includegraphics[width=0.48\textwidth]{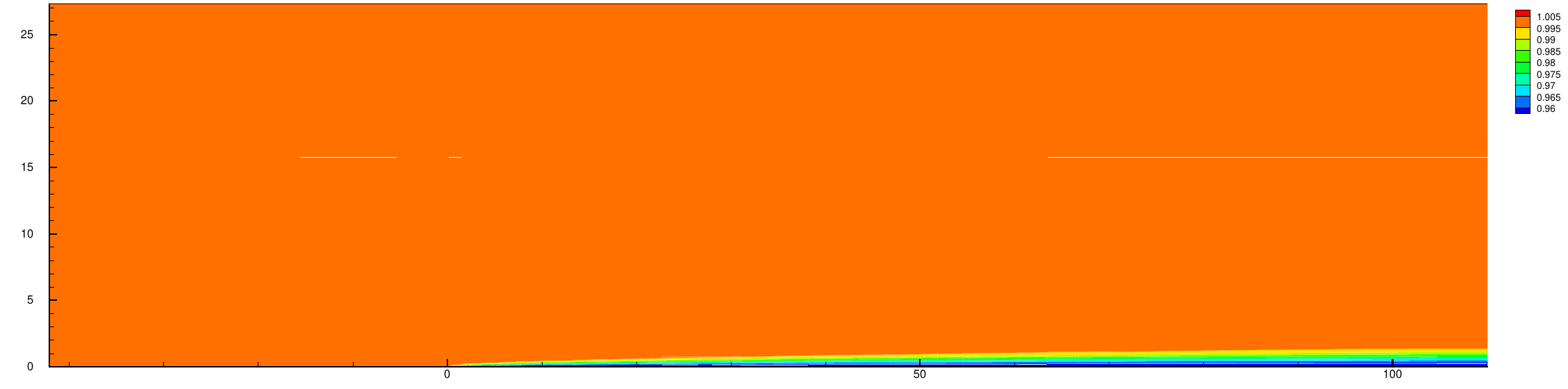}{b}
\includegraphics[width=0.48\textwidth]{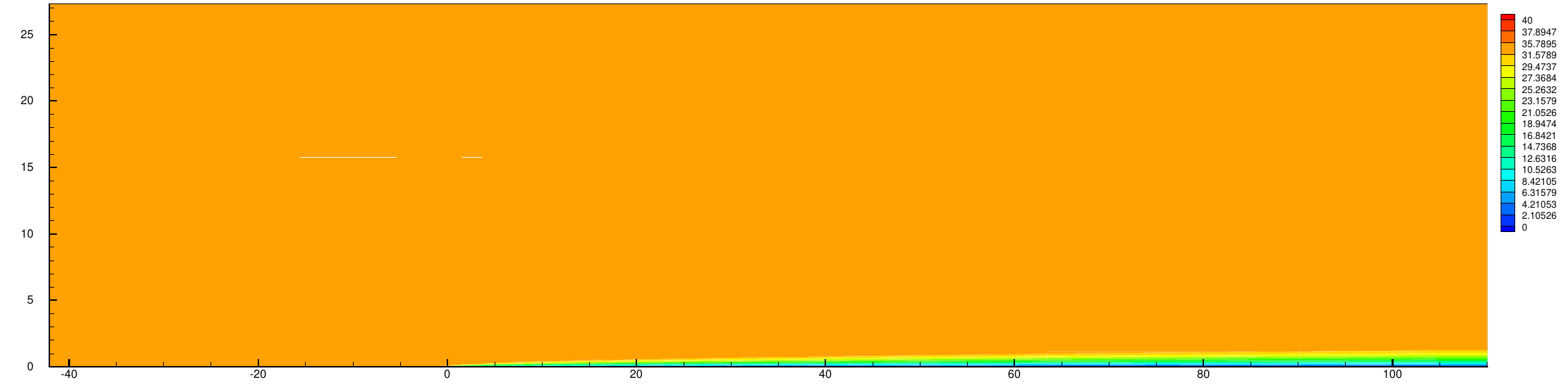}{c}
\includegraphics[width=0.48\textwidth]{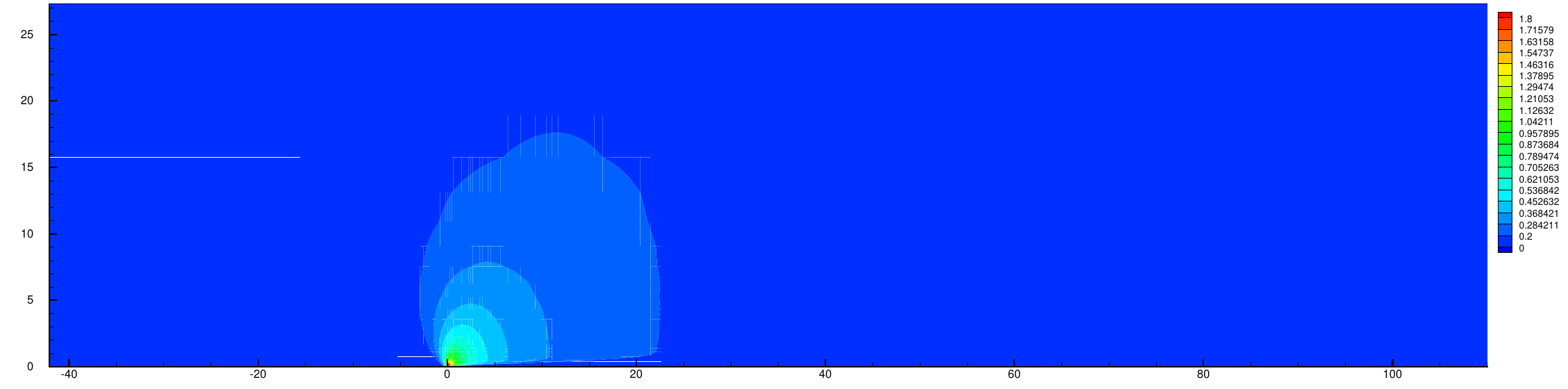}{d}
\caption{Laminar boundary layer computation using UGKWP method at $\mathrm{M}=0.3$ and $\mathrm{Re}=10^{5}$.
(a) mesh distribution; (b) density contours; (c) U velocity contours; (d) V velocity contours.}
\label{layer1}
\end{figure}

\begin{figure}
\centering
\includegraphics[width=0.48\textwidth]{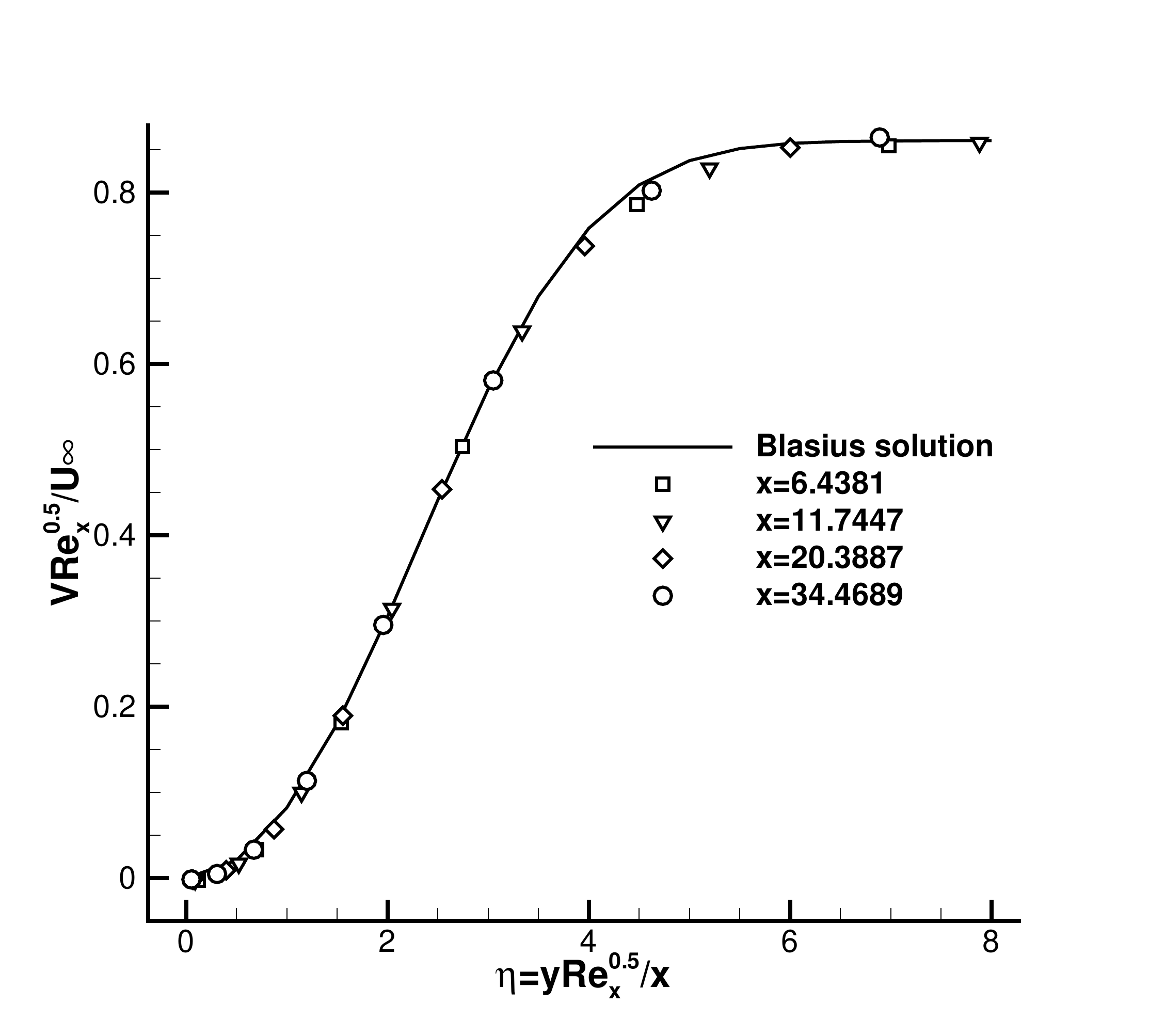}{a}
\includegraphics[width=0.48\textwidth]{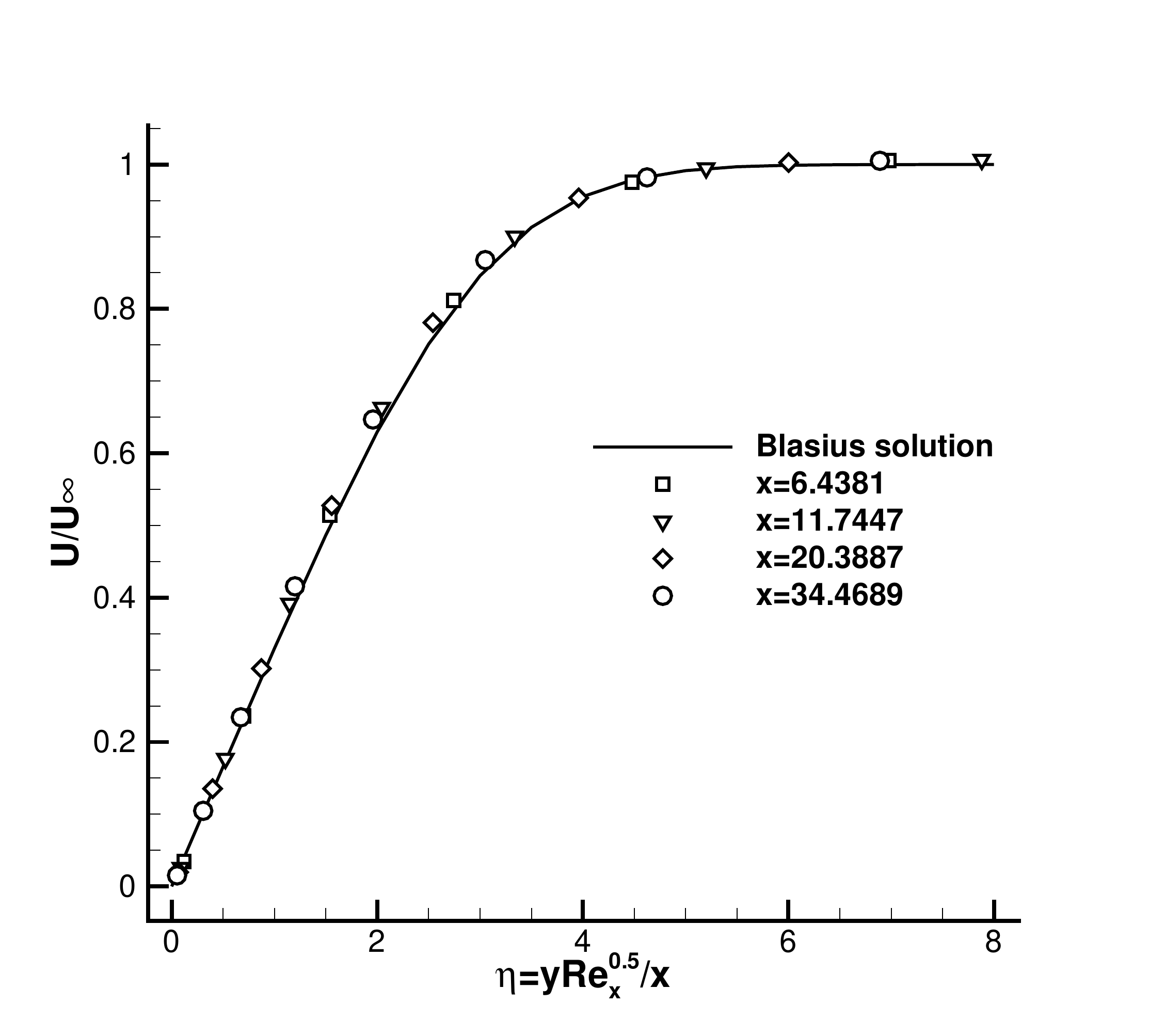}{b}
\caption{Velocity profile of UGKWP method comparing to Navier-Stokes reference solution. (a) U-velocity distribution at different locations; (b) V-velocity distribution at different locations. Symbols: solution of UGKWP method, lines: Blasius solution.}
\label{layer2}
\end{figure}

It is challenging for a particle method to calculate the Navier-Stokes solution under a cell size much larger than the mean free path and time step much larger than the collision time.
To show the ability of the UGKWP method in capturing the Navier-Stokes solution under such condition, we calculate the Navier-Stokes boundary layer.
A gas flow with density $\rho_0=1.0$ and temperature $T_0=5.56\times10^{-2}$ passes over a flat plate at speed $U_0=0.1$. The Reynolds number is set to be $\mathrm{Re}=10^{5}$, and the viscosity is fixed at $\mu=1.05\times10^{-4}$.
The computational domain is $[-44.16,112.75]\times[0,29.8]$, a rectangular mesh with $120\times30$ nonuniform grid points is used as shown in Fig. \ref{layer1}(a).
The CFL condition is chosen as 0.95. At steady state, the velocity profile of UGKWP method agrees well with the Blasius solution as shown in Fig. \ref{layer2}.
For this calculation, the simulation particle mass is $6.32\times10^{-18}$.
Since the time step is much larger than the local collision time, only about 1-2 particles will be stored in each cell, and the computational time is less than 2 minutes.

\section{Conclusion}\label{conclusion}

In this paper, a multiscale multi-efficiency preserving unified gas-kinetic wave-particle method is proposed under the UGKS framework.
The UGKWP method is highly efficient in the simulation of multiscale gas flows from hypersonic to low speed microflows.
In the UGKWP method, both probability density distribution and simulation particles are used to describe the gas particles,
and the simulation particles are sampled only for capturing local non-equilibrium caused by the particle free transport.
The evolution of microscopic quantities is coupled with the evolution of macroscopic quantities in the mesh size scale,
where the flow physics has been directly modeled.

The multiscale modeling or the multiple physical-numerical modeling requires the inclusion of numerical cell size and time step scale into the construction of numerical models instead of a direct discretization of partial differential equations.
UGKS and UGKWP method model gas evolution on the scales of cell size and time step.
According to the ratio between the time step and the particle collision time,
the schemes capture flow dynamics in all flow regimes efficiently.
The UGKWP method is a multi-efficiency preserving scheme, which means the computational cost of the scheme is on the scale of particle methods in the rarefied regime,
and comparable to hydrodynamic solvers in the continuum regime.
Specially, the UGKWP method converges to the gas-kinetic scheme in the continuum regime and does not suffer from stochastic noise.
Due to the implementation of simulation particles in UGKWP method, for hypersonic flows the computational cost has been reduced by hundreds of times in comparison with the UGKS with a direct discretization of particle velocity space.
The current method is based on the BGK-type model for the particle collision and more realistic collisional model can be included in UGKWP method as well.

The methodology of UGKWP method is important for the theoretical fluid dynamics study as well.
On the mesh size scale, even in the near continuum flow regime the use of particles can capture the non-equilibrium effect through the particle streaming or fluid element penetration when a cell size is much larger than the turbulent eddies, which cannot be described by the hydrodynamic equations with averaged flow variables only. The fundamental difficulties associated with the Navier-Stokes equations for the description of separation flow and turbulence may come from the continuum mechanics assumption for continuous connection of fluid elements without breakdown.
The direct modeling equations in UGKWP method release such a constraint and may help in the study of non-equilibrium turbulent flow.
 The non-equilibrium transports in other system, such as chemical reaction, plasma, multiphase, and
radiation, can be solved efficiently by UGKWP method as well.

\section*{Acknowledgments}
The current research is supported by Hong Kong research grant council (16207715, 16206617)
and National Science Foundation of China (11772281, 91530319).

\section*{References}
\bibliography{waveparticle}

\begin{thebibliography}{10}
\expandafter\ifx\csname url\endcsname\relax
  \def\url#1{\texttt{#1}}\fi
\expandafter\ifx\csname urlprefix\endcsname\relax\def\urlprefix{URL }\fi
\expandafter\ifx\csname href\endcsname\relax
  \def\href#1#2{#2} \def\path#1{#1}\fi

\bibitem{bird1963approach}
G.~Bird, Approach to translational equilibrium in a rigid sphere gas, The
  Physics of Fluids 6~(10) (1963) 1518--1519.

\bibitem{bird1994molecular}
G.~Bird, Molecular gas dynamics and the direct simulation {Monte} {Carlo} of
  gas flows, Clarendon, Oxford 508 (1994) 128.

\bibitem{shen2001information}
C.~Shen, J.~Jiang, J.~Fan, Information preservation method for the case of
  temperature variation, in: AIP Conference Proceedings, Vol. 585, AIP, 2001,
  pp. 185--192.

\bibitem{sun2002direct}
Q.~Sun, I.~D. Boyd, A direct simulation method for subsonic, microscale gas
  flows, Journal of Computational Physics 179~(2) (2002) 400--425.

\bibitem{baker2005variance}
L.~L. Baker, N.~G. Hadjiconstantinou, Variance reduction for {Monte} {Carlo}
  solutions of the {Boltzmann} equation, Physics of Fluids 17~(5) (2005)
  051703.

\bibitem{homolle2007low}
T.~M. Homolle, N.~G. Hadjiconstantinou, A low-variance deviational simulation
  {Monte} {Carlo} for the {Boltzmann} equation, Journal of Computational
  Physics 226~(2) (2007) 2341--2358.

\bibitem{pareschi2000asymptotic}
L.~Pareschi, G.~Russo, Asymptotic preserving {Monte} {Carlo} methods for the
  {Boltzmann} equation, Transport Theory and Statistical Physics 29~(3-5)
  (2000) 415--430.

\bibitem{ren2014asymptotic}
W.~Ren, H.~Liu, S.~Jin, An asymptotic-preserving {Monte} {Carlo} method for the
  {Boltzmann} equation, Journal of Computational Physics 276 (2014) 380--404.

\bibitem{degond2011moment}
P.~Degond, G.~Dimarco, L.~Pareschi, The moment-guided {Monte} {Carlo} method,
  International Journal for Numerical Methods in Fluids 67~(2) (2011) 189--213.

\bibitem{burt2008low}
J.~M. Burt, I.~D. Boyd, A low diffusion particle method for simulating
  compressible inviscid flows, Journal of Computational Physics 227~(9) (2008)
  4653--4670.

\bibitem{schwartzentruber2006hybrid}
T.~E. Schwartzentruber, I.~D. Boyd, A hybrid particle-continuum method applied
  to shock waves, Journal of Computational Physics 215~(2) (2006) 402--416.

\bibitem{macrossan2001particle}
M.~N. Macrossan, A particle simulation method for the {BGK} equation, in: AIP
  Conference Proceedings, Vol. 585, AIP, 2001, pp. 426--433.

\bibitem{fei2018particle}
F.~Fei, J.~Zhang, J.~Li, Z.~Liu, A unified stochastic particle
  {Bhatnagar}--{Gross}--{Krook} method for multiscale gas flows, arXiv preprint
  arXiv:1808.03801.

\bibitem{jenny2010solution}
P.~Jenny, M.~Torrilhon, S.~Heinz, A solution algorithm for the fluid dynamic
  equations based on a stochastic model for molecular motion, Journal of
  computational physics 229~(4) (2010) 1077--1098.

\bibitem{chu1965}
C.~Chu, Kinetic theoretic description of the formation of a shock wave, Phys.
  Fluids 8~(1) (1965) 12--22.

\bibitem{JCHuang1995}
J.~Yang, J.~Huang, Rarefied flow computations using nonlinear model {Boltzmann}
  equations, J. Comput. Phys. 120~(2) (1995) 323--339.

\bibitem{Mieussens2000}
L.~Mieussens, Discrete-velocity models and numerical schemes for the
  {Boltzmann-BGK} equation in plane and axisymmetric geometries, J. Comput.
  Phys. 162~(2) (2000) 429--466.

\bibitem{tcheremissine2005direct}
F.~Tcheremissine, Direct numerical solution of the boltzmann equation, Tech.
  rep., DTIC Document (2005).

\bibitem{Kolobov2007}
V.~Kolobov, R.~Arslanbekov, V.~Aristov, A.~Frolova, S.~Zabelok, Unified solver
  for rarefied and continuum flows with adaptive mesh and algorithm refinement,
  J. Comput. Phys. 223~(2) (2007) 589--608.

\bibitem{LiZhiHui2009}
Z.~Li, H.~Zhang, Gas-kinetic numerical studies of three-dimensional complex
  flows on spacecraft re-entry, J. Comput. Phys. 228~(4) (2009) 1116--1138.

\bibitem{ugks2010}
K.~Xu, J.-C. Huang, A unified gas-kinetic scheme for continuum and rarefied
  flows, Journal of Computational Physics 229~(20) (2010) 7747--7764.

\bibitem{wu2015fast}
L.~Wu, J.~Zhang, J.~M. Reese, Y.~Zhang, A fast spectral method for the
  {Boltzmann} equation for monatomic gas mixtures, Journal of Computational
  Physics 298 (2015) 602--621.

\bibitem{aristov2012direct}
V.~Aristov, Direct methods for solving the {Boltzmann} equation and study of
  nonequilibrium flows.

\bibitem{huang2013}
J.-C. Huang, K.~Xu, P.~Yu, A unified gas-kinetic scheme for continuum and
  rarefied flows {III}: Microflow simulations, Communications in Computational
  Physics 14~(5) (2013) 1147--1173.

\bibitem{wu2014solving}
L.~Wu, J.~M. Reese, Y.~Zhang, Solving the {Boltzmann} equation
  deterministically by the fast spectral method: application to gas microflows,
  Journal of Fluid Mechanics 746 (2014) 53--84.

\bibitem{jin1999efficient}
S.~Jin, Efficient asymptotic-preserving {(AP)} schemes for some multiscale
  kinetic equations, SIAM Journal on Scientific Computing 21~(2) (1999)
  441--454.

\bibitem{chen2017unified}
S.~Chen, C.~Zhang, L.~Zhu, Z.~Guo, A unified implicit scheme for kinetic model
  equations. part {I}. memory reduction technique, Science bulletin 62~(2)
  (2017) 119--129.

\bibitem{degond2010multiscale}
P.~Degond, G.~Dimarco, L.~Mieussens, A multiscale kinetic--fluid solver with
  dynamic localization of kinetic effects, Journal of Computational Physics
  229~(13) (2010) 4907--4933.

\bibitem{xu-book}
K.~Xu, Direct Modeling for Computational Fluid Dynamics: Construction and
  Application of Unified Gas-kinetic Scheme, World Scientic, 2015.

\bibitem{huang2012}
J.-C. Huang, K.~Xu, P.~Yu, A unified gas-kinetic scheme for continuum and
  rarefied flows {II}: multi-dimensional cases, Communications in Computational
  Physics 12~(3) (2012) 662--690.

\bibitem{liu2016}
C.~Liu, K.~Xu, Q.~Sun, Q.~Cai, A unified gas-kinetic scheme for continuum and
  rarefied flows {IV}: Full {Boltzmann} and model equations, Journal of
  Computational Physics 314 (2016) 305--340.

\bibitem{zhu2016implicit}
Y.~Zhu, C.~Zhong, K.~Xu, Implicit unified gas-kinetic scheme for steady state
  solutions in all flow regimes, Journal of Computational Physics 315 (2016)
  16--38.

\bibitem{zhu2017implicit}
Y.~Zhu, C.~Zhong, K.~Xu, Unified gas-kinetic scheme with multigrid convergence
  for rarefied flow study, Physics of Fluids 29~(9) (2017) 096102.

\bibitem{zhu2018implicit}
Y.~Zhu, C.~Zhong, K.~Xu, An implicit unified gas-kinetic scheme for unsteady
  flow in all knudsen regimes, arXiv preprint arXiv:1801.02022.

\bibitem{jiang}
D.~Jiang, Study of the gas-kinetic scheme based on the analytic solution of
  model equations, Ph.D. Thesis, China Aerodynamics Research and Development
  Center.

\bibitem{sun2015asymptotic}
W.~Sun, S.~Jiang, K.~Xu, S.~Li, An asymptotic preserving unified gas kinetic
  scheme for frequency-dependent radiative transfer equations, Journal of
  Computational Physics 302 (2015) 222--238.

\bibitem{sun2017}
W.~Sun, S.~Jiang, K.~Xu, A multidimensional unified gas-kinetic scheme for
  radiative transfer equations on unstructured mesh, Journal of Computational
  Physics 351 (2017) 455--472.

\bibitem{sun2018}
W.~Sun, S.~Jiang, K.~Xu, An asymptotic preserving implicit unified gas kinetic
  scheme for requency-dependent radiative transfer equations., International
  Journal of Numerical Analysis and Modeling 15 (2018) 134--153.

\bibitem{li2018ugkp}
W.~Li, C.~Liu, Y.~Zhu, J.~Zhang, K.~Xu, A unified gas-kinetic particle method
  for multiscale photon transport, arXiv preprint arXiv:1810.05984.

\bibitem{liu2017}
C.~Liu, K.~Xu, A unified gas kinetic scheme for continuum and rarefied flows
  {V}: Multiscale and multi-component plasma transport, Communications in
  Computational Physics 22~(5) (2017) 1175--1223.

\bibitem{liu2018unified}
C.~Liu, Z.~Wang, K.~Xu, A unified gas-kinetic scheme for continuum and rarefied
  flows {VI}: Dilute disperse gas-particle multiphase system, arXiv preprint
  arXiv:1802.04961.

\bibitem{wang2015unified}
R.~Wang, Unified gas-kinetic scheme for the study of non-equilibrium flows,
  Ph.D. thesis, Hong Kong University of Science and Technology (2015).

\bibitem{shakhov1968generalization}
E.~Shakhov, Generalization of the {Krook} kinetic relaxation equation, Fluid
  Dynamics 3~(5) (1968) 95--96.

\bibitem{xu2001gas}
K.~Xu, A gas-kinetic {BGK} scheme for the {Navier}--{Stokes} equations and its
  connection with artificial dissipation and {Godunov} method, Journal of
  Computational Physics 171~(1) (2001) 289--335.

\bibitem{chen2015comparative}
S.~Chen, K.~Xu, A comparative study of an asymptotic preserving scheme and
  unified gas-kinetic scheme in continuum flow limit, Journal of Computational
  Physics 288 (2015) 52--65.

\bibitem{BGK1954}
P.~Bhatnagar, E.~Gross, M.Krook, A model for collision processes in gases.
  \uppercase\expandafter{\romannumeral1}. small amplitude processes in charged
  and neutral one-component systems, Phys. Rev. 94~(3) (1954) 511--525.

\bibitem{chapman1990mathematical}
S.~Chapman, T.~G. Cowling, D.~Burnett, The mathematical theory of non-uniform
  gases: an account of the kinetic theory of viscosity, thermal conduction and
  diffusion in gases, Cambridge university press, 1990.

\end{thebibliography}
\bibliographystyle{elsarticle-num}
\biboptions{numbers,sort&compress}

\clearpage

\clearpage

\begin{figure}
\centering
\includegraphics[width=0.6\textwidth]{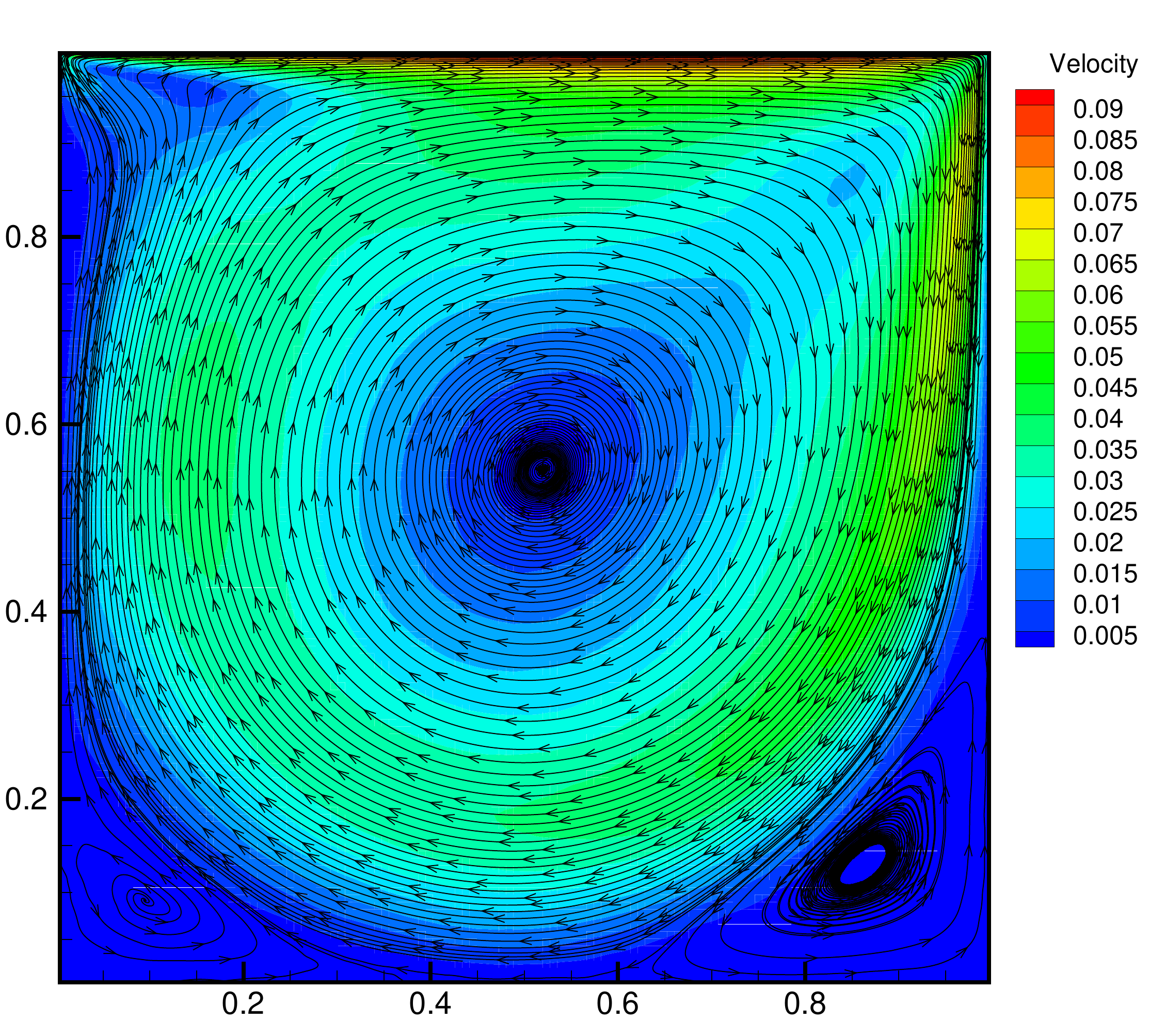}
\caption{Stream line and velocity contour of the lid-driven cavity flow at $\mathrm{Re}=1000$.}
\label{cavity31}
\end{figure}

\end{document}